\documentclass[11pt,oneside,reqno]{article}
\usepackage[width=15.7cm,height=24cm,centering]{geometry}
\usepackage{amsmath,amsthm,amssymb,color,bm}
\usepackage[authoryear]{natbib}
\usepackage{booktabs}

\newcommand{\tr}{\triangle}

\usepackage{tikz}
\RequirePackage{amsthm,amsmath,amsfonts,amssymb,mathrsfs,dsfont,mathtools,thmtools}
\RequirePackage[colorlinks,citecolor=blue,urlcolor=blue,linkcolor=blue]{hyperref}
\usepackage[capitalize]{cleveref}
\RequirePackage{graphicx}

\crefname{equation}{}{}
\crefformat{equation}{\textup{(#2#1#3)}}
\crefrangeformat{equation}{\textup{(#3#1#4)--(#5#2#6)}}
\crefdefaultlabelformat{#2\textup{#1}#3}
\crefname{lemma}{Lemma}{Lemmas}
\crefname{page}{p.}{pp.}
\usepackage[normalem]{ulem}

\usepackage{bbm}
\usepackage{mathrsfs}
\usepackage{graphicx}
\usepackage{tikz}
\usepackage{enumitem}
\usetikzlibrary{arrows, automata}
\usetikzlibrary{calc}
\usetikzlibrary{positioning}

\numberwithin{equation}{section}
\allowdisplaybreaks[4]

\theoremstyle{plain}
\newtheorem{theorem}{Theorem}[section]
\newtheorem{proposition}{Proposition}[section]
\newtheorem{lemma}{Lemma}[section]
\newtheorem{corollary}{Corollary}[section]

\newtheorem{remark}{Remark}[section]

\theoremstyle{definition}

\makeatletter

\newcount\minute
\newcount\hour
\newcount\hourMins
\def\now{%
\minute=\time%
\hour=\time \divide \hour by 60%
\hourMins=\hour \multiply\hourMins by 60%
\advance\minute by -\hourMins%
\zeroPadTwo{\the\hour}:\zeroPadTwo{\the\minute}%
}
\def\zeroPadTwo#1{\ifnum #1<10 0\fi#1}

\renewcommand{\cite}{\citet}

\def\^#1{\ifmmode {\mathaccent"705E #1} \else {\accent94 #1} \fi}
\def\~#1{\ifmmode {\mathaccent"707E #1} \else {\accent"7E #1} \fi}

\def\*#1{#1^\ast}
\edef\-#1{\noexpand\ifmmode {\noexpand\bar{#1}} \noexpand\else \-#1\noexpand\fi}
\def\>#1{\vec{#1}}
\def\.#1{\dot{#1}}

\def\atop{\@@atop}
\def\*#1{\mathscr{#1}}

\renewcommand{\le}{\leqslant}

\renewcommand{\ge}{\geqslant}

\newcommand{\Var}{\mathop{\mathrm{Var}}\nolimits}
\newcommand{\Cov}{\mathop{\mathrm{Cov}}}

\def\mid{\vert}

\def\beqn#1\eeqn{\begin{align}#1\end{align}}
\def\beq#1\eeq{\begin{align*}#1\end{align*}}

\usepackage{graphicx}
\usepackage{latexsym}
\usepackage{amsmath,amsthm,amssymb,amscd}
\usepackage{epsf,amsmath}

\def\E{{\IE}}

\renewcommand\section{\@startsection {section}{1}{\z@}%
{-2.8ex \@plus -.8ex \@minus -.2ex}%
{1.0ex \@plus.2ex}%
{\center\small\sc\mathversion{bold}}}

\renewcommand\subsection{\@startsection{subsection}{2}{0pt}%
{-2.8ex \@plus -.8ex \@minus -.2ex}%
{.8ex \@plus.2ex}%
{\bfseries\mathversion{bold}}}

\def\subsubsection#1{\@startsection{subsubsection}{3}{0pt}%
{\medskipamount}%
{-10pt}%
{\normalsize\itshape}{\kern-2.2ex. #1.}}

\def\blfootnote{\xdef\@thefnmark{}\@footnotetext}

\makeatother

\usepackage{array,needspace}
\renewcommand{\E}{\mathbb E}
\newcommand{\Pp}{\mathbb P}

\newcommand{\R}{\mathbb R}
\newcommand{\ind}{\mathbf1}
\newcommand{\cR}{\mathcal R}

\let\Cov\relax
\let\Var\relax

\let\tr\relax

\DeclareMathOperator{\Cov}{Cov}
\DeclareMathOperator{\Var}{Var}

\DeclareMathOperator{\tr}{tr}

\hypersetup{
 pdftitle={Optimal High-Dimensional Martingale CLTs over Rectangles},
 pdfauthor={Zi-Yao Su and Xiuqi Sun}}

\begin{document}

\title{Optimal High-Dimensional Martingale CLTs over Rectangles}
\author{Zi-Yao Su$^*$, Xiuqi Sun$^\dagger$}
\date{\it The Chinese University of Hong Kong$^*$, Tongji University$^\dagger$}
\maketitle

\noindent{\bf Abstract:}
We prove Gaussian approximation bounds over rectangles for
high-dimensional martingale sums, with lower bounds showing optimality
up to powers of \(\log n\). The results cover bounded increments,
conditional and marginal sub-Gaussian and sub-exponential tails, and
uniform \(q\)th-moment bounds on the largest absolute coordinate for
\(2<q\le3\). The predictable quadratic variation \(Q\) may be random;
the comparison Gaussian has deterministic covariance \(\Sigma=\E Q\).
We assume \(\Sigma_{jj}=1\) and
\(\max_j\sum_k|\Sigma_{jk}|^s\le M\) for fixed \(s>0\) and \(M\ge1\),
allowing \(\Sigma\) to be singular, and quantify the error due to
\(Q-\Sigma\). Conditional tail bounds yield the same dimension
threshold as boundedness, whereas marginal tail bounds give smaller
thresholds. We also establish necessary and sufficient conditions for
uniform convergence in all these classes; under the uniform moment
assumption, the criterion holds for every fixed \(q>2\).
When \(Q=I_d\), the bounded and marginal light-tail rates have no additional factors in
\(\log n\) if \(\log(ed)\ge\{\log(en)\}^2\); the same holds for the
finite-moment rates with \(2<q\le3\) in every dimension, for all
sufficiently large \(n\).

\medskip

\noindent{\bf AMS 2020 subject classification:}
60F05; 60G42; 62E17

\noindent{\bf Keywords and phrases:}
martingale central limit theorem; Gaussian
approximation; high dimension; rectangles; predictable quadratic variation.

\section{Introduction}\label{main:introduction}

Gaussian approximation over rectangles describes the joint accuracy of
coordinatewise probability statements, including the distribution of a
maximum. For martingale sums, its accuracy depends both on the increment
tails and on fluctuations of the sum of conditional covariance matrices.
Even when this sum is deterministic, a bound on unconditional tails can
permit larger errors than the corresponding conditional bound. We study
these distinctions for bounded increments, sub-Gaussian and
sub-exponential tails, and uniform moment bounds on the largest absolute
coordinate. Our aim is to determine the maximal approximation error over
each class and the necessary and sufficient conditions for this error
to vanish.

Let \(X_1,\ldots,X_n\) be \(\R^d\)-valued martingale differences
with respect to a filtration \((\mathcal F_i)_{i=0}^n\). Write
\(\E_i=\E[\,\cdot\mid\mathcal F_i]\), so that
\(\E_{i-1}X_i=0\), and define
\begin{equation}\label{main:eq:normalization}
 S_n=\frac1{\sqrt n}\sum_{i=1}^nX_i,\qquad
 Q=\frac1n\sum_{i=1}^n\E_{i-1}(X_iX_i^\top),\qquad
 \Sigma=\E Q.
\end{equation}
The matrix \(Q\) is the predictable quadratic variation at time \(n\)
of the normalized martingale. Orthogonality of martingale differences
gives \(\Sigma=\Cov(S_n)\), so the Gaussian target
\(G_\Sigma\sim N(0,\Sigma)\) has the same deterministic covariance as
\(S_n\), even when \(Q\) is random. For fixed \(s>0\), \(M\ge1\) and
\(1\le p\le\infty\), our common assumptions are
\begin{equation}\label{main:eq:setting}
 \Sigma_{jj}=1,\qquad
 \max_j\sum_k|\Sigma_{jk}|^s\le M,\qquad
 \big\|\|Q-\Sigma\|_{\max}\big\|_p\le\delta,
\end{equation}
where \(\|A\|_{\max}=\max_{j,k}|A_{jk}|\).
The \(L^\infty\) norm is the essential supremum, and we use
\(1/\infty=0\). The filtration, its initial sigma-field, and the
individual conditional covariance matrices are otherwise arbitrary.
The condition on \(\Sigma\) permits singular matrices. It limits the
number of coordinates strongly correlated with any given coordinate:
for \(t>0\),
\[
 \#\{k:|\Sigma_{jk}|>t\}\le Mt^{-s}.
\]

Let \(\mathcal R_d\) be the axis-parallel rectangles in \(\R^d\),
allowing infinite endpoints. We measure error by
\[
 \rho_{\mathcal R}(U,V)=
 \sup_{A\in\mathcal R_d}|\Pp(U\in A)-\Pp(V\in A)|.
\]
This distance controls coordinate maxima and simultaneous coordinate
constraints. Optimality below refers to the supremum of
\(\rho_{\mathcal R}(S_n,G_\Sigma)\) over a specified class of arrays;
each sum is compared with a centered Gaussian vector having the same
covariance.

The increment assumptions distinguish conditional from marginal tail
control. Bounded increments satisfy \(\|X_i\|_\infty\le B\) almost
surely, with fixed \(B\ge2\). For fixed \(K\ge4\), the conditional
and marginal exponential assumptions are, respectively,
\begin{align}
 &\E_{i-1}\exp\{(|X_{ij}|/K)^\alpha\}\le2
       \quad\text{almost surely for every }i,j,
       \label{main:eq:conditional-tails}\\
 &\E\exp\{(|X_{ij}|/K)^\alpha\}\le2
       \quad\text{for every }i,j,
       \label{main:eq:marginal-tails}
\end{align}
where \(\alpha=2\) gives sub-Gaussian tails and \(\alpha=1\)
gives sub-exponential tails. Marginal control is unconditional.
The finite-moment assumption is
\begin{equation}\label{main:eq:moment-envelope}
 \max_{i\le n}\E\|X_i\|_\infty^q\le B^q,\qquad q>2,
\end{equation}
with fixed \(q,B\). This controls the \(q\)th moment of
\(\max_j|X_{ij}|\), uniformly in the increment index, sample size, and
dimension.

The error bounds have two terms: one depends on the increment
assumption, and the other on the deviation of \(Q\) from \(\Sigma\).
Table~\ref{main:tab:scales} lists the first term. The second is
\[
 \{\delta\log(ed)\}^{1/(2+1/p)}
\]
Add the two terms and take the minimum with one. On each stated
class, this gives a lower bound for the maximal rectangle error up to
a fixed constant; the upper bound differs by a power of \(\log n\).
Theorem~\ref{main:joint} specifies this power, the dependence of the
constants, and the range of \(n\).

\Needspace{14\baselineskip}
\begin{center}
\begin{minipage}{\textwidth}
\refstepcounter{table}\label{main:tab:scales}
\centering
\small
\textbf{Table \thetable.} Error terms and necessary and sufficient dimension conditions.
\par\medskip
\setlength{\tabcolsep}{3pt}
\renewcommand{\arraystretch}{1.32}
\begin{tabular}{@{}>{\raggedright\arraybackslash}p{.34\textwidth}
 >{\centering\arraybackslash}p{.37\textwidth}
 >{\centering\arraybackslash}p{.25\textwidth}@{}}
\toprule
Increment assumption & Term due to the increments & Dimension condition\\
\midrule
Bounded; conditional sub-Gaussian or sub-exponential
 & \(n^{-1/2}\{\log(ed)\}^{3/2}\)
 & \(\log(ed)=o(n^{1/3})\)\\
Marginal sub-Gaussian
 & \(n^{-1/2}\{\log(ed)\}^{2}\)
 & \(\log(ed)=o(n^{1/4})\)\\
Marginal sub-exponential
 & \(n^{-1/2}\{\log(ed)\}^{5/2}\)
 & \(\log(ed)=o(n^{1/5})\)\\
Uniform moment bound \eqref{main:eq:moment-envelope}, \mbox{\(2<q\le3\)}
 & \(\big[n^{1-q/2}\{\log(ed)\}^{q/2}\big]^{1/(q+1)}\)
 & \(\log(ed)=o(n^{1-2/q})\)\\
\bottomrule
\end{tabular}
\end{minipage}
\end{center}

The first three rows distinguish conditional from marginal tail
control. Under a conditional exponential-moment bound, large increments
remain rare after conditioning on the past, and the dimension threshold
is the same as under boundedness. A marginal bound permits the increment
law to change substantially on a rare event determined by earlier
increments. The constructions in Proposition~\ref{low:cubic} use this freedom even
when each increment has conditional covariance \(I_d\), and hence
\(Q=I_d\). The difference between the thresholds therefore persists
without fluctuations of \(Q\).

The dimension condition in each row, together with
\(\delta_n\log(ed_n)\to0\), is necessary and sufficient for
uniform convergence. For every fixed \(q>2\), the finite-moment
criterion has dimension exponent \(\min\{1/3,1-2/q\}\);
the quantitative upper and lower bounds cover \(2<q\le3\)
and sufficiently large \(n\).
These full convergence criteria are stated separately in
Theorem~\ref{main:validity}. Their sufficiency requires a separate
argument, since a bound containing powers of \(\log n\) need not cover
every slowly subcritical sequence.

When \(Q=I_d\) almost surely, we can also remove the additional powers
of \(\log n\) from several quantitative bounds. The bounded and marginal
light-tail rates then hold up to constants for
\(\log(ed)\ge\{\log(en)\}^2\). The finite-moment rate has the
same optimal order for \(2<q\le3\), in every dimension and for all
sufficiently large \(n\). Theorem~\ref{main:identity} gives these
refinements after the common rate and convergence results.

\paragraph{Relation to earlier work.}
The independent theory provides a direct benchmark.
The rectangle CLTs of \citet{CCK2013,CCK2017} initiated a line
of work on dimension dependence and covariance geometry.
For fixed bounded increments, bounding the Gaussian expected
maximum in \citet[arXiv v1, Theorem~1]{CCKK2026} by
\(C\sqrt{\log(ed)}\) gives
\(n^{-1/2}\{\log(ed)\}^{3/2}\sqrt{\log(en)}\) under unrestricted
covariance. Their Corollary~4 removes the sample
logarithm under a fixed absolute pairwise correlation gap.
Here we allow martingale dependence and random \(Q\), under
the correlation condition in \eqref{main:eq:setting}. We compare the
two correlation assumptions after the main statements.

For martingales, \citet[Theorem~2.1]{belloni2018high} compare smooth
test functions when the quadratic variation is close to a deterministic
matrix; their assumption also requires a deterministic positive
semidefinite upper bound with high probability. They then obtain
approximation results for maxima. Related multivariate bounds include
\citet{kojevnikov2022berry}; coupling and Wasserstein bounds are developed
by \citet{cattaneo2025yurinskii} and \citet{fang2026martingale},
respectively. Our covariance error is measured directly by the
\(L^p\) norm of \(\|Q-\Sigma\|_{\max}\).

A close comparison under a third-moment bound is
\citet[arXiv v3, Theorem~2.2]{wu2026berry}. When \(Q=I_d\) and
\(\max_i\E\|X_i\|_\infty^3\) is bounded by a fixed constant, their
rectangle bound is \(Cn^{-1/8}\{\log(ed)\}^{9/8}\).
Our upper and lower bounds have order
\(1\wedge n^{-1/8}\{\log(ed)\}^{3/8}\) in this class.
The \(n^{-1/4}\) conclusion of their Theorem~2.5 additionally
requires a positive lower bound on the smallest eigenvalue of every
conditional covariance matrix, as well as its stated moment and
conditional variance bounds. This distinction matters because our
individual conditional covariance matrices may be singular.
The scalar exponent \(1/8\) is already known
\citep{haeusler1988rate,Joos1993,rinott2000bounds}; here the
finite-moment bounds identify the dependence on dimension.

The scalar exponent
\((2+1/p)^{-1}\) is classical
\citep{haeusler1988rate,mourrat2013rate}; the bounds here determine
its dimension factor and combine it with the increment-dependent
terms in Table~\ref{main:tab:scales}. The lower bounds adapt the scalar constructions
of \citet{Joos1993}.

\paragraph{Proof ideas.}
The first difficulty is to sum local replacement errors when
conditional variances vary over time. We retain the conditional
covariance of each increment as a weight in the local Gaussian
comparison. A small-variance increment then incurs a correspondingly
small cost, and the weights can be summed along the predictable
covariance path. This develops the residual-variance replacement
of \citet[Theorem~2.1]{Rollin2018} and the covariance-matched
convolution method of \citet[Proposition~2.1]{fang2026martingale}.

When \(Q\) is random, the matrix needed to complete a partial sum
to covariance \(\Sigma\) need not be positive semidefinite. It therefore
cannot in general serve as the covariance of an added Gaussian vector.
We instead convolve with a Gaussian whose covariance is positive
semidefinite and account for the remaining symmetric matrix by a finite
sum of derivatives of the smoothed test function. When the two matrices
change while their sum stays fixed, their contributions cancel through
the chosen derivative order. Estimates on rectangle boundaries control
the remainder and yield the term involving \(\delta\).

To prove sufficiency throughout the convergence range, we need error
bounds that remain useful for rare events. We first compare the
probabilities that a fixed number of coordinate blocks violate their
rectangle constraints. We bound these errors relative to products of
selected one-dimensional Gaussian tail probabilities, using exponential
tilting for light tails and positive smoothing under moment bounds.
The correlation condition in \eqref{main:eq:setting} then controls the
sum over choices of blocks. Using one-dimensional tails avoids dividing
by a possibly zero joint probability when a block covariance is singular.

Section~\ref{main:results} states all main results.
Sections~\ref{core:section} and~\ref{thr:sec:geometry} develop the
predictable comparisons and Gaussian geometry.
Sections~\ref{quant:section} and~\ref{quant:identity-section} prove
the light-tail quantitative bounds and the identity refinement;
Section~\ref{thr:section} proves the complete light-tail convergence
range. Section~\ref{poly:section} proves the finite-moment results,
including their lower bound for the increment-dependent term. The remaining lower
constructions and the proofs of the three main theorems are in
Sections~\ref{low:section} and~\ref{main:proofs}.
Appendices~\ref{app:logs} and~\ref{low:scalar-section} give
refinements of the powers of \(\log n\) and a scalar example showing
that some logarithmic factors are necessary.

\section{Main results}\label{main:results}

Throughout this section the arrays satisfy
\eqref{main:eq:normalization}--\eqref{main:eq:setting}, with
fixed \(s>0\), \(M\ge1\) and \(1\le p\le\infty\).
The increment constants \(B\ge2\), \(K\ge4\), and the moment
order \(q>2\), when used, are fixed as well.
We first state upper and lower bounds that differ only by powers of
\(\log n\), then necessary and sufficient conditions for convergence,
and finally bounds optimal up to constants when \(Q=I_d\).

For bounded increments, put
\begin{equation}\label{main:eq:D}
 D_{n,d}(\delta)=
 \sup_{\substack{\text{arrays satisfying
           \eqref{main:eq:normalization}--\eqref{main:eq:setting}}\\
                    \|X_i\|_\infty\le B\ \mathrm{a.s.}}}
 \rho_{\mathcal R}(S_n,G_\Sigma).
\end{equation}
Replacing boundedness by \eqref{main:eq:conditional-tails},
\eqref{main:eq:marginal-tails}, or
\eqref{main:eq:moment-envelope} defines, respectively,
\[
 D^{\mathrm{cond},\alpha}_{n,d}(\delta),\qquad
 D^{\mathrm{marg},\alpha}_{n,d}(\delta),\qquad
 D^q_{n,d}(\delta).
\]
All these suprema allow arbitrary filtered probability spaces and
initial sigma-fields. The statements hold for any consistent choice
of open or closed finite rectangle endpoints. Monotone approximation
and continuity of the Gaussian marginals give the same supremum.

\Needspace{9\baselineskip}
\subsection{Optimal approximation rates}

\begin{theorem}[Optimal rates up to powers of \(\log n\)]\label{main:joint}
The following bounds hold.
\begin{enumerate}[label=\textup{(\roman*)},leftmargin=*]
\item For bounded increments, there are \(c,C>0\), depending only
on \(B,s,M,p\), such that for every \(n\ge2\), \(d\ge1\) and
\(\delta\ge0\),
\begin{equation}\label{main:eq:rates}
\begin{aligned}
 D_{n,d}(\delta)
 &\ge c\left[1\wedge\left\{
 \frac{\{\log(ed)\}^{3/2}}{\sqrt n}
 +\{\delta\log(ed)\}^{1/(2+1/p)}\right\}\right],\\
 D_{n,d}(\delta)
 &\le 1\wedge C\{\log(en)\}^{2+3\max\{2,s\}}
 \left\{\frac{\{\log(ed)\}^{3/2}}{\sqrt n}
 +\{\delta\log(ed)\}^{1/(2+1/p)}\right\}.
\end{aligned}
\end{equation}
For either \(\alpha\in\{1,2\}\), both inequalities hold with
\(D^{\mathrm{cond},\alpha}_{n,d}(\delta)\) in place of
\(D_{n,d}(\delta)\), with constants depending only on \(K,s,M,p\).

\item For marginal exponential moments and either
\(\alpha\in\{1,2\}\), there are \(c,C>0\), depending only on
\(K,s,M,p\), such that for every \(n\ge2\), \(d\ge1\) and
\(\delta\ge0\),
\begin{equation}\label{main:eq:marginal-rates}
\begin{aligned}
 D^{\mathrm{marg},\alpha}_{n,d}(\delta)
 &\ge c\left[1\wedge\left\{
 \frac{\{\log(ed)\}^{3/2+1/\alpha}}{\sqrt n}
 +\{\delta\log(ed)\}^{1/(2+1/p)}\right\}\right],\\
 D^{\mathrm{marg},\alpha}_{n,d}(\delta)
 &\le 1\wedge C\{\log(en)\}^{2+3\max\{2,s\}}
 \left\{\frac{\{\log(ed)\}^{3/2+1/\alpha}}{\sqrt n}
 +\{\delta\log(ed)\}^{1/(2+1/p)}\right\}.
\end{aligned}
\end{equation}

\item For \(2<q\le3\), there are \(c,C>0\) and \(n_0\ge2\),
depending only on \(q,B,s,M,p\), such that for every \(n\ge n_0\),
\(d\ge1\) and \(\delta\ge0\),
\begin{equation}\label{main:eq:moment-rates}
\begin{aligned}
 D^q_{n,d}(\delta)
 &\ge c\left[1\wedge\left\{
 [n^{1-q/2}\{\log(ed)\}^{q/2}]^{1/(q+1)}
 +\{\delta\log(ed)\}^{1/(2+1/p)}\right\}\right],\\
 D^q_{n,d}(\delta)
 &\le 1\wedge C\{\log(en)\}^{2+3\max\{2,s\}}\\[-2pt]
 &\qquad\times\left\{[n^{1-q/2}\{\log(ed)\}^{q/2}]^{1/(q+1)}
 +\{\delta\log(ed)\}^{1/(2+1/p)}\right\}.
\end{aligned}
\end{equation}
\end{enumerate}
\end{theorem}

Separate examples give the lower bounds for the increment-dependent
term and for the term involving \(\delta\). The maximum of these two
bounds is comparable to their sum, and all the examples have
\(\Sigma=I_d\). Thus the powers of \(\log d\) and \(\delta\), and
the polynomial rate in \(n\), cannot be improved uniformly over the
stated classes. We do not claim that the displayed power of
\(\log n\) is optimal; Appendix~\ref{app:logs} reduces it in several
regimes.

\Needspace{9\baselineskip}
\subsection{Necessary and sufficient conditions for convergence}

\begin{theorem}[Uniform convergence]\label{main:validity}
For every pair of deterministic sequences \(d_n\ge1\) and
\(\delta_n\ge0\), the following equivalences hold as \(n\to\infty\).
\begin{enumerate}[label=\textup{(\roman*)},leftmargin=*]
\item For bounded increments,
\begin{equation}\label{main:eq:validity}
 D_{n,d_n}(\delta_n)\to0
 \quad\Longleftrightarrow\quad
 \log(ed_n)=o(n^{1/3}),\qquad
 \delta_n\log(ed_n)\to0.
\end{equation}
The same equivalence holds with
\(D^{\mathrm{cond},\alpha}_{n,d_n}(\delta_n)\) in place of
\(D_{n,d_n}(\delta_n)\), for either \(\alpha\in\{1,2\}\).

\item For marginal exponential moments and either
\(\alpha\in\{1,2\}\),
\begin{equation}\label{main:eq:marginal-validity}
 \begin{split}
 D^{\mathrm{marg},\alpha}_{n,d_n}(\delta_n)\to0
 \quad\Longleftrightarrow\quad&
 \log(ed_n)=o(n^{1/(3+2/\alpha)}),\\[-2pt]
 &\delta_n\log(ed_n)\to0.
 \end{split}
\end{equation}

\item For every fixed \(q>2\), under the moment condition
\eqref{main:eq:moment-envelope},
\begin{equation}\label{main:eq:moment-validity}
 \begin{split}
 D^q_{n,d_n}(\delta_n)\to0
 \quad\Longleftrightarrow\quad&
 \log(ed_n)=o\!\left(n^{\min\{1/3,1-2/q\}}\right),\\[-2pt]
 &\delta_n\log(ed_n)\to0.
 \end{split}
\end{equation}
\end{enumerate}
\end{theorem}

The lower bounds in Theorem~\ref{main:joint} give necessity in
its quantitative ranges; the bounded construction supplies the
finite-moment dimension obstruction when \(q>3\).
Sufficiency throughout the stated ranges follows from the relative
comparisons in Sections~\ref{thr:section} and~\ref{poly:section}.
These arguments control sequences approaching the dimension
threshold arbitrarily slowly.

\subsection{Optimal rates when \texorpdfstring{\(Q=I_d\)}{Q equals Id}}

Restrict the bounded, marginal exponential, and finite-moment
classes to \(Q=I_d\) almost surely, and write their maximal errors as
\[
 D^{I,\mathrm{bd}}_{n,d},\qquad
 D^{I,\mathrm{marg},\alpha}_{n,d},\qquad
 D^{I,q}_{n,d},
\]
respectively. The condition fixes the sum of conditional covariance
matrices. The individual matrices may still be
random, singular, and noncommuting.

\begin{theorem}[Optimal rates when \(Q=I_d\)]\label{main:identity}
The following comparisons hold with constants independent of \(n,d\).
\begin{enumerate}[label=\textup{(\roman*)},leftmargin=*]
\item For every \(n\ge2\) and \(d\ge1\) such that
\(\log(ed)\ge\{\log(en)\}^2\),
\begin{align}
 D^{I,\mathrm{bd}}_{n,d}
 &\asymp_B 1\wedge
       \frac{\{\log(ed)\}^{3/2}}{\sqrt n},
       \label{main:eq:identity-bounded}\\
 D^{I,\mathrm{marg},\alpha}_{n,d}
 &\asymp_K 1\wedge
       \frac{\{\log(ed)\}^{3/2+1/\alpha}}{\sqrt n},
       \qquad \alpha\in\{1,2\}.
       \label{main:eq:identity-marginal}
\end{align}

\item For \(2<q\le3\), every \(d\ge1\), and all sufficiently
large \(n\),
\begin{equation}\label{main:eq:identity-polynomial}
 D^{I,q}_{n,d}\asymp_{q,B}
 1\wedge[n^{1-q/2}\{\log(ed)\}^{q/2}]^{1/(q+1)}.
\end{equation}
The threshold on \(n\) depends only on \(q,B\).
\end{enumerate}
\end{theorem}

These upper bounds retain the conditional Gaussian probability that
the remaining coordinates satisfy their interval constraints after selected
coordinates have been fixed at their boundary values. This removes the additional
powers of \(\log n\) in the stated ranges; setting \(\delta=0\) in
Theorem~\ref{main:joint} alone does not give these conclusions.
The coordinates of \(S_n\) need not be independent.

\subsection{The covariance assumption and related results}

\begin{remark}\label{main:rem:geometry}
A fixed operator-norm bound \(\|\Sigma\|_{\mathrm{op}}\le M_0\)
implies the correlation condition in \eqref{main:eq:setting} with \(s=2\):
\(\Sigma^2\preceq M_0\Sigma\) gives
\(\sum_k|\Sigma_{jk}|^2\le M_0\).
The same condition also permits singular fixed-size all-one
blocks. For \(s>1\), the matrices
\[
 \Sigma=(1-d^{-1/s})I_d+d^{-1/s}\ind\ind^\top
\]
have sums \(\sum_k|\Sigma_{jk}|^s\) at most two and largest eigenvalue
\(1-d^{-1/s}+d^{1-1/s}\).
Thus the condition permits some growing largest eigenvalues.
The convergence criteria also hold on the classes restricted by
a fixed operator-norm bound or by \(\Sigma=I_d\), since the
lower-bound examples remain in those classes.

As structural assumptions, the correlation condition in
\eqref{main:eq:setting} and a fixed
absolute pairwise correlation gap are incomparable. Fixed-size
all-one blocks satisfy the former and violate the latter.
A fixed positive equicorrelation satisfies the latter but has
\(\sum_k|\Sigma_{jk}|^s\to\infty\) as \(d\to\infty\), for every
fixed \(s>0\).
\end{remark}

For independent sums, \citet[Theorem~2(i)]{CCKK2026} treat every
fixed \(q>2\) under a moment bound on the largest absolute coordinate and
allow its scale to grow. Their part~(ii) adds a coordinate
fourth-moment condition. Our results keep the moment bound fixed and allow martingale dependence.
Earlier improvements in the independent rectangle CLT include
\citet{DengZhang2020,Koike2021,FangKoike2021,CCKK2022,Lopes2022,CCK2023};
recent treatments of degenerate or unrestricted covariance include
\citet{FKLZ2026,GaoZhang2026,LopesMasini2026,CCKK2026}.
The proofs below cite the particular Gaussian estimates and smoothing
arguments at their points of use.

In the proofs, \(C,c>0\) may change from line to line. Their
dependence is restricted to the fixed parameters in the statement
under consideration. We use normalized increments
\(\xi_i=X_i/\sqrt n\), with
\(C_i=\E_{i-1}\xi_i\xi_i^\top\) and
\(P_i=\sum_{j\le i}C_j\), so \(P_n=Q\).
Stopped and truncated arrays are compared with a specified
deterministic target covariance. It need not equal the covariance of
the modified sum.

\section{Gaussian comparison at predictable times}\label{core:section}

The local comparison error must decrease with the conditional covariance
of the increment being replaced. We first construct a transport whose
error bound is linear in a trace involving that covariance. Applying
these transports backwards through the array then reduces the total
error to a telescoping sum of log-determinants. We also prove a summation
lemma for predictable weights, a covariance completion lemma, and
estimates under exponential tilting. These will be used to compare
probabilities of rare events.

Covariance-matched Gaussian convolution and residual-variance
replacement appear in \citet[Proposition~2.1]{fang2026martingale}
and \citet[Theorem~2.1]{Rollin2018}. A coupling that is measurable in
its parameters is constructed in
\citet[Proposition~3.1]{FangKoike2024}. We prove the trace bound below
for bounded increments and construct a jointly Borel transport, so that
the same map can be applied to conditional distributions.

For a positive definite matrix \(A\), let \(\phi_A\) denote the density of
\(N(0,A)\). We write \(\phi=\phi_1\) for the standard normal
density and \(\Phi\) for its distribution function.
For \(r\ge1\), \(W_{r,\infty}\) is the infimum of
\(\{\E\|U-V\|_\infty^r\}^{1/r}\) over couplings of the two laws.
All norms of vectors without a subscript
in this section are Euclidean.

\subsection{A Gaussian convolution estimate with a covariance weight}

\begin{lemma}\label{core:lem:convolution}
Let \(\nu\) be a centered probability measure supported on
\([-b,b]^k\), where \(b>0\), and put \(V=\int zz^\top\nu(dz)\).
For every \(\Gamma\succeq kb^2I_k\), there is a map
\(T_{\nu,\Gamma}:\R^k\to\R^k\) such that
\[
 T_{\nu,\Gamma}(\xi+G_\Gamma)\sim N(0,\Gamma+V),
 \qquad \xi\sim\nu,\quad \xi\ \text{independent of }G_\Gamma.
\]
The map is jointly Borel in \((\nu,\Gamma,x)\), with the weak topology
on probability measures on the cube. If \(r\ge1\) and
\(\Gamma\succeq rkb^2I_k\), then the same map satisfies
\begin{equation}\label{core:eq:local}
 \left\{\E\|\xi+G_\Gamma-T_{\nu,\Gamma}(\xi+G_\Gamma)\|_\infty^r
 \right\}^{1/r}
 \le 17rb\,\tr\bigl((\Gamma+V)^{-1}V\bigr).
\end{equation}
\end{lemma}

\begin{proof}
We interpolate between two distributions with the same covariance.
For \(0\le t\le1\), define
\[
 R_t=\Gamma+(1-t)V,\qquad
 p_t(x)=\E\phi_{R_t}(x-\sqrt t\,\xi).
\]
Thus \(p_0=\phi_{\Gamma+V}\), and \(p_1\) is the density of
\(\xi+G_\Gamma\). Since \(V\preceq kb^2I_k\preceq\Gamma\), every
Gaussian covariance \(R_t\) is positive definite.

Differentiation of the Gaussian density, followed by Taylor expansion,
gives
\begin{align}
 \partial_t p_t&=-\operatorname{div}J_t,\notag\\
 J_t(x)&=\frac1{2\sqrt t}\E\{\xi\phi_{R_t}(x-\sqrt t\,\xi)\}
                   +\frac12V\nabla p_t(x)\notag\\
 &=\frac{\sqrt t}{2}\E\int_0^1
       \bigl\{(1-u)\xi\xi^\top-V\bigr\}
       \nabla^2\phi_{R_t}(x-u\sqrt t\,\xi)\xi\,du.
 \label{core:eq:flux}
\end{align}
Indeed, the term of first order in the first expectation is
\(-\sqrt t\,V\nabla\phi_{R_t}\), and it cancels the constant term in
\(V\nabla p_t\). Bounded support and the positive lower bound on \(R_t\)
justify differentiation locally uniformly in \(x\). The last expression
in \eqref{core:eq:flux} also extends continuously to \(t=0\).

Suppose now that \(\Gamma\succeq rkb^2I_k\). Write
\(\sigma_t=\tr(R_t^{-1}V)\). Jensen's inequality gives
\begin{equation}\label{core:eq:density-lower}
 \frac{p_t(x)}{\phi_{R_t}(x)}
 =\E\exp\left\{\sqrt t\,x^\top R_t^{-1}\xi
                 -\frac t2\xi^\top R_t^{-1}\xi\right\}
 \ge e^{-t\sigma_t/2}\ge e^{-1/(2r)}.
\end{equation}
To estimate the flux, fix \(R=R_t\), let \(\xi'\) be an independent
copy of \(\xi\), and put \(a=R^{-1/2}\xi\), \(a'=R^{-1/2}\xi'\).
Both have norm at most \(r^{-1/2}\). For fixed \(a,a'\), set
\[
 H_{a',a}(z)=(a'\cdot z)(a\cdot z)-a'\cdot a.
\]
If \(Z\sim N(0,I_k)\) and \(0\le s\le1\), completing the square yields
\[
 \left\|e^{sZ\cdot a-s^2|a|^2/2}H_{a',a}(Z-sa)\right\|_r
 =
 e^{(r-1)s^2|a|^2/2}
 \left\|H_{a',a}(Z+(r-1)sa)\right\|_r.
\]
The one-dimensional Gaussian moment bound
\(\{\E|N(0,1)|^{2r}\}^{1/(2r)}\le2\sqrt r\) and H\"older's
inequality therefore bound the right-hand side by
\[
 e^{(r-1)/(2r)}(9r+1)|a'||a|.
\]
The Hessian formula, with \(z=R^{-1/2}x\), gives
\begin{align*}
 \xi\xi^\top\nabla^2\phi_R(x-s\xi)\xi
 &=\xi\phi_R(x)e^{sz\cdot a-s^2|a|^2/2}H_{a,a}(z-sa),\\
 V\nabla^2\phi_R(x-s\xi)\xi
 &=\E_{\xi'}\{
 \xi'\phi_R(x)e^{sz\cdot a-s^2|a|^2/2}H_{a',a}(z-sa)\}.
\end{align*}
After averaging over \(\xi\) and taking the \(L^r\) norm under
\(\phi_R\), the first expression is bounded by a constant times
\(rb\,\E|a|^2\). For the second, the corresponding factor is
\[
 b\,\E|a'||a|=b(\E|a|)^2\le b\,\E|a|^2.
\]
In particular, the covariance enters linearly in both terms.
Combining these estimates with \eqref{core:eq:density-lower} and
\(\int_0^1(1-u)\,du=1/2\), we obtain
\begin{equation}\label{core:eq:velocity}
 \left\{\int\left\|\frac{J_t(x)}{p_t(x)}\right\|_\infty^r
                      p_t(x)\,dx\right\}^{1/r}
 \le \frac34\,17rb\sqrt t\,\tr(R_t^{-1}V).
\end{equation}
Here the combined exponential factor is at most
\(\exp\{(r^2-1)/(2r^2)\}\), and
\(10e^{1/2}<17\). This also proves the estimate when \(r=1\).

We next verify that this velocity produces a measurable transport.
Set \(v_t=J_t/p_t\) for \(t>0\), and \(v_0=0\). The last expression
in \eqref{core:eq:flux} shows that \(v\) and its spatial derivative
are continuous at time zero. Define
\[
 m_t(x)=
 \frac{\E\{\xi\phi_{R_t}(x-\sqrt t\,\xi)\}}{p_t(x)}.
\]
Differentiating this posterior mean gives
\[
 Dm_t(x)=\sqrt t\,\Cov(\xi\mid x)R_t^{-1},
\]
where the conditional distribution is the probability measure
proportional to \(\phi_{R_t}(x-\sqrt t\,\xi)\nu(d\xi)\). Hence
\[
 Dv_t(x)=
 \frac12(I_k+tVR_t^{-1})\Cov(\xi\mid x)R_t^{-1}
 -\frac12VR_t^{-1}.
\]
The posterior distribution is supported on the same cube, so
\(\|\Cov(\xi\mid x)\|_{\mathrm{op}}\le kb^2\). Under the weaker
condition \(\Gamma\succeq kb^2I_k\), it follows that
\(\|Dv_t(x)\|_{\mathrm{op}}\le3/2\), uniformly in \(t,x\).

The equation
\[
 \frac{d}{dt}\Phi_t(z)=v_t(\Phi_t(z)),\qquad \Phi_0(z)=z,
\]
has a unique global solution on \([0,1]\). The uniform Lipschitz
bound gives uniqueness and local existence by Picard iteration, while
\[
 |\Phi_t(z)|\le e^{3t/2}
       \left\{|z|+\int_0^t|v_s(0)|\,ds\right\}
\]
rules out escape in finite time. Solving backwards gives the inverse
flow. The spatial derivative satisfies
\[
 \det D\Phi_t(z)=
 \exp\left\{\int_0^t\operatorname{div}v_s(\Phi_s(z))\,ds\right\}>0.
\]
The continuity equation in \eqref{core:eq:flux}, applied along the
flow and then continued to time zero, gives
\[
 p_t(\Phi_t(z))\det D\Phi_t(z)=\phi_{\Gamma+V}(z).
\]
Thus \(\Phi_t\) sends \(N(0,\Gamma+V)\) to density \(p_t\).
This proves the distributional assertion with
\(T_{\nu,\Gamma}=\Phi_1^{-1}\).

For \(Z_{\Gamma+V}\sim N(0,\Gamma+V)\), Minkowski's integral
inequality and \eqref{core:eq:velocity} give
\begin{align*}
 \|\|\Phi_1(Z_{\Gamma+V})-Z_{\Gamma+V}\|_\infty\|_r
 &\le \frac34\,17rb
       \int_0^1\sqrt t\,\tr(R_t^{-1}V)\,dt\\
 &\le17rb\,\tr((\Gamma+V)^{-1}V).
\end{align*}
The last step uses \(R_t\succeq\Gamma\succeq(\Gamma+V)/2\) and
\(\int_0^1\sqrt t\,dt=2/3\).

Finally, \(\nu\mapsto V\) is continuous on the probability measures
on the compact cube. The last expression in \eqref{core:eq:flux} and
the positive lower bound on the eigenvalues of \(R_t\) imply that
\((\nu,\Gamma,t,x)\mapsto v_t(x)\) is continuous,
locally uniformly, including at \(t=0\). The common Lipschitz bound
and Gronwall's inequality give continuous dependence of the forward
and backward flows on the parameters on every compact set.
Consequently \(T_{\nu,\Gamma}(x)\) is jointly Borel. If \(V=0\),
the law is concentrated at zero and the map is the identity.
\end{proof}

\subsection{Backward coupling and summation of covariance weights}

\begin{proposition}\label{core:prop:transport}
Let \((\xi_i,\mathcal F_i)_{i=1}^N\) be martingale differences in
\(\R^k\), on any filtered probability space, such that
\[
 \|\xi_i\|_\infty\le b,\qquad
 C_i=\E_{i-1}\xi_i\xi_i^\top,\qquad
 \sum_{i=1}^NC_i=\Gamma\quad\text{a.s.},
 \qquad \tr\Gamma\le2k.
\]
For every \(r\ge1\),
\begin{align}
 W_{r,\infty}\left(\sum_{i=1}^N\xi_i,G_\Gamma\right)
 &\le
 17rbk\log\left(1+\frac{2}{rkb^2}\right)
 +Cb\sqrt{rk\{r+\log(2k)\}}\notag\\
 &\le C_k rb\log(e+1/b).
 \label{core:eq:transport}
\end{align}
Here \(\Gamma\) is deterministic. The constants do not depend on
the length \(N\), the filtration or the smallest eigenvalue of
\(\Gamma\).
\end{proposition}

\begin{proof}
Put
\[
 P_i=\sum_{j\le i}C_j,\qquad
 \Lambda=rkb^2I_k,\qquad A_i=\Gamma-P_i+\Lambda.
\]
The matrix \(\Gamma-P_i\) is \(\mathcal F_{i-1}\)-measurable and
equals \(\sum_{j>i}C_j\succeq0\) almost surely. In particular,
\(A_i\succeq\Lambda\) and \(A_N=\Lambda\).
All the following identities hold on a common event of probability one.
For the pointwise definition of the maps, replace
\(\Gamma-P_i\) by its positive part on the exceptional set.

Let \(\nu_i\) be a conditional law of \(\xi_i\) given
\(\mathcal F_{i-1}\), chosen centered everywhere. Such a kernel exists
without any regularity assumption on the underlying probability
space. Indeed, take conditional expectations of a countable dense
unital rational vector space in \(C([-b,b]^k)\), containing the
coordinate functions and their products. Outside a common measurable
null set, they define a positive norm-one functional. Its continuous
extension and the Riesz representation theorem give \(\nu_i\).
Uniform approximation followed by a monotone class argument verifies
the conditional-law identity for bounded Borel functions. On the
exceptional set use \(\delta_0\).

Enlarge the space by an independent \(Y_N\sim N(0,\Lambda)\).
Starting from \(i=N\), define backwards
\[
 U_i=\xi_i+Y_i,\qquad
 Y_{i-1}=T_{\nu_i,A_i}(U_i),\qquad
 D_i=U_i-Y_{i-1}.
\]
We prove by backward induction that
\begin{equation}\label{core:eq:backward}
 \mathcal L(Y_i\mid\mathcal F_i)=N(0,A_i).
\end{equation}
To check the induction, for any bounded Borel \(f\),
\begin{align*}
 \E_{i-1}f(\xi_i+Y_i)
 &=\E_{i-1}\int f(\xi_i+y)\,N(0,A_i)(dy)\\
 &=\iint f(x+y)\,N(0,A_i)(dy)\nu_i(dx).
\end{align*}
The first equality conditions on \(\mathcal F_i\), and the second uses
the predictability of \(A_i\). The same identity holds for
parameter-dependent Borel tests by a monotone class argument.
Lemma~\ref{core:lem:convolution} therefore proves
\eqref{core:eq:backward} at \(i-1\) and gives
\begin{equation}\label{core:eq:conditional-cost}
 \E_{i-1}\|D_i\|_\infty^r\le a_i^r,\qquad
 a_i=17rb\,\tr(A_{i-1}^{-1}C_i).
\end{equation}
In particular, \(Y_0\sim N(0,\Gamma+\Lambda)\), also conditionally
on the possibly nontrivial initial sigma-field.

The weights \(a_i\) have a deterministic total bound. Since
\(0\preceq A_{i-1}^{-1/2}C_iA_{i-1}^{-1/2}\prec I_k\),
\[
 \tr(A_{i-1}^{-1}C_i)
 \le \log\det A_{i-1}-\log\det A_i.
\]
This follows from \(-\log(1-x)\ge x\) and does not require the
matrices to commute. Concavity of \(\log\), applied to the eigenvalues
of \(\Gamma\), now yields
\begin{equation}\label{core:eq:weight-budget}
 \sum_{i=1}^Na_i
 \le17rb\log\frac{\det(\Gamma+\Lambda)}{\det\Lambda}
 \le17rbk\log\left(1+\frac{2}{rkb^2}\right).
\end{equation}
The same log-determinant summation appears in
\citet[arXiv v3, Appendix~B.2, equation~(74), and Corollary~3.4]{wu2026uncertainty}.
The local conditional cost used here is \eqref{core:eq:conditional-cost}.

For \(r>1\), weighted H\"older gives
\[
 \left(\sum_i\|D_i\|_\infty\right)^r
 \le
 \left(\sum_i a_i\right)^{r-1}
 \sum_{i:a_i>0}\frac{\|D_i\|_\infty^r}{a_i^{r-1}}.
\]
Every weight is predictable, and \(D_i=0\) almost surely on
\(\{a_i=0\}\). Taking expectations and using
\eqref{core:eq:conditional-cost} proves that the \(L^r\) norm of
\(\sum_i\|D_i\|_\infty\) is at most the right-hand side of
\eqref{core:eq:weight-budget}. For \(r=1\), sum the conditional
first moments. The recursion telescopes:
\[
 \sum_iD_i=\sum_i\xi_i+Y_N-Y_0.
\]
It remains to remove the independent Gaussian terms added to the two
marginal laws.
For \(Z\sim N(0,I_k)\), the union bound and integration of its tail give
\[
 \|\|Z\|_\infty\|_r\le C\sqrt{r+\log(2k)}.
\]
The two independent \(N(0,\Lambda)\) regularizations therefore
cost at most \(Cb\sqrt{rk\{r+\log(2k)\}}\).
This gives the first bound in \eqref{core:eq:transport}.
The second follows for fixed \(k\).
\end{proof}

\subsection{Integration along predictable weights}

\begin{lemma}\label{core:lem:clock}
Let \(w_i\ge0\) be \(\mathcal F_{i-1}\)-measurable and suppose that
\(\sum_{i=1}^Nw_i\le W\) almost surely, where \(W\) is deterministic.
For nonnegative random variables \(K_i\),
\begin{equation}\label{core:eq:clock}
 \E\sum_{i=1}^Nw_iK_i
 \le
 W\sup_\tau\E K_\tau,
\end{equation}
where the supremum is over indices with
\(\{\tau=i\}\in\mathcal F_{i-1}\), allowing a cemetery value with
\(K_\infty=0\).
\end{lemma}

\begin{proof}
For \(0<t\le W\), let
\(\tau(t)=\inf\{i:\sum_{j\le i}w_j\ge t\}\), with value \(\infty\)
if the level is never reached. These are predictable indices. Tonelli's
theorem and the pathwise lengths of their level sets give
\[
 \E\sum_iw_iK_i
 =\int_0^W\E K_{\tau(t)}\,dt
 \le W\sup_\tau\E K_\tau.
\]
\end{proof}

We will apply this lemma to Gaussian densities evaluated on fixed
coordinate faces of a rectangle. The required expectation bound must
hold uniformly over predictable indices before the random weights are
summed. The weights and the densities need not be independent.

\subsection{Predictable completion under an entrywise covariance bound}

Predictable stopping followed by Gaussian completion is used in
\citet[Theorem~2.1 and Appendix~A.3]{belloni2018high} and
\citet[supplement, Lemma~SA.8]{cattaneo2025yurinskii}
(Lemma~B.8 in the long preprint). Here the completed covariance is
chosen within a prescribed entrywise distance of the target. The
original array is preserved whenever its predictable quadratic
variation already satisfies that bound.

\begin{lemma}\label{thr:lem:feasible-completion}
Let \(\xi_1,\ldots,\xi_n\) be martingale differences in \(\R^d\),
and let
\[
 P_i=\sum_{j\le i}\E_{j-1}\xi_j\xi_j^\top,\qquad P_0=0.
\]
Fix a deterministic positive semidefinite \(\Sigma\) with unit
diagonal and \(\eta\ge0\). On an extension of the probability space
there is a martingale array with terminal sum \(Y\) and predictable
quadratic variation \(A_*\) such that
\begin{equation}\label{thr:eq:feasible-agreement}
 \|A_*-\Sigma\|_{\max}\le\eta\quad\text{a.s.},\qquad
 Y=\sum_i\xi_i\quad\text{on }
       \{\|P_n-\Sigma\|_{\max}\le\eta\}.
\end{equation}
It is obtained by retaining a predictably chosen initial portion of the
original array and appending centered increments that are independent
conditionally on the information available at the end of the original
array.

If \(\|\xi_i\|_\infty\le b\), this can be done with \(n+nd\)
increments, each bounded by
\(\max\{b,\sqrt{(1+\eta)/n}\}\).
Alternatively the added noise can be split into \(n\) conditionally
Gaussian increments, each with coordinate variance at most
\((1+\eta)/n\). In particular, a conditional coordinate
\(\psi_\alpha\) bound at scale \(K/\sqrt n\), \(\alpha\in\{1,2\}\),
is preserved with a changed constant when \(\eta\le1\).
\end{lemma}

\begin{proof}
Let
\[
 \mathcal K_\eta=\{A\succeq0:\|A-\Sigma\|_{\max}\le\eta\},
 \qquad
 \mathcal H_\eta=\{P\succeq0:P\preceq A
                       \text{ for some }A\in\mathcal K_\eta\}.
\]
Both sets are compact. Indeed, their defining matrix inequalities
are closed, their diagonal entries are bounded by \(1+\eta\),
and positive semidefiniteness bounds every off-diagonal entry by
the corresponding diagonal entries. The second set is the projection
of a compact set. It contains zero and is downward closed in the
positive semidefinite order.

For each \(P\in\mathcal H_\eta\), choose
\(A(P)\in\mathcal K_\eta\) with \(A(P)\succeq P\). We require this
choice to be Borel and to satisfy \(A(P)=P\) when
\(P\in\mathcal K_\eta\). To construct it,
for each positive integer \(m\), minimize over \(A\in\mathcal K_\eta\)
\[
 \|A-P\|_{\mathrm F}^2+
       m\|(P-A)_+\|_{\mathrm F}^2 .
\]
The minimizer is unique by strict convexity and depends continuously
on \(P\), by compactness and continuity of the objective.
For \(P\in\mathcal H_\eta\), its limit as \(m\to\infty\) is the
unique minimizer of \(\|A-P\|_{\mathrm F}^2\) over
\(A\in\mathcal K_\eta\) with \(A\succeq P\). To verify this, compare
with any such matrix:
the penalty tends to zero along every convergent subsequence, and
the limit minimizes the first term over the feasible set.
Uniqueness gives convergence of the entire sequence.
Thus \(A(P)\) is Borel and \(A(P)=P\) on \(\mathcal K_\eta\).

Retain increment \(i\) if \(P_j\in\mathcal H_\eta\) for every
\(j\le i\), and discard it and all subsequent increments at the first
failure. The decision is \(\mathcal F_{i-1}\)-measurable. The sum
\(P\) of the retained conditional covariance matrices belongs to
\(\mathcal H_\eta\).
Put \(R=A(P)-P\succeq0\). With fresh independent signs append,
for each column \(j\), the \(n\) increments
\[
 n^{-1/2}R^{1/2}e_j\varepsilon_{j,\ell},
 \qquad 1\le\ell\le n .
\]
Their conditional covariance matrices sum to \(R\), and every
coordinate is bounded by
\(\sqrt{(1+\eta)/n}\).
For the Gaussian version append
\(n^{-1/2}R^{1/2}Z_\ell\), \(1\le\ell\le n\), with fresh
independent standard Gaussian vectors.

On \(\{P_n\in\mathcal K_\eta\}\), every partial sum \(P_i\) is
dominated by \(P_n\), so all original increments are retained,
\(P=P_n\), and \(R=0\). This proves
\eqref{thr:eq:feasible-agreement}.
The conditional tail assertion follows from the coordinate variance
bound for the Gaussian increments. An arbitrary initial sigma-field
does not affect any part of the construction.
\end{proof}

The set \(\mathcal H_\eta\) records whether a positive semidefinite
completion remains possible within the entrywise bound. This is why the
stopping rule tests membership in \(\mathcal H_\eta\): failure of
\(P_i\preceq\Sigma\) alone would not imply
\(\Sigma\preceq P_i\).

If
\(\|\|Q-\Sigma\|_{\max}\|_p\le\delta_n\) and
\(\delta_n\log(ed)\to0\), choose a deterministic \(\eta_n\) with
\begin{equation}\label{thr:eq:localization-choice}
 \eta_n\log(ed)\to0,\qquad
 \Pp(\|Q-\Sigma\|_{\max}>\eta_n)\to0.
\end{equation}
For example, for finite \(p\) one may take
\[
 \eta_n=\frac{\{\delta_n\log(ed)+n^{-1}\}^{1/2}}{\log(ed)}.
\]
Markov's inequality verifies the second condition; for \(p=\infty\)
take \(\eta_n=\delta_n\).
The same localization is available whenever
\(\log(ed)\|Q-\Sigma\|_{\max}\to0\) in probability.
We shall prove the comparison with the deterministic target \(\Sigma\)
under an almost sure bound on the covariance discrepancy. The
covariance of the completed sum need not equal \(\Sigma\).

\subsection{Exponential tilting in fixed dimension}

For the next two lemmas the dimension \(k\) is fixed. Later we will
apply them to projections onto a fixed number of coordinates. We use a
predictable exponential change of measure, as in the martingale
Cram\'er method \citep{fan2013cramer}, and derive the normal approximation
under the changed measure from Lemma~\ref{core:lem:convolution}.

\begin{lemma}\label{thr:lem:sparse-tilt}
Let \(x_1,\ldots,x_N\in\R^k\) be martingale differences, with any
finite \(N\), such that
\[
 \|x_i\|_\infty\le b,\qquad
 C_i=\E_{i-1}x_ix_i^\top,\qquad
 \sum_{i=1}^NC_i=\Gamma ,
 \qquad \max_j\Gamma_{jj}\le2,
\]
where \(\Gamma\succeq0\) is deterministic.
Let \(n\ge2\), \(b\ge n^{-1/2}\), \(T\ge1\), and
\(\|\lambda\|_\infty\le T\), with \(kbT\le1/4\).
Put
\[
 z_i=\lambda^\top x_i,\quad v_i=\lambda^\top C_i\lambda,\quad
 \ell_i=\frac{e^{z_i}}{\E_{i-1}e^{z_i}},\quad
 \Psi=\sum_i\log\E_{i-1}e^{z_i}.
\]
The product of the \(\ell_i\)'s defines a probability
\(\Pp_\lambda\). If
\(\mu_i=\E_\lambda[x_i\mid\mathcal F_{i-1}]\), then almost surely
\begin{align}
 |\Psi-\tfrac12\lambda^\top\Gamma\lambda|
 &\le Cb\|\lambda\|_1\lambda^\top\Gamma\lambda,
       \label{thr:eq:bounded-cumulant}\\
 \left\|\sum_i\mu_i-\Gamma\lambda\right\|_\infty
 &\le Cb\lambda^\top\Gamma\lambda.
       \label{thr:eq:bounded-drift}
\end{align}
For every fixed \(A>0\), under \(\Pp_\lambda\) there is a coupling
of \(\sum_i x_i\) with \(Z\sim N(0,\Gamma)\) such that
\begin{equation}\label{thr:eq:bounded-tilted-coupling}
 \Pp_\lambda\left(
 \left\|\sum_i x_i-\Gamma\lambda-Z\right\|_\infty
 >C_{k,A}b\{\log(en)^2+T\log(en)+T^2\}\right)
 \le n^{-A}.
\end{equation}
Constants do not depend on \(N\), on the smallest positive eigenvalue
of \(\Gamma\), or on the initial sigma-field.
\end{lemma}

\begin{proof}
Since \(|z_i|\le1/4\) and \(\E_{i-1}z_i=0\), Taylor's formula
gives
\[
 \E_{i-1}e^{z_i}=1+\tfrac12v_i+
       O(b\|\lambda\|_1v_i),\qquad
 \|\mu_i-C_i\lambda\|_\infty\le Cbv_i .
\]
The squared term when taking the logarithm is also bounded by
\(Cb\|\lambda\|_1v_i\), because
\(v_i\le b^2\|\lambda\|_1^2\).
Summation and \(\sum_i v_i=\lambda^\top\Gamma\lambda\) prove
\eqref{thr:eq:bounded-cumulant}--\eqref{thr:eq:bounded-drift}.
Moreover,
\[
 \E_{i-1}(\ell_i-1)^2\le Cv_i,\qquad
 V_i:=\E_\lambda[(x_i-\mu_i)(x_i-\mu_i)^\top
                         \mid\mathcal F_{i-1}]
 \preceq CC_i .
\]

Fix \(r\ge2\). Add a positive definite regularization to each remaining
covariance sum:
\[
 D=4rkb^2I_k,\qquad
 A_i=D+\Gamma-\sum_{j\le i}C_j,\qquad
 h_i=\tr(A_i^{-1}C_i).
\]
Apply Lemma~\ref{core:lem:convolution}, conditionally under
\(\Pp_\lambda\), to the centered increment \(x_i-\mu_i\), whose
coordinate bound is \(2b\), and a Gaussian vector with covariance \(A_i\).
The conditional \(L^r\) norm of the transport displacement is at most
\(Crb h_i\).
The resulting Gaussian covariance is \(A_i+V_i\); the backward
recursion requires \(A_i+C_i\).

To compare these two Gaussians, interpolate
\[
 B_t=A_i+(1-t)V_i+tC_i,\qquad 0\le t\le1.
\]
The linear velocity
\((C_i-V_i)B_t^{-1}y/2\) transports their Gaussian laws. Indeed,
\(B_t\) solves its covariance equation:
\[
 \frac12(C_i-V_i)B_t^{-1}B_t
 +\frac12B_tB_t^{-1}(C_i-V_i)
 =C_i-V_i=\frac{d}{dt}B_t.
\]
Thus no commutation of the covariance matrices is needed. Since
\[
 V_i-C_i=\E_{i-1}\{(\ell_i-1)x_ix_i^\top\}-\mu_i\mu_i^\top,
\]
Gaussian moments, Minkowski, and conditional Cauchy--Schwarz
bound the first term's \(L^r\) velocity by
\[
 Cb\sqrt r\,
 \{\E_{i-1}(\ell_i-1)^2\}^{1/2}
 \{\E_{i-1}|B_t^{-1/2}x_i|^2\}^{1/2}
 \le Cb\sqrt{rv_i h_i}.
\]
For the second term, use \(\|\mu_i\|_\infty\le b\) and
\[
 |B_t^{-1/2}\mu_i|
 =|B_t^{-1/2}\E_{i-1}\{(\ell_i-1)x_i\}|
 \le C\sqrt{v_i h_i}.
\]
Integrating this linear velocity gives the same bound on the
transport displacement.
Its coefficients are Borel functions of the conditional moments,
so it can be composed with the conditional convolution map.

Starting with a fresh \(N(0,D)\), run the backward construction in
the proof of Proposition~\ref{core:prop:transport}, composing these
two maps at each step. The conditional output covariance is
\(A_{i-1}=A_i+C_i\).
The centered increments are martingale differences under
\(\Pp_\lambda\), while the original \(C_i\)'s remain predictable
and sum to the same deterministic \(\Gamma\).
Since \(C_i\preceq kb^2I_k\preceq A_i\),
\[
 \sum_i h_i
 \le2\log\frac{\det(D+\Gamma)}{\det D}
 \le C_k\log(en),\qquad
 \sum_i\sqrt{v_i h_i}
 \le\{\lambda^\top\Gamma\lambda\,C_k\log(en)\}^{1/2}.
\]
The first inequality follows by applying
\(x\le2\log(1+x)\) to the eigenvalues of
\(A_i^{-1/2}C_iA_i^{-1/2}\), which are at most one.
The last bound uses \(\tr\Gamma\le2k\) and \(b\ge n^{-1/2}\).

As in Proposition~\ref{core:prop:transport}, apply weighted H\"older
before taking conditional expectations. The total displacement then
has \(L^r\) norm at most
\[
 C_k b\{r\log(en)+T\sqrt{r\log(en)}\}.
\]
Removing the two endpoint Gaussian regularizations costs at most
\(C_kb(r+\sqrt r)\). Add the pathwise drift
\eqref{thr:eq:bounded-drift} and take
\(r=\lceil(A+2)\log(en)\rceil\).
Markov's inequality at \(e\) times the \(L^r\) bound gives
\eqref{thr:eq:bounded-tilted-coupling}.
\end{proof}

\Needspace{20\baselineskip}
\begin{lemma}\label{thr:lem:conditional-tilt}
Let \(N\le3n\), and retain the deterministic covariance sum
\(\sum_i C_i=\Gamma\) and the diagonal bound from
Lemma~\ref{thr:lem:sparse-tilt}. Replace boundedness by
\[
 \E_{i-1}\exp\{(\sqrt n\,|x_{ij}|/K)^\alpha\}\le2,
 \qquad \alpha\in\{1,2\}.
\]
Fix \(k,K,A\) and \(T\ge1\). If \(\|\lambda\|_\infty\le T\) and
\(T\{\log(en)\}^{1/\alpha}/\sqrt n\) is sufficiently small,
the predictable exponential tilt is defined and
\begin{align}
 |\Psi-\tfrac12\lambda^\top\Gamma\lambda|
       &\le C_{k,K}n^{-1/2}T^3,
       \label{thr:eq:conditional-cumulant}\\
 \left\|\sum_i\E_\lambda[x_i\mid\mathcal F_{i-1}]
                      -\Gamma\lambda\right\|_\infty
       &\le C_{k,K}n^{-1/2}T^2 .
       \label{thr:eq:conditional-drift}
\end{align}
There is a coupling under this tilt with \(Z\sim N(0,\Gamma)\)
whose failure probability is at most \(n^{-A}\) and whose
maximum-coordinate displacement from
\(\sum_i x_i-\Gamma\lambda\) is at most
\begin{equation}\label{thr:eq:conditional-tilted-coupling}
 C_{k,K,A}n^{-1/2}
 \{\{\log(en)\}^{2+1/\alpha}+T\log(en)+T^2\}.
\end{equation}
\end{lemma}

\begin{proof}
Write \(M_i=\|x_i\|_\infty\). Integrating the conditional tail bound
gives, for each fixed \(r\ge1\),
\[
 \E_{i-1}M_i^r\le C_{r,k,K}n^{-r/2}.
\]
The corresponding conditional \(L^r\) norm is bounded by
\(C_{r,k,K}n^{-1/2}\).
For \(T/\sqrt n\) small, conditional H\"older also gives
\[
 \E_{i-1}|z_i|^3e^{|z_i|}\le C_{k,K}n^{-3/2}T^3,\qquad
 \E_{i-1}|x_{ij}||z_i|^2e^{|z_i|}
       \le C_{k,K}n^{-3/2}T^2 .
\]
For example, Jensen's inequality applied to the average of the absolute
coordinates bounds \(\E_{i-1}e^{c|z_i|}\) uniformly for every
fixed \(c\), after reducing the allowed value of \(T/\sqrt n\).
Taylor expansion of the conditional moment generating function and
its first derivative proves
\eqref{thr:eq:conditional-cumulant}--\eqref{thr:eq:conditional-drift}.
No truncation is used in these two estimates.

For the normal approximation under the tilt, fix a sufficiently
large constant \(C_1=C_1(k,K,A)\), and truncate at
\[
 b=C_1\{\log(en)\}^{1/\alpha}/\sqrt n,\qquad
 x_i'=x_i\ind_{\{M_i\le b\}} .
\]
The assumed smallness implies \(kbT\le1/4\), after adjusting the
fixed smallness constant. For any prescribed \(a>0\), increasing
\(C_1\) gives
\begin{equation}\label{thr:eq:tilted-truncation}
 \E_{i-1}\big[
 (1+M_i^4+|z_i|^4)e^{2|z_i|}
                    \ind_{\{M_i>b\}}\big]\le n^{-a}.
\end{equation}
Indeed
\(\Pp(M_i>u\mid\mathcal F_{i-1})
\le2k\exp\{-(\sqrt n\,u/K)^\alpha\}\).
On \(u\ge b\), the polynomial weights and \(e^{2kTu}\) are
absorbed into a fixed fraction of this exponential.
The tilt denominator is at least one by Jensen.
Thus the probability under \(\Pp_\lambda\) of any truncation is
at most \(3n^{1-a}\).

Let \(\mu_i'\) and \(V_i'\) be the tilted conditional mean and
centered covariance of \(x_i'\).
On the retained event \(\ell_i\le e^{1/4}\), so
\[
 V_i'\preceq CC_i,\qquad
 C_i\preceq C_{k,K}n^{-1}I_k\preceq kb^2I_k.
\]
Use the same regularized remaining covariance sums
\(A_i=4rkb^2I_k+\Gamma-\sum_{j\le i}C_j\) and weights
\(h_i=\tr(A_i^{-1}C_i)\).
The convolution estimate bounds the conditional \(L^r\) displacement
for the centered truncated increment by \(Crb h_i\).
For the covariance correction, decompose
\begin{align*}
 V_i'-C_i
 &=\E_{i-1}\{(\ell_i-1)x_ix_i^\top\ind_{\{M_i\le b\}}\}\\
 &\quad-\E_{i-1}\{x_ix_i^\top\ind_{\{M_i>b\}}\}
           -\mu_i'(\mu_i')^\top .
\end{align*}
The Gaussian transport for the first term has \(L^r\) displacement
bounded by
\[
 C\sqrt r
 \{\E_{i-1}[M_i^2(\ell_i-1)^2]\}^{1/2}\sqrt{h_i}
 \le C_{k,K}n^{-1}T\sqrt{r h_i}.
\]
For the last inequality use conditional fourth moments,
\(|e^z-1|\le |z|e^{|z|}\), and the uniform exponential moment
just established.
Furthermore
\[
 \|\mu_i'\|_\infty\le C_{k,K}n^{-1}T+n^{-a},
 \qquad
 |B_t^{-1/2}\mu_i'|
 \le C_{k,K}n^{-1/2}T\sqrt{h_i}+C_k\sqrt n\,n^{-a},
\]
where \(B_t=A_i+(1-t)V_i'+tC_i\).
The second estimate follows from
\[
 \mu_i'=
 \E_{i-1}\{(\ell_i-1)x_i\ind_{\{M_i\le b\}}\}
             -\E_{i-1}\{x_i\ind_{\{M_i>b\}}\}
\]
and \(B_t\succeq A_i\succeq n^{-1}I_k\).
Since \(T/\sqrt n\) is small, the rank-one term obeys the
preceding transport bound plus an error smaller than any prescribed
negative power of \(n\), by increasing \(a\).
The discarded covariance has entries at most \(n^{-a}\); its transport
displacement is at most \(C_k\sqrt{rn}\,n^{-a}\) in \(L^r\).

The backward construction and weighted H\"older now bound the
\(L^r\) norm of the total centered displacement by
\[
 C_{k,K}\{rb\log(en)
       +n^{-1/2}T\sqrt{r\log(en)}\}
       +n^{3-a},
\]
because
\(\sum_i h_i\le C_k\log(en)\) and
\(\sum_i\sqrt{h_i}\le\sqrt{3n\sum_i h_i}\).
Removing the regularizations is absorbed in the first term.
The difference between the truncated and original tilted drifts is
at most \(Cn^{1-a}\).
Take \(r=\lceil(A+3)\log(en)\rceil\), choose \(a\) large enough,
apply Markov's inequality, and add
\eqref{thr:eq:conditional-drift} and the truncation probability.
This proves \eqref{thr:eq:conditional-tilted-coupling}.
\end{proof}

\section{Gaussian correlation estimates}\label{thr:sec:geometry}

The summability assumption on the correlations has two consequences
used below. Coordinates with correlations close to one in absolute
value can be grouped into blocks of bounded size. The same assumption
also bounds sums of Gaussian densities over rectangle faces. We prove
both facts in this section.
Throughout this section, \(\Sigma\) is a correlation matrix and
\begin{equation*}\label{thr:eq:row-power}
 \max_j\sum_i|\Sigma_{ij}|^s\le M
 \tag{G}
\end{equation*}
for fixed \(s>0\) and \(M\ge1\). All assertions allow singular
\(\Sigma\). A signed coordinate means \(\varepsilon_jG_j\), where
\(\varepsilon_j\in\{-1,1\}\) and \(G\sim N(0,\Sigma)\).
Changing coordinate signs preserves \eqref{thr:eq:row-power}.
We write \(\bar\Phi=1-\Phi\), where \(\Phi\) is the standard normal
distribution function.

\subsection{Correlation blocks and Gaussian rectangle probabilities}

\begin{lemma}\label{thr:block-partition}
There are constants \(b_0<\infty\) and \(\kappa>0\), depending only
on \(s,M\), and a partition of the coordinates into blocks of size
at most \(b_0\), such that
\begin{equation}
 r_{BC}:=\max_{i\in B,j\in C}|\Sigma_{ij}|\le1-\kappa
 \quad(B\ne C),\qquad
 \sup_B\sum_{C\ne B}r_{BC}^s\le b_0M .
 \label{thr:eq:block-correlation}
\end{equation}
\end{lemma}

\begin{proof}
Represent \(\Sigma\) as the Gram matrix of unit vectors \(v_j\).
Choose an integer \(b_0>2^sM+1\) and
\(\kappa<(4b_0^2)^{-1}\). Join \(i,j\) when
\(|\langle v_i,v_j\rangle|>1-\kappa\).
Along a path of length \(h\), choose successive signs so that the
Euclidean distance between its signed endpoint vectors is at most
\(h\sqrt{2\kappa}\). If a component had more than \(b_0\) vertices,
it would contain \(b_0\) distinct vertices within distance \(b_0\)
of one vertex. Their absolute correlations with that vertex would
be at least \(1-b_0^2\kappa>1/2\), contradicting
\eqref{thr:eq:row-power}. The components give the partition.
Finally,
\[
 \sum_{C\ne B}r_{BC}^s
 \le\sum_{i\in B}\sum_{j\notin B}|\Sigma_{ij}|^s\le b_0M .
 \qedhere
\]
\end{proof}

For a rectangle, an exit means that at least one coordinate violates
its interval constraint. Let \(q_B\) be the probability of such an
exit among the Gaussian coordinates in block \(B\). Omit unconstrained
blocks. Writing each lower exit as an upper exit of the negated
coordinate, let \(b_B\) be the smallest threshold in the block. Then
\begin{equation}\label{thr:eq:block-tail}
 \bar\Phi(b_B)\le q_B\le2b_0\bar\Phi(b_B).
\end{equation}

\begin{lemma}\label{thr:gaussian-cylinder}
For the blocks in Lemma~\ref{thr:block-partition},
\begin{equation}
 \Pp(G\in A)\le
 C_{s,M}\left(1+\sum_Bq_B\right)^{-1/(s+1)} .
 \label{thr:eq:qualitative-cylinder}
\end{equation}
In particular, if \(\sum_Bq_B>a+1\), the rectangle is contained
in a cylinder, obtained by retaining entire blocks of constraints,
for which the sum of the block exit probabilities belongs to
\([a,a+1]\).
\end{lemma}

\begin{proof}
Choose one exit at threshold \(b_B\) in each block. Write \(u_i\)
for these signed thresholds, \(a_i=\bar\Phi(u_i)\), and
\(w_i=1+(u_i)_+\). Let \(N\) count the selected Gaussian exits.
By \eqref{thr:eq:block-tail},
\(\E N=\sum_i a_i\) is comparable to \(\sum_Bq_B\).

Fix \(0<\epsilon<1/4\). If a pair has signed correlation \(r\le0\),
its indicator covariance is nonpositive. If
\(0<r\le\epsilon/(w_iw_j)\), interpolation of the bivariate normal
distribution gives
\[
 \Cov(\ind_{\{G_i>u_i\}},\ind_{\{G_j>u_j\}})
 \le C\epsilon a_i a_j .
\]
Indeed, throughout the interpolation its boundary density is at most
\(C\phi((u_i)_+)\phi((u_j)_+)\). For nonnegative thresholds
this follows from \(r u_i u_j\le\epsilon\); if one threshold is
negative, the cross term has the favorable sign, and if both are
negative the density is bounded. Integration over the correlation
interval and the normal hazard bound prove the display. These are
the bivariate interpolation identities underlying
\citet[Theorem~2.1]{LiShao2002}.

For the remaining positive pairs, order the thresholds so that
\(u_i\le u_j\), and put \(x=(u_i)_+\), \(y=(u_j)_+\).
The correlation is at most \(1-\kappa\).
If \(u_i>0\), exponential Markov inequality with vector
\((x,y)/(2-\kappa)\) gives
\[
 \Pp(G_i>u_i,G_j>u_j)
 \le \exp\left\{-\frac{x^2+y^2}{2(2-\kappa)}\right\}.
\]
The exponent gains at least \(c_\kappa y^2\) over \(x^2/2\).
If \(u_i\le0\), discard the first constraint and use the marginal
bound for the second. In either case, the condition
\(r w_iw_j>\epsilon\) and Mills' inequalities imply
\[
 \Pp(G_i>u_i,G_j>u_j)
 \le C_{s,\kappa}\epsilon^{-s}r^s a_i .
\]
For example, when both thresholds are positive, the exponential gain
absorbs \((1+x)^{s+1}(1+y)^s\); the other cases are weaker.
Assign each pair to the coordinate with the smaller threshold and use
\eqref{thr:eq:block-correlation}. We obtain
\[
 \Var(N)\le C_{s,M}
 \{\E N+\epsilon(\E N)^2+\epsilon^{-s}\E N\}.
\]
Chebyshev's inequality, with
\(\epsilon=(\E N)^{-1/(s+1)}\) when the mean is sufficiently large,
bounds \(\Pp(N=0)\) as asserted. For bounded means increase the
constant. For the cylinder assertion, add blocks until
\(\sum_B q_B\) first reaches \(a\). Each summand is a probability
and is therefore at most one.
\end{proof}

\subsection{Weighted Gaussian comparisons in fixed dimension}

The next lemma bounds changes in a block event by a product of
one-dimensional Gaussian tail probabilities. It therefore remains useful
when those probabilities tend to zero. The covariance within a block
need not be invertible.

\begin{lemma}\label{thr:lem:weighted-gaussian-shifts}
Fix a total dimension \(k\). Partition the coordinates into one group
with arbitrary rectangle constraints and finitely many further blocks.
We call the first group the core.
Let \(V_0\) be block diagonal for this partition, with unit diagonal.
In each further block \(B\), impose a nonempty collection of upper exits
whose thresholds belong to \([b_B,2b_B]\), where \(b_B\ge1\). We call
these the rare blocks. Let \(E\) be the intersection of all these exits
and the rectangle constraints on the first group, and put
\[
 x_B=\max\{\text{selected thresholds in }B\},\qquad
 g=\prod_B\bar\Phi(x_B).
\]
Let \(D_w\) have diagonal entries \(1+b_B\) in rare block \(B\)
and one in the first group. For \(0\le c\le1\), set
\(V_c=V_0+cD_w^{-2}\).
Write \(E^{+z}\) and \(E^{-z}\) for expansion and contraction of
each constraint by \(z\) divided by its coordinate weight. Then
\begin{align}
 \Pp(G_{V_c}\in E^{+z})&\le C_k g e^{C_kz},
       &&z\ge0,\label{thr:eq:gaussian-expansion}\\
 \Pp(G_{V_c}\in E^{+z}\setminus E^{-z})
     &\le C_k g z e^{C_kz},
       &&z\ge0,\label{thr:eq:gaussian-weighted-strip}\\
 |\Pp(G_{V_c}\in E)-\Pp(G_{V_0}\in E)|
     &\le C_k\sqrt c\,g .
       \label{thr:eq:weighted-ridge}
\end{align}
The bounds also hold when there are no rare blocks, with \(g=1\).
\end{lemma}

\begin{proof}
The Gaussian groups remain independent under \(V_c\).
Within a rare block all marginal variances equal
\(\sigma_B^2=1+c/(1+b_B)^2\). Choose a coordinate whose threshold
is \(x_B\) and tilt that coordinate with parameter
\(x_B/\sigma_B^2\). Leave the first group untilted.
The likelihood prefactor in block \(B\) is
\(\exp\{-x_B^2/(2\sigma_B^2)\}\). If only the constraint on the
selected coordinate is retained, the remaining weighted integral is
at most
\[
 \frac{C e^{Cz}}{1+x_B}
\]
after a threshold movement of size \(z/(1+b_B)\).
This follows by integrating the one-dimensional Gaussian density
against the exponential weight. Since \(x_B/(1+b_B)\le2\),
the effect of the movement is at most \(e^{Cz}\).
Moreover,
\[
 \frac{\exp\{-x_B^2/(2\sigma_B^2)\}}{1+x_B}
 \le C\bar\Phi(x_B),
\]
uniformly in \(c\le1\).
Discarding the other constraints proves
\eqref{thr:eq:gaussian-expansion}.

A change of a constraint in a rare block requires a coordinate
to lie in a slab of width \(2z/(1+b_B)\). Under the tilted law,
every coordinate has a density bounded by an absolute constant.
This uses only the marginal variance \(\sigma_B^2\ge1\), so
no joint density of the possibly singular block is needed.
On the enlarged event for the selected coordinate, the likelihood
weight is at most \(e^{Cz}\). Thus this block contributes at most
\(C_k z e^{C_kz}/(1+b_B)\), whereas every unchanged rare block
retains a one-dimensional integral of order \((1+x_B)^{-1}\).
The ratio \((1+x_B)/(1+b_B)\) is at most two.
Slabs in the first group have probability at most \(C_kz\).
Summing over the finitely many faces proves
\eqref{thr:eq:gaussian-weighted-strip}. This argument also covers
empty contractions.

Finally realize
\(G_{V_c}=G_{V_0}+\sqrt c\,D_w^{-1}Z\), with an independent
standard Gaussian \(Z\). The membership disagreement is contained
in the weighted boundary strip of width
\(\sqrt c\,\|Z\|_\infty\). Integrating
\eqref{thr:eq:gaussian-weighted-strip}, with \(c=0\), proves
\eqref{thr:eq:weighted-ridge}.
\end{proof}

\begin{lemma}\label{thr:lem:strong-gaussian}
Fix \(k\) and \(\kappa>0\). Let \(V\) be a \(k\)-dimensional
correlation matrix with \(|V_{ij}|\le1-\kappa\) for \(i\ne j\).
Let \(x_i\ge1\), and let a graph on these coordinates satisfy
\[
 |V_{ij}|x_ix_j\le\epsilon
 \quad\text{whenever \(i,j\) belong to different components.}
\]
Put \(D_h=1+(h-1)(1-\kappa)\).
If \(0\le\nu\max_i x_i^2\le1\) and \(\epsilon\le c_k\), then
\begin{equation}\label{thr:eq:strong-gaussian}
 \Pp(G_{V+\nu I}\ge x)
 \le C_{k,\kappa}
 \prod_{i\ {\rm isolated}}\bar\Phi(x_i)
 \prod_{\substack{C\ {\rm component}\\|C|\ge2}}
 \exp\left\{-\frac{\sum_{i\in C}x_i^2}{2D_{|C|}}\right\}.
\end{equation}
For every \(t\ge0\), the same bound for thresholds \(x-t\mathbf1\)
holds after multiplication by
\[
 \exp\left\{t\sum_i\frac{x_i}{D_{|C(i)|}}\right\},
\]
where \(C(i)\) is the component containing \(i\).
In particular, the extra factor is bounded when \(t\max_i x_i\le1\),
and is at most \(e^{kz}\) for \(t=z/T\), \(T\ge\max_i x_i\).
\end{lemma}

\begin{proof}
For \(i\) in a component of size \(h\), take
\(\lambda_i=x_i/D_h\). Since the absolute row sums of \(V_{CC}\)
are at most \(D_h\), its largest eigenvalue is at most \(D_h\).
Consequently
\[
 -\lambda^\top x+\tfrac12\lambda^\top(V+\nu I)\lambda
 \le-\sum_C\frac{\sum_{i\in C}x_i^2}{2D_{|C|}}+C_k .
\]
The error consists of covariances between different components and
the added diagonal covariance \(\nu I\).
Under the Gaussian exponential tilt the inverse likelihood, after
extracting this prefactor, is
\(\exp\{-\lambda^\top(G-x)\}\ind_{\{G\ge x\}}\).
Discard its nonisolated factors, which are at most one.
The covariance of the isolated coordinates has smallest eigenvalue
at least \(1/2\), by the cross-component assumption and \(x_i\ge1\).
Its density is therefore bounded by \(C_k\), regardless of its mean.
Integration over the isolated coordinates gives the factor
\(C_k\prod_{i\ {\rm isolated}}x_i^{-1}\).
Mills' inequalities give \eqref{thr:eq:strong-gaussian}.
A downward shift contributes exactly the additional upper factor
\(\exp(t\sum_i\lambda_i)\). The isolated density is still bounded
independently of its mean, and integration from the shifted boundary
still gives the factor \(\prod_{i\ {\rm isolated}}\lambda_i^{-1}\).
This proves all the shifted versions.
\end{proof}

\subsection{A finite forest estimate}

\begin{lemma}\label{thr:lem:face-sum}
Let \(J\) have at most three elements, let \(B\) be a correlation
matrix on \(J\), and let \(\|u\|_\infty\le T\). For every \(s>0\),
\begin{equation}\label{thr:eq:forest}
 \inf_{\|\theta\|_\infty\le T}
 e^{-\theta^\top u+\theta^\top B\theta/2}
 \le C_s\sum_{F\ {\rm forest\ on}\ J}
 \left(\prod_{ij\in F}|B_{ij}|^s\right)
 \prod_{C\ {\rm component\ of}\ F}
       \max_{j\in C}e^{-u_j^2/2}.
\end{equation}
Isolated vertices are included as components.
For a rectangle whose finite signed thresholds satisfy \(|b|\le T\),
put
\[
 W(A)=\sum_{\text{finite faces }b}e^{-b^2/2}.
\]
Under \eqref{thr:eq:row-power}, summing the left side of
\eqref{thr:eq:forest} over all choices of \(j\) distinct coordinate
faces, \(1\le j\le3\), gives at most
\begin{equation}\label{thr:eq:face-sum}
 C_{s,M}\sum_{\ell=1}^j W(A)^\ell .
\end{equation}
The same conclusion holds for covariance \(\Sigma+\nu I\) when
\(\nu T^2\le C_0\), with a constant depending also on \(C_0\).
\end{lemma}

\begin{proof}
The assertion for one vertex follows by taking \(\theta=u\).
If \(|B_{ij}|>1/100\), set \(\theta_i=0\), apply induction to the
remaining vertices, and attach \(i\) as a leaf of \(j\).
The constant increases by at most \(100^s\), since
\(|B_{ij}|^s>100^{-s}\), and the maximum weight in the enlarged
component cannot decrease.

Otherwise every off-diagonal entry has modulus at most \(1/100\).
For \(0\le r\le1/100\),
\[
 e^{r|xy|}\le C_s\{1+r^s e^{(x^2+y^2)/16}\}.
\]
When \(r|xy|\le1\) this is immediate. On its complement,
using \(x^2+y^2\ge2|xy|\), the ratio of the left side to
the second term is at most
\(r^{-s}\exp\{1-1/(8r)\}\), uniformly bounded in \(r\).
Take \(\theta=u\), multiply these inequalities over the at most
three pairs, and expand. Each isolated vertex retains
\(e^{-u_i^2/2}\). Every nonisolated vertex has degree at most two,
so its remaining factor is at most \(e^{-3u_i^2/8}\).
On a component with at least two vertices their product is no
larger than \(\max_i e^{-u_i^2/2}\).
If the selected edges contain a cycle, remove an edge; this only
increases the edge product and preserves the components.
This proves \eqref{thr:eq:forest}.

In each component choose a maximizing vertex as root, or bound its
maximum by the sum over possible roots. Sum from leaves toward roots.
Summation over each leaf gives a factor at most \(2M\), accounting for
the two endpoint choices of a coordinate. Summation over each root
gives \(W(A)\).
There are finitely many forests on at most three vertices, proving
\eqref{thr:eq:face-sum}.
For covariance \(\Sigma+\nu I\), divide the covariance by
\(1+\nu\) and the thresholds
by \(\sqrt{1+\nu}\). Correlations decrease in modulus, and
\[
 e^{-b^2/(2(1+\nu))}
 \le e^{C_0/2}e^{-b^2/2}.
\]
After returning to the original coordinates, the tilting parameters
still have absolute value at most \(T\), so the same proof applies.
\end{proof}

\subsection{Reduction to rectangles with controlled boundary sums}

\begin{proposition}\label{thr:prop:quantitative-cylinder}
Let \(n\ge2\), and assume
\(\log(2end)\le C_0n\).
Consider a rectangle all of whose finite signed thresholds satisfy
\(|b|\le D_0\sqrt{\log(2end)}\). In each constrained coordinate
choose the smaller signed threshold \(b_j\), and put
\[
 V(A)=\sum_j\bar\Phi(b_j).
\]
Set \(u=\max\{2,s\}\). Then either
\[
 V(A)\le C\{\log(en)\}^{u},
\]
or \(A\) is contained in a cylinder \(A_0\), obtained by retaining
selected faces, with
\[
 V(A_0)\le C\{\log(en)\}^{u},
 \qquad \Pp(G\in A_0)\le Cn^{-8}.
\]
For the original rectangle in the first case, or for \(A_0\) in the
second case, the sum of the face weights satisfies
\begin{equation}\label{thr:eq:quantitative-weight}
 W\le C\{\log(en)\}^{u}\sqrt{\log(ed)}.
\end{equation}
Constants depend only on \(s,M,C_0,D_0\).
\end{proposition}

\begin{proof}
For high thresholds we compare with a Gaussian vector having a common
factor. For low thresholds we first obtain a second-moment estimate
within each threshold band and then apply Gaussian isoperimetry.

\paragraph{High thresholds.}
Consider first the selected thresholds
\[
 b_j>A_0\sqrt{\log(en)},
\]
where \(A_0\) is a sufficiently large constant. We show that either
their tail probabilities sum to at most \(C\log(en)\), or a subset
of these constraints defines a cylinder of Gaussian probability at
most \(Cn^{-8}\).
Join two coordinates when their signed correlation exceeds \(1/2\).
The graph has degree at most \(2^sM\), so it can be colored with a
bounded number of colors.
If the sum of the tail probabilities is large, choose a color class
with the largest such sum and retain coordinates until the sum lies
between \(C_1\log(en)\) and \(C_1\log(en)+1\).
Within this collection all signed correlations are at most \(1/2\).
Set
\[
 \theta=\{8\log(2end)\log(en)\}^{-1},
\]
and compare its covariance \(R\) to the equicorrelation matrix
\(R_\theta\).

The Gaussian interpolation formula used in
\citet[Theorem~2.1]{LiShao2002} gives
\[
 \Pp_R(G_j\le b_j\ \forall j)
 \le \Pp_{R_\theta}(G_j\le b_j\ \forall j)
 +C\sum_{i<j:R_{ij}>\theta}
                  e^{-(b_i^2+b_j^2)/3}.
\]
Indeed, along \(tR+(1-t)R_\theta\), the derivative for one
off-diagonal entry is its bivariate boundary density times a
conditional probability between zero and one. Entries at most
\(\theta\) contribute nonpositively. The other correlations lie
between \(\theta\) and \(1/2\), giving the displayed density bound.
The interpolation is positive definite for \(t<1\); the singular
endpoint follows by continuity.

Assign each pair to the coordinate with the smaller threshold. Each
coordinate is assigned at most \(M\theta^{-s}\) neighbors, and
\[
 e^{-(b_i^2+b_j^2)/3}
 \le C(1+b_i)e^{-b_i^2/6}\bar\Phi(b_i)
 \quad\text{when }b_i\le b_j.
\]
Since \(\log(2end)\le C_0n\), the factor \(\theta^{-s}\) grows
at most polynomially in \(n\). Choosing \(A_0\) large enough
makes the total interpolation error at most \(Cn^{-10}\).

Represent the comparison vector as
\(\sqrt\theta\,Z_0+\sqrt{1-\theta}\,Z_j\), with independent
standard normal variables. Except on \(Z_0<-\sqrt{24\log(en)}\),
the conditional threshold exceeds \(b_j\) by at most
\(C_{D_0}/\sqrt{\log(2end)}\).
The Gaussian hazard bound therefore shows that each conditional exit
probability is at least a fixed multiple of \(\bar\Phi(b_j)\).
Conditional independence bounds the no-exit probability by
\[
 Cn^{-10}+\exp\{-cC_1\log(en)\}\le Cn^{-8}
\]
for sufficiently large \(C_1\).

\paragraph{Low thresholds.}
We next prove the exponent \(u=\max\{2,s\}\) for every \(s>0\).
The row sums of \(|\Sigma_{ij}|^u\) are at most \(M\).
Separate the central group \(b_j\le1\), indexed below by \(x=0\),
and dyadic groups
\[
 x<b_j\le2x,\qquad x=1,2,4,\ldots,
\]
up to \(A_0\sqrt{\log(en)}\).
Within each group, join coordinates whose positive signed correlation
exceeds \(1/2\), and color the resulting graph. At most
\(1+\lceil2^uM\rceil\) colors are needed.
Define
\[
 t_x=(2+2x)^{-1},\qquad
 f_x=(1+x)e^{-x^2/6},
\]
and the cutoff for the sum of tail probabilities
\[
 K_x=C_2\left\{\log(en)(1+x)^2+
       \{\log(en)\}^u(1+x)^{2u}f_x\right\}.
\]
If the sum of the tail probabilities in a group exceeds
\((1+\lceil2^uM\rceil)(K_x+1)\), choose a color class with the
largest such sum and retain a collection with
\(v=\sum_jq_j\in[K_x,K_x+1]\), where \(q_j=\bar\Phi(b_j)\).
Let
\[
 I=\{j:G_j>b_j+t_x\},\qquad Y=|I|,
\]
in signed coordinates.
The hazard bound gives \(cv\le\E Y\le v\).

Pairs with signed correlation at most \(c_0/(1+x)^2\) have joint
shifted exit probability at most \(Cq_iq_j\).
For negative correlations use comparison with independence.
For positive correlations, integrate the density bound
\(C\exp\{-(z_i^2+z_j^2)/(2(1+r))\}\);
the resulting tail ratios are bounded because
\(r(1+2x)^2\) is bounded.
In the central group a negative threshold has exit probability
at least one half, which gives the same conclusion.

Every pair in the chosen color also satisfies
\[
 \Pp(i,j\in I)\le C f_x\bar\Phi(\min(b_i,b_j)).
\]
For a positive threshold band, both exits force the sum of the two
signed Gaussian coordinates
above twice the smaller threshold, and the variance of that sum is
at most three. Hence its probability is at most \(e^{-2b^2/3}\),
which has the stated bound by Mills' inequalities.
In the central group use a marginal bound.
Assign each pair to its smaller threshold. The number of neighbors
whose signed correlation exceeds \(c_0/(1+x)^2\) is at most
\(CM(1+x)^{2u}\).
Since \((1+x)^{2u}f_x\) is bounded, these estimates imply
\[
 \E Y^2\le C(v^2+v),\qquad
 \E\sum_{\substack{i,j\in I\\i\ne j}}|\Sigma_{ij}|^u
       \le CM f_xv .
\]
Paley--Zygmund, followed by Markov's inequality with sufficiently
large constants, yields an event of probability at least \(c>0\)
on which
\[
 cv\le Y\le Cv,\qquad
 \sum_{\substack{i,j\in I\\i\ne j}}|\Sigma_{ij}|^u\le C f_xv.
\]

Represent the signed Gaussian coordinates as
\(\langle v_j,Z\rangle\) on the
Gaussian support. For any point of the no-exit set, a displacement
from \(Z\) must lower each coordinate indexed by \(I\) by at least
\(t_x\).
Taking its inner product with \(\sum_{j\in I}v_j\), and then
using H\"older's inequality, gives on the preceding event
\begin{align*}
 \operatorname{dist}(Z,\text{no-exit set})^2
 &\ge
 \frac{t_x^2Y^2}
 {Y+\sum_{i\ne j\in I}|\Sigma_{ij}|}\\
 &\ge c t_x^2\min\{v,v^{1/u}f_x^{-1/u}\}.
\end{align*}
For the second inequality use
\[
 \sum_{i\ne j\in I}|\Sigma_{ij}|
 \le Y^{2-2/u}
       \left(\sum_{i\ne j\in I}|\Sigma_{ij}|^u\right)^{1/u}.
\]
An empty no-exit set needs no estimate.
The definition of \(K_x\), with \(C_2\) large enough, makes the
distance exceed \(5\sqrt{\log(en)}+C\).
By the Gaussian isoperimetric inequality \citep{Borell1975},
a set whose neighborhood of this radius misses a set of fixed
positive Gaussian probability has Gaussian probability at most
\(\Phi(-5\sqrt{\log(en)})\le Cn^{-8}\).

The dyadic sum of \(\log(en)(1+x)^2\) is \(O(\log(en)^2)\).
The sum of \((1+x)^{2u}f_x\) over the dyadic groups is finite.
Thus the bounds for all groups, including the high thresholds, sum
to \(C\{\log(en)\}^u\). If a group exceeds its bound, the selected
constraints in that group give the required cylinder. Otherwise
\(V(A)\le C\{\log(en)\}^u\).

\paragraph{Sum of the face weights.}
Faces with \(b>2\sqrt{\log(ed)}\) contribute at most
\(2d e^{-2\log(ed)}\le1\) to \(W\).
Every other face satisfies
\[
 e^{-b^2/2}\le C\sqrt{\log(ed)}\,\bar\Phi(b),
\]
including negative \(b\). In each coordinate the smaller signed
threshold has the larger tail probability, so the sum over both faces
is at most twice \(V\). This proves
\eqref{thr:eq:quantitative-weight}.
\end{proof}

\section{Quantitative bounds with random quadratic variation}
\label{quant:section}

We prove the light-tail upper bounds in Theorem~\ref{main:joint}.
We first assume an almost-sure bound on \(\|Q-\Sigma\|_{\max}\),
then use Lemma~\ref{thr:lem:feasible-completion} to obtain a bound
in terms of its \(L^p\) norm. The main obstacle is that
\(\Sigma-P_i\) need not be positive semidefinite and hence cannot
serve as the covariance of an added Gaussian. We decompose it into
a positive Gaussian covariance and a symmetric error matrix, and
apply a finite differential correction to the smoothed test function.
The correction cancels changes of the decomposition through the
chosen order. For each fixed coordinate face, we bound the resulting
density uniformly over predictable indices and then sum the
conditional covariance weights.

For a positive semidefinite matrix \(A\), write
\[
 P_Af(x)=\E f(x+G_A),\qquad G_A\sim N(0,A).
\]
These operators include \(P_0f=f\). For a symmetric matrix \(E\), put
\(D_E=\sum_{j,k}E_{jk}\partial_{jk}\).

\subsection{Product smoothing and a finite covariance correction}

Let
\[
 p(z)=\frac1{4\cosh^2(z/2)},\qquad p_h(z)=h^{-1}p(z/h).
\]
For an interval \(I_j\), possibly unbounded, set
\(f_j(x)=\int_{I_j}p_h(y-x)\,dy\). For a rectangle
\(A=\prod_{j=1}^d I_j\), let \(f=\prod_j f_j\). Thus \(f\) is the
probability of membership after adding independent logistic variables of
scale \(h\). Positive smoothing and weighted Gaussian boundary
estimates also enter \citet[Appendix~B]{CCKK2026} and
\citet[Sections~3--4 and Appendix~H]{GaoZhang2026}. Those works
use a multinomial logistic construction. Here the coordinate-product
kernel gives the following relative derivative bounds, including
their dependence on the derivative order.

\begin{lemma}\label{quant:lem:product-derivatives}
For every integer \(\ell\ge0\),
\begin{equation}\label{quant:eq:product-global}
 \sum_{\alpha\in[d]^\ell}|\partial_\alpha f(x)|
 \le \ell!(C/h)^\ell.
\end{equation}
If \(\ell\ge1\), the left side is also at most
\begin{equation}\label{quant:eq:product-local}
 \ell!(C/h)^\ell
 \sum_{j=1}^d\min\{f_j(x_j),1-f_j(x_j)\}.
\end{equation}
Both estimates pass through Gaussian convolution, with the right
side of \eqref{quant:eq:product-local} convolved by the same Gaussian.
The constant \(C\) is absolute.
\end{lemma}

\begin{proof}
The function \(p\) is analytic on the strip \(|\operatorname{Im}z|<\pi\).
For real \(x\) and complex \(z\) with \(|z|\le1\),
\(|p(x+z)|\le Cp(x)\). Indeed,
\[
|\cosh((x+z)/2)|
 \ge \cos(1/2)\cosh((x+\operatorname{Re}z)/2)
 \ge e^{-1/2}\cos(1/2)\cosh(x/2)
\]
by the elementary identity for the modulus of a complex hyperbolic
cosine.
Cauchy's formula therefore gives
\(|p^{(\ell)}(x)|\le \ell!C^{\ell+1}p(x)\).
Integrating either over \(I_j\) or over its complement yields, for
\(\ell\ge1\),
\[
 |f_j^{(\ell)}(x)|
 \le \ell!(C/h)^\ell\min(f_j(x),1-f_j(x))
 \le 2\ell!(C/h)^\ell f_j(x)(1-f_j(x)).
\]
In the sum of the ordered derivatives, a support \(J\) of size \(k\)
and positive multiplicities summing to \(\ell\) have multinomial
coefficient \(\ell!/\prod_{j\in J}\ell_j!\). The factorials in the
one-dimensional bounds cancel these denominators. If
\(v=\sum_j(1-f_j(x_j))\), summation over supports gives
\[
 \Bigl(\prod_jf_j(x_j)\Bigr)
 \sum_{|J|=k}\prod_{j\in J}(1-f_j(x_j))
 \le e^{-v}\frac{v^k}{k!}\le1.
\]
There are \(\binom{\ell-1}{k-1}\) choices of the positive
multiplicities. Summing over \(k\) proves
\eqref{quant:eq:product-global}.

For the localized estimate, mark a differentiated coordinate \(j\),
let \(a\in\{1,\ldots,\ell\}\) be its multiplicity, and apply
\eqref{quant:eq:product-global} to the product over the other
coordinates. This overcounts each ordered derivative at least once
and bounds the resulting sum by
\[
 \sum_j\sum_{a=1}^{\ell}
 \binom{\ell}{a}
 a!(C/h)^a\min(f_j,1-f_j)
 (\ell-a)!(C/h)^{\ell-a}.
\]
The factor \(\ell\) is absorbed into \(C^\ell\). Differentiation under
a Gaussian integral is justified by the bounded derivatives, and
gives the final assertion.
\end{proof}

\begin{lemma}\label{quant:lem:finite-correction}
Let \(R\ge1\) be an integer and define
\begin{equation}\label{quant:eq:corrected-function}
 \mathcal U_R(U,E;x)
 =\sum_{r=0}^R\frac{D_E^rP_Uf(x)}{2^rr!},
 \qquad U\succeq0,\quad E=E^\top.
\end{equation}
There is an absolute \(c>0\) such that the following assertions hold
whenever \(R\eta/h^2\le c\) and \(\|E\|_{\max}\le\eta\).
For \(0\le j\le4\),
\begin{equation}\label{quant:eq:corrected-derivatives}
 \sum_{\alpha\in[d]^j}
 |\partial_\alpha\mathcal U_R(U,E;x)|\le C_jh^{-j}.
\end{equation}
Moreover,
\begin{equation}\label{quant:eq:terminal-correction}
 |\mathcal U_R(U,E;x)-P_Uf(x)|
 \le C P_U\!\left[
       \sum_j\min(f_j,1-f_j)\right](x).
\end{equation}
If \(U,U'\succeq0\), \(\|E\|_{\max},\|E'\|_{\max}\le\eta\), and
\(U+E=U'+E'\), then
\begin{equation}\label{quant:eq:heat-switch}
 |\mathcal U_R(U,E;x)-\mathcal U_R(U',E';x)|
 \le (CR\eta/h^2)^{R+1}.
\end{equation}
\end{lemma}

\begin{proof}
The sum in \eqref{quant:eq:corrected-derivatives} is at most
\[
 h^{-j}\sum_{r=0}^R
 \frac{(2r+j)!}{2^rr!}C^{2r+j}(\eta/h^2)^r.
\]
For fixed \(j\le4\), its \(r\)-th summand is bounded by
\(C_j(1+r)^j(Cr\eta/h^2)^r\). Choosing \(c\) sufficiently small
makes the series summable uniformly in \(R\). Applying
\eqref{quant:eq:product-local} instead gives
\eqref{quant:eq:terminal-correction}.

For the last assertion, set \(U_t=U+t(U'-U)\) and
\(E_t=E-t(U'-U)\). Suppose first that the Gaussian bases are positive
definite. Price's covariance differentiation identity
\citep[Theorem~1 and equation~(1.8)]{Voigtlaender2021},
\(\frac{d}{dt}P_{U_t}f=\frac12D_{U'-U}P_{U_t}f\), gives
the exact cancellation
\begin{equation}\label{quant:eq:finite-cancellation}
 \frac{d}{dt}\mathcal U_R(U_t,E_t;x)
 =\frac{D_{U'-U}D_{E_t}^RP_{U_t}f(x)}{2^{R+1}R!}.
\end{equation}
Indeed, differentiating the coefficient \(D_{E_t}^r\) cancels the
Gaussian derivative of the preceding term. Since
\(\|U'-U\|_{\max}=\|E-E'\|_{\max}\le2\eta\), the absolute value of
the right side is at most
\[
 \frac{(2R+2)!}{2^RR!}(C/h)^{2R+2}\eta^{R+1}
 \le (CR\eta/h^2)^{R+1}.
\]
Integrate over \(t\). For semidefinite bases, add \(\varepsilon I\)
to \(U,U'\), apply the proved bound, and let \(\varepsilon\downarrow0\).
All derivatives of \(f\) appearing here are bounded and continuous,
so Gaussian weak convergence justifies this limit.
\end{proof}

\subsection{A localized third derivative}

The following estimate retains the covariance weight needed for
summation over the increments. We use the Gaussian face integration in
\citet[arXiv v1, Lemmas~2.2--2.3]{FKLZ2026}; we spell out the coefficients and
the additional logistic derivatives.

\begin{lemma}\label{quant:lem:localized-faces}
Let \(U\succ0\), \(F_U=P_Uf\), and \(v\in\R^d\). For each
\(J\subset[d]\), \(1\le |J|\le3\), let \(V_J\succeq U_{JJ}\) have
diagonal entries at least one. For every integer \(\ell\ge0\),
\begin{align}
 \sum_{\alpha\in[d]^\ell}|\partial_\alpha D_v^3F_U(x)|
 &\le \ell!(C/h)^\ell
 \sum_{\substack{J\subset[d]\\1\le |J|\le3}}\ \sum_{b\in\partial A_J}
 C\|v_J\|_\infty(v_J^\top U_{JJ}^{-1}v_J)\sqrt{\det V_J}
 \notag\\[-2pt]
 &\quad{}\times\int_{\R^{|J|}}
 (1+|U_{JJ}^{-1/2}y|)^{3-|J|}
 \phi_{U_{JJ}}(y)\prod_{j\in J}p(z_j)\,dz_J,
 \quad y=b-hz_J-x_J .
 \label{quant:eq:localized-faces}
\end{align}
Here \(\partial A_J\) chooses one finite endpoint from each interval
indexed by \(J\).
\end{lemma}

\begin{proof}
Fix the logistic shifts first. Differentiating a Gaussian rectangle
once in direction \(v\) sets coordinate \(j\) at an endpoint and
contributes the signed factor \(v_j\) times its marginal Gaussian density. Having fixed
coordinates \(I\), a new face \(l\) contributes
\[
 v_{l\mid I}=v_l-U_{lI}U_{II}^{-1}v_I.
\]
A derivative that stays on the current Gaussian density instead
contributes a Hermite factor. After three derivatives the final face
has \(k=1,2,\) or \(3\) fixed coordinates, and the Hermite degree is
at most \(3-k\). There are at most \(3!\) orders for a fixed final
face and a fixed endpoint choice.

To bound a term, whiten the Gaussian in these coordinates. Every
normalized normal or successive conormal coefficient is at most
\((v_J^\top U_{JJ}^{-1}v_J)^{1/2}\); equivalently,
\[
 \frac{v_{l\mid I}^2}
 {U_{ll}-U_{lI}U_{II}^{-1}U_{Il}}
 \le v_J^\top U_{JJ}^{-1}v_J.
\]
Keep the first coefficient as \(|v_j|/\sqrt{U_{jj}}\), and bound the
other two by their quadratic energy. Integration over the affine span
of the face contributes
\(\sqrt{\det U_{JJ}}\phi_{U_{JJ}}(b-x_J)\), multiplied by the conditional
probability that the remaining coordinates satisfy their intervals. The Schur determinant
is monotone under positive semidefinite domination, so
\[
 \frac{\det U_{JJ}}{U_{jj}}
 \le\frac{\det V_J}{(V_J)_{jj}}\le\det V_J.
\]
These observations give \eqref{quant:eq:localized-faces} for
\(\ell=0\), including singular limits of conditional faces after
positive Gaussian smoothing.

For higher derivatives, retain the signed face expression before
taking absolute values. Average over the logistic shifts on the
coordinates outside \(J\). The conditional probability for those
coordinates is a Gaussian expectation of their product of interval
functions. Derivatives in a coordinate outside \(J\) act on this
product. Derivatives in a coordinate in \(J\) can be transferred by
integration by parts to its logistic density: the dependence on
\(x_J,z_J\) is through \(b-hz_J-x_J\), including the conditional
Gaussian means. If \(a\) of the \(\ell\) derivatives act outside
\(J\), their ordered absolute sum is bounded by \(a!(C/h)^a\),
by Lemma~\ref{quant:lem:product-derivatives}. The other derivatives
contribute at most \((\ell-a)!(C/h)^{\ell-a}\), with at most
\(3^{\ell-a}\) choices of coordinates in \(J\). Summing the \(\binom{\ell}{a}\)
assignments proves the stated factorial bound. Gaussian decay and
the integrable logistic derivatives justify all integrations by
parts. Empty or unbounded intervals follow by limits.
\end{proof}

\subsection{A pathwise covariance bound}

\begin{proposition}\label{quant:prop:pathwise}
Let \(\Sigma\) be a deterministic correlation matrix with
\(\max_j\sum_k|\Sigma_{jk}|^s\le M\), and suppose
\[
 \|Q-\Sigma\|_{\max}\le\eta\quad\text{almost surely}.
\]
Here \(\Sigma\) need not equal \(\E Q\). Define the increment scale
\[
 a_n=
 \begin{cases}
 B/\sqrt n,&\|X_i\|_\infty\le B,\\[2pt]
 C_K\{\log(en)\}^{1/\alpha}/\sqrt n,
 &\E_{i-1}\exp((|X_{ij}|/K)^\alpha)\le2.
 \end{cases}
\]
The constant \(C_K\) in this definition is sufficiently large.
Set \(u=\max\{2,s\}\).
Then
\begin{align}
 \rho_{\cR}(S_n,G_\Sigma)
 \le C\min\bigl\{1,\,
 &a_n\{\log(en)\}^{1+3u}\{\log(ed)\}^{3/2}\notag\\
 &+\sqrt\eta\,\{\log(en)\}^{u+1/2}
                       \sqrt{\log(ed)}\bigr\}.
 \label{quant:eq:pathwise}
\end{align}
The same estimate holds if one conditionally Gaussian increment is
appended after the original increments, provided the stated covariance
bound holds for the resulting array.
\end{proposition}

\begin{proof}
We give the proof for the original increments first. The assertion about an
appended Gaussian will be checked below. Write
\(\xi_i=X_i/\sqrt n\), \(S_i=\sum_{\ell\le i}\xi_\ell\),
\(C_i=\E_{i-1}\xi_i\xi_i^\top\), and \(P_i=\sum_{\ell\le i}C_\ell\).
It suffices to treat the expression inside the minimum in
\eqref{quant:eq:pathwise} smaller than a fixed sufficiently small
constant. Choose
\begin{align}
 T&=D_0\sqrt{\log(2end)},\qquad R=\lceil100\log(en)\rceil,\notag\\
 g&=C_1\{a_n\log(2end)\log(en)+\sqrt\eta\},\notag\\
 h&=C_2\{a_n\log(2end)\log(en)+\sqrt{R\eta}\}.
 \label{quant:eq:widths}
\end{align}
The constants are chosen in this order, with \(C_2\) large relative
to \(C_1\). The assumed smallness implies
\[
 hT\le c,\qquad (g^2+\eta)T^2\le c,\qquad a_nT^3\le c,\qquad
 R\eta/h^2\le c,
\]
where the last constant is small enough for
Lemma~\ref{quant:lem:finite-correction}. Also \(g\ge n^{-1/2}\),
\(g\le h\), and \(\log(2end)\le Cn\). These implications follow
from \(\log(2end)\le C\log(en)\log(ed)\) and \(u\ge1\).

\paragraph{Positive covariance bases.}
For the product-smoothed test of a rectangle, use
\[
 U_i=\E_{i-1}Q-P_i+g^2I,\qquad E_i=\Sigma-\E_{i-1}Q.
\]
Since \(Q-P_i\succeq0\) and \(P_i\) is
\(\mathcal F_{i-1}\)-measurable, \(U_i\succeq g^2I\).
Conditional Jensen gives \(\|E_i\|_{\max}\le\eta\).
The pre-step pair is \((U_i+C_i,E_i)\), and the post-step pair is
\((U_i,E_i)\). When passing from step \(i\) to step \(i+1\), the two
decompositions have the same sum, \(\Sigma-P_i+g^2I\). Consequently,
\eqref{quant:eq:heat-switch} bounds the error from each change by
\(n^{-12}\), after increasing \(C_2\). The initial pair can be compared
with \((\Sigma+g^2I,0)\), and the final pair with
\((g^2I,\Sigma-Q)\), with the same error bound. The total error from
these changes is at most \(Cn^{-11}\).

Within step \(i\), condition on \(\mathcal F_{i-1}\), and let
\(\zeta_i\sim N(0,C_i)\). Taylor expansion to order two compares
\[
 \mathcal U_R(U_i,E_i;S_{i-1}+\xi_i)
 \quad\text{and}\quad
 \mathcal U_R(U_i,E_i;S_{i-1}+\zeta_i).
\]
The constant, linear and quadratic terms have identical conditional
expectations. The Gaussian expectation of the second expression is
\(\mathcal U_R(U_i+C_i,E_i;S_{i-1})\). Thus only integrated third
derivatives remain.

\paragraph{Summing the covariance weights.}
Fix a set \(J\) of \(k\le3\) coordinates and define
\[
 \Gamma=\Sigma_{JJ}+3\eta I_k,\quad
 V=\Gamma+g^2I_k,\quad B_i=V-P_{i,JJ},\quad
 w_i=\tr(B_i^{-1}C_{i,JJ}).
\]
The entrywise error on a matrix of size at most three implies
\[
 P_{i,JJ}\preceq Q_{JJ}\preceq\Gamma,\qquad
 g^2I_k\preceq U_{i,JJ}\preceq V,\qquad
 B_i\preceq U_{i,JJ}+6\eta I_k\preceq C U_{i,JJ}.
\]
In either increment model, \(\|C_{i,JJ}\|_{\mathrm{op}}\le Ca_n^2\).
The eigenvalues of \(B_i^{-1/2}C_{i,JJ}B_i^{-1/2}\) are therefore
bounded by an absolute constant smaller than one. Applying
\(x\le2\log(1+x)\) and telescoping determinants gives
\begin{equation}\label{quant:eq:clock-budget}
 \sum_i w_i
 \le C\log\frac{\det V}{\det B_n}
 \le C\log(en).
\end{equation}
The last bound uses the bounded diagonal of \(V\) and
\(B_n\succeq g^2I_k\succeq n^{-1}I_k\).

Apply Lemma~\ref{quant:lem:localized-faces} to the third derivatives
of the corrected function. Summing the terms of order \(2r\) with
the coefficients of \(E_i\) contributes a factor at most
\(\sum_{r\le R}(2r)!(C\eta/h^2)^r/r!\le C\).
In a face integral, discard the region
\(|U_{i,JJ}^{-1/2}y|>D_1\sqrt{\log(2end)}\).
Its polynomial times Gaussian density is at most
\(Cg^{-k}(2end)^{-cD_1^2}\). The expected remaining coefficient is
bounded by \(Cg^{-2}\E\|v_J\|_\infty^3\), where \(v\) is the real
or Gaussian increment. There are at most \(n(2d)^3\) terms,
\(\sqrt{\det V}\le C\), and \(g^{-1}\le\sqrt n\). Choosing \(D_1\)
large makes the total discarded contribution at most \(n^{-10}\).
On the retained region the Hermite factor is bounded by
\(C\{\log(2end)\}^{(3-k)/2}\).

Let
\[
 Z_i(b)=
 (\phi_{B_{i-1}}*p_h^{\otimes k})(b-S_{i-1,J}).
\]
The marginal density at a Taylor displacement \(tv_J\) is at most
\(C\exp(C\|v_J\|_\infty/h)Z_i(b)\). To see this, translate the
logistic variables and use \(|(\log p)'|\le1\).
We also compare the two Gaussian covariance matrices. The matrices
\(U_{i,JJ}\) and \(B_{i-1}\) differ entrywise by at most
\(C(\eta+a_n^2)\). For their convolved logistic density, each second
derivative is at most \(Ch^{-2}\) times the density. The heat
equation along the positive covariance segment therefore bounds the
logarithm of the density ratio by
\(C(\eta+a_n^2)/h^2\le C\).

The conditional coefficient multiplying \(Z_i(b)\) is bounded by
\begin{equation}\label{quant:eq:conditional-cubic}
 \E_{i-1}\!\left[
 \|v_J\|_\infty(v_J^\top U_{i,JJ}^{-1}v_J)
 e^{C\|v_J\|_\infty/h}\right]
 \le Ca_n(w_i+n^{-10}).
\end{equation}
For a bounded real increment, remove its supremum norm before
taking the conditional quadratic moment. For the Gaussian
increment, write \(v_J=C_{i,JJ}^{1/2}Z\), with standard \(Z\).
Then \(\|v_J\|_\infty\le Ca_n|Z|\), and
\[
 Z^\top C_{i,JJ}^{1/2}U_{i,JJ}^{-1}C_{i,JJ}^{1/2}Z
 \le Cw_i|Z|^2.
\]
The fixed-dimensional moment of \(|Z|^3e^{C|Z|}\) is finite.
For a real conditional-tail increment, split at
\(\|v_J\|_\infty=a_n\). The lower part has the same covariance
bound. Conditional tail integration, with the constant in \(a_n\)
sufficiently large, gives
\[
 \E_{i-1}\!\left[
 \|v_J\|_\infty^3 e^{C\|v_J\|_\infty/h}
 \ind_{\{\|v_J\|_\infty>a_n\}}\right]\le n^{-25}.
\]
Here \(h\ge Ca_n\), and the fixed-dimensional conditional tail is
at most \(6\exp[-(t\sqrt n/K)^\alpha]\).
Since \(U_{i,JJ}^{-1}\preceq g^{-2}I\), this proves
\eqref{quant:eq:conditional-cubic}. The small absolute tail is
multiplied by the same marginal density before summing over coordinates;
it is not multiplied by the number of faces at this stage.
The enlarged predictable weights \(w_i+n^{-10}\) still have total
at most \(C\log(en)\).

\paragraph{Densities at predictable indices.}
For a signed endpoint vector \(b\), \(\|b\|_\infty\le T\), put
\[
 \omega_J(b)=
 \inf_{\|\lambda\|_\infty\le T}
 \exp\{-\lambda^\top b+\lambda^\top\Sigma_{JJ}\lambda/2\}.
\]
We claim that, uniformly over predictable
\(\sigma\in\{1,\ldots,n+1\}\),
\begin{equation}\label{quant:eq:stopped-density}
 \sqrt{\det V}\,\E Z_\sigma(b)
 \le C\omega_J(b)+(2nd)^{-30}.
\end{equation}
At \(\sigma=n+1\), interpret the kernel using the full prefix
\(P_n,S_n\). Stop the selected martingale just before \(\sigma\),
and complete its covariance to \(\Gamma\). Denote the completed sum
by \(U\). The added Gaussian has conditional covariance
\(\Gamma-P_{\sigma-1,JJ}\).

For bounded increments this completion is the weak limit of bounded sign
completions: split each column of the residual square root into
arbitrarily many independent signed steps with coordinate bound
at most \(a_n\). The total conditional covariance remains \(\Gamma\), its trace
is bounded by a constant depending only on \(k\), and
Proposition~\ref{core:prop:transport} is uniform in the number of
steps. Passing to the Gaussian-completion limit by lower
semicontinuity gives
\[
 W_{r,\infty}(U,G_\Gamma)
 \le C_k r a_n\log(e+1/a_n).
\]
Choose \(r=C\log(2end)\). The width \(g\) in
\eqref{quant:eq:widths} permits a coupling of \(U\) and
\(Z\sim N(0,\Gamma)\) with
\(\|U-Z\|_\infty\le g\), except on an event of probability
at most \((2nd)^{-100}\). On the good event, completing a square
gives the pointwise bound
\[
 \phi_{g^2I_k}(y-U)
 \le 2^{k/2}e^{k/2}\phi_{2g^2I_k}(y-Z).
\]
Take expectations and use Gaussian convolution. On the complement,
the kernel is at most \(Cg^{-k}\). These two estimates give
\begin{equation}\label{quant:eq:kernel-coupling}
 \E\phi_{g^2I_k}(y-U)
 \le C\phi_{\Gamma+2g^2I_k}(y)+(2nd)^{-80}.
\end{equation}

For conditional tails, truncate each original projected increment at
supremum norm \(a_n\) and center conditionally. Keep all the
completion signs, whose bound can be chosen at most \(a_n\).
For every \(\|\lambda\|_\infty\le T\), conditional tail integration
yields
\begin{equation}\label{quant:eq:tilted-tail}
 \E_{i-1}\!\left[
 (1+\|x\|_\infty^4+|\lambda^\top x|^4)e^{2|\lambda^\top x|}
 \ind_{\{\|x\|_\infty>a_n\}}\right]\le n^{-40}.
\end{equation}
The lost covariance is positive semidefinite and has trace at most
\(Cn^{-39}\) after summation. Append independent signed columns
of its square root. The resulting bounded array has total conditional
covariance \(\Gamma\); its sum differs from \(U\) by at most \(Cn^{-18}\)
when no original increment was truncated. The bounded comparison,
lifted conditionally on this bounded sum, gives
\eqref{quant:eq:kernel-coupling} on that event with
\(\Gamma+4g^2I_k\) in place of \(\Gamma+2g^2I_k\).

For the discarded event use the predictable likelihood
\[
 L_N(\lambda)=
 \prod_i\frac{\exp(\lambda^\top x_i)}
                 {\E_{i-1}\exp(\lambda^\top x_i)}
\]
of the untruncated completed array. The number of original increments is
at most \(n\). Formula \eqref{quant:eq:tilted-tail}, together with
\(\E_{i-1}e^{\lambda^\top x_i}\ge1\), bounds its probability under
this tilt by \(Cn^{-39}\). Taylor expansion of each conditional
exponential gives the pathwise compensator bound
\[
 \sum_i\log\E_{i-1}e^{\lambda^\top x_i}
 \le \tfrac12\lambda^\top\Gamma\lambda+Ca_nT^3+Cn^{-39}.
\]
For the original increments this follows from their conditional
third moments and exponential tails. For the arbitrarily many
bounded completion signs the total cubic remainder is at most
\(Ca_n\|\lambda\|_1\lambda^\top\Gamma\lambda\).
The identity
\[
 e^{-\lambda^\top u}\phi_{g^2I_k}(y-u)
 =e^{-\lambda^\top y+g^2|\lambda|^2/2}
       \phi_{g^2I_k}(y-u-g^2\lambda)
\]
places the discarded event on the same exponential scale as the
Gaussian density. Under the tilted law its probability is at most
\(Cn^{-39}\). Use the pathwise compensator bound above and bound
the remaining kernel by \(Cg^{-k}\). Its contribution is therefore
at most
\[
 Cg^{-k}n^{-39}
 \exp\{-\lambda^\top y+
           \lambda^\top(\Gamma+g^2I_k)\lambda/2\}.
\]
This is absorbed into the same exponential bound as the main
Gaussian term, uniformly in the sign subdivision. All bounds pass
to the Gaussian-completion limit.

For any positive definite \(A\),
\[
 \phi_A(y)\le
 (2\pi)^{-k/2}(\det A)^{-1/2}
 \exp\{-\lambda^\top y+\lambda^\top A\lambda/2\}.
\]
Multiplication by \(\sqrt{\det V}\) cancels the normalizing factor
in the preceding bounds up to a constant. Integrate at \(y=b-hz\)
against the product logistic density. Its exponential moment is
bounded because \(hT\) is small; the diagonal covariance inflation
costs a bounded factor because \((\eta+g^2)T^2\) is small.
Optimizing over \(\lambda\) proves
\eqref{quant:eq:stopped-density}. This argument never takes an
inverse of the possibly singular target covariance.

Lemma~\ref{core:lem:clock}, \eqref{quant:eq:clock-budget} and
\eqref{quant:eq:conditional-cubic} now bound the sum of the
within-step errors by
\begin{equation}\label{quant:eq:drift}
 Ca_n\log(en)
 \sum_{k=1}^3\{\log(2end)\}^{(3-k)/2}
 \sum_{\substack{|J|=k\\b\in\partial A_J}}\omega_J(b)
 +O(n^{-6}).
\end{equation}

\paragraph{The terminal correction and unsmoothing.}
For every scalar endpoint \(|b|\le T\),
\begin{align}
 &\E e^{-|S_{n,j}+gG_j-b|/h}\notag\\
 &\quad+
 \Pp\{\text{the segment from }S_{n,j}
          \text{ to }S_{n,j}+gG_j+hZ_j\text{ contains }b\}
 \notag\\[-2pt]
 &\hspace{12mm}\le
 C(h+g+\sqrt\eta)e^{-b^2/2}+(2nd)^{-25},
 \label{quant:eq:scalar-strip}
\end{align}
where the normal and logistic variables are independent of the
martingale. To prove this, complete the scalar predictable quadratic variation
to \(1+\eta\) by a fresh conditional Gaussian and add an independent
normal of variance \(g^2\). Call the resulting sum \(U'\).
The total added conditional variance is at most \(2\eta+g^2\).
Conditionally on the original terminal data, this added normal has
absolute value at most \(2\sqrt{2\eta+g^2}\) with probability at
least \(1/2\). Hence
\[
 \Pp(|S_{n,j}-b|\le z)
 \le2\Pp(|U'-b|\le z+2\sqrt{2\eta+g^2}).
\]
The proof of \eqref{quant:eq:stopped-density}, before logistic
integration and with the fixed scalar tilt \(\lambda=b\), bounds
the density of \(U'\) at \(y\) by
\[
 C\exp\{-by+(1+\eta+4g^2)b^2/2\}+(2nd)^{-70}.
\]
It follows that a strip of radius \(z\) around \(b\) has
probability at most
\[
 C(z+g+\sqrt\eta)e^{Tz}e^{-b^2/2}
       +C(z+g+\sqrt\eta)(2nd)^{-70}.
\]
For the crossing event in \eqref{quant:eq:scalar-strip}, integrate
this estimate at \(z=g|G_j|+h|Z_j|\). For the exponential
expectation, use \(z=g|G_j|+hE\), with \(E\) an independent
mean-one exponential variable. The required exponential moments
are bounded by \(hT,gT\le c\), which proves
\eqref{quant:eq:scalar-strip}.

Let
\[
 W_A=\sum_{\text{finite faces }b}e^{-b^2/2}.
\]
If smoothing changes membership in the rectangle, at least one
coordinate crosses an endpoint. Also
\[
 \min(f_j(x),1-f_j(x))
 \le C\sum_{b\in\partial I_j}e^{-|x-b|/h}.
\]
Thus \eqref{quant:eq:terminal-correction} at
\((g^2I,\Sigma-Q;S_n)\), together with
\eqref{quant:eq:scalar-strip}, bounds both the terminal correction
and the final unsmoothing by
\begin{equation}\label{quant:eq:boundary-charge}
 C(h+g+\sqrt\eta)W_A+O(n^{-7}).
\end{equation}
The initial smoothing at \(G_\Sigma\) satisfies the same bound,
directly from its unit marginal variances.

We may restrict all finite endpoints to \([-T,T]\). For bounded
increments, the scalar Bernstein consequence of
\citet{Freedman1975}, and \(Q_{jj}\le1+\eta\), bound the error from this restriction
by \(O(n^{-8})\) after increasing \(D_0\).
For conditional tails, conditional exponential expansion at
\(\lambda=\pm T\) gives
\[
 \log\E e^{\lambda S_{n,j}}
 \le (1+\eta)\lambda^2/2+Ca_nT^3.
\]
Here the bound follows by multiplying the conditional exponential
inequalities; it does not replace the martingale by independent
increments. Exponential Markov and a union bound give the same clipping
error. The Gaussian tail bound gives the same estimate for the target.

\paragraph{Summing over the rectangle faces.}
Lemma~\ref{thr:lem:face-sum} bounds the inner sums in
\eqref{quant:eq:drift} by \(C\sum_{\ell=1}^k W_A^\ell\).
Proposition~\ref{thr:prop:quantitative-cylinder} allows us to
retain the original rectangle or a containing cylinder with
\[
 W_A\le C\{\log(en)\}^u\sqrt{\log(ed)}.
\]
In the cylinder alternative its Gaussian probability is at most
\(Cn^{-8}\). Since
\(\sqrt{\log(2end)}\le
C\{\log(en)\}^u\sqrt{\log(ed)}\), substitution in
\eqref{quant:eq:drift} gives the first term of
\eqref{quant:eq:pathwise}. Formula \eqref{quant:eq:boundary-charge}
gives the same sampling term and the displayed square-root
covariance term. A containing cylinder bounds both probabilities
of the original rectangle by its own Gaussian probability plus
the proved comparison error. This proves the assertion for every
rectangle. Atoms and open or closed endpoints are covered by
closed-strip estimates and monotone limits.

Finally suppose a conditional Gaussian increment is appended.
Its within-step comparison is an exact heat convolution and has
zero Taylor remainder. There is one additional switch of
covariance decomposition. In every stopped completion its
Gaussian covariance combines with the newly added Gaussian;
alternatively both can be subdivided into bounded sign steps
before applying Proposition~\ref{core:prop:transport}.
The exponential compensator of the appended increment is exactly
quadratic. The scalar-strip and clipping arguments are therefore
unchanged, and only the original increments contribute to
\eqref{quant:eq:drift}. This proves the final assertion.
\end{proof}

\subsection{Bounds under an \texorpdfstring{$L^p$}{Lp} condition on quadratic variation}

\begin{corollary}\label{quant:cor:joint-rates}
Under the assumptions on the deterministic target in
Proposition~\ref{quant:prop:pathwise}, suppose
\[
 \big\|\|Q-\Sigma\|_{\max}\big\|_{L^p}\le\delta,
 \qquad 1\le p\le\infty.
\]
Put \(\gamma_p=(2+1/p)^{-1}\), with \(1/\infty=0\), and set
\(u=\max\{2,s\}\). In the bounded and
conditional exponential-moment models,
\begin{align}
 \rho_{\cR}(S_n,G_\Sigma)
 \le C\min\bigl\{1,\,
 &a_n\{\log(en)\}^{1+3u}\{\log(ed)\}^{3/2}\notag\\
 &+\{\delta\{\log(en)\}^{1+2u}\log(ed)\}^{\gamma_p}\bigr\}.
 \label{quant:eq:joint-rates}
\end{align}
For marginal coordinate exponential moments, the same estimate
holds with
\[
 a_n=C_K\{\log(2end)\}^{1/\alpha}/\sqrt n .
\]
In particular, the common exponent of \(\log(en)\) can be taken
to be \(2+3\max\{2,s\}\).
\end{corollary}

\begin{proof}
For any \(\eta>0\), use the stopping and completion construction
of Lemma~\ref{thr:lem:feasible-completion}.
It retains the original increments until their cumulative
conditional covariance can no longer be dominated by a positive
semidefinite matrix within entrywise distance \(\eta\) of
\(\Sigma\). At the end, choose a measurable feasible cover and
append a conditional Gaussian with the missing covariance.
The resulting predictable quadratic variation is within \(\eta\) of \(\Sigma\)
on every path, and on
\(\{\|Q-\Sigma\|_{\max}\le\eta\}\) all original increments are retained
and the appended Gaussian is zero. The martingale sum therefore
changes with probability at most
\(\Pp(\|Q-\Sigma\|_{\max}>\eta)\).
The last assertion of Proposition~\ref{quant:prop:pathwise} gives
\begin{align}
 \rho_{\cR}(S_n,G_\Sigma)
 \le{}&
 Ca_n\{\log(en)\}^{1+3u}\{\log(ed)\}^{3/2}\notag\\
 &+C\sqrt\eta\,\{\log(en)\}^{u+1/2}\sqrt{\log(ed)}
 +\Pp(\|Q-\Sigma\|_{\max}>\eta).
 \label{quant:eq:probability-budget}
\end{align}
For \(p<\infty\), balance the last two terms by taking
\[
 \eta=
 \delta^{2p/(2p+1)}
 \{\{\log(en)\}^{1+2u}\log(ed)\}^{-1/(2p+1)}.
\]
Markov's inequality proves \eqref{quant:eq:joint-rates}.
For \(p=\infty\), take \(\eta=\delta\); for \(\delta=0\),
apply the pathwise bound directly at \(\eta=0\).

Now assume only marginal coordinate exponential moments.
Let \(Y_i=\|X_i\|_\infty\), choose
\(v=K\{C\log(2end)\}^{1/\alpha}\), and set
\[
 X_i'=X_i\ind_{\{Y_i\le v\}}
       -\E_{i-1}[X_i\ind_{\{Y_i\le v\}}].
\]
The new increments satisfy \(\|X_i'\|_\infty\le2v\). The union bound
\(\Pp(Y_i>t)\le2d e^{-(t/K)^\alpha}\) and tail integration imply,
with \(C\) sufficiently large,
\begin{equation}\label{quant:eq:marginal-truncation}
 \E\|S_n-S_n'\|_\infty+
 \E\|Q-Q'\|_{\max}\le C_K(2nd)^{-80}.
\end{equation}
For clarity, the first term is bounded by
\(2\sqrt n\max_i\E[Y_i\ind_{\{Y_i>v\}}]\).
For the second term, the covariance removed at each step is the sum of
the omitted second moment and the outer product of the conditional
truncation mean. Its maximum entry is at most twice the conditional
tail second moment, and hence
\[
 \E\|Q-Q'\|_{\max}
 \le2\max_i\E[Y_i^2\ind_{\{Y_i>v\}}].
\]

Set \(\varepsilon=(2nd)^{-8}\). Except on an event of probability
at most \(C_K(2nd)^{-72}\), both differences in
\eqref{quant:eq:marginal-truncation} are at most \(\varepsilon\).
Use \eqref{quant:eq:probability-budget} for the bounded primed
array, with the original target \(\Sigma\) and tolerance
\(\eta+\varepsilon\). Its exceptional covariance probability is
at most
\[
 \Pp(\|Q-\Sigma\|_{\max}>\eta)+C_K(2nd)^{-72}.
\]
Comparing a rectangle expanded or contracted by \(\varepsilon\)
with the same Gaussian changes its probability by at most
\(Cd\varepsilon\), using the unit marginal variances. The
\(\sqrt\varepsilon\) term and all exceptional probabilities are
absorbed by the sampling term. Thus the same optimization applies.
This argument controls the probability that \(Q'\) differs from
\(\Sigma\) by more than \(\eta+\varepsilon\). Thus it also applies
when \(p=\infty\), without requiring an almost-sure bound on
\(\|Q'-\Sigma\|_{\max}\).

Finally, in the conditional model
\(\{\log(en)\}^{1/\alpha}\le\log(en)\). In the marginal model,
\[
 \{\log(2end)\}^{1/\alpha}
 \le C\{\log(en)\}^{1/\alpha}\{\log(ed)\}^{1/\alpha}.
\]
Consequently the extra factors multiplying the first terms in the
main bounds are at most
\(\{\log(en)\}^{2+3u}\). The extra factor in the covariance term is
\(\{\log(en)\}^{(1+2u)\gamma_p}\), whose exponent is smaller.
Taking \(u=\max\{2,s\}\) proves the assertion.
\end{proof}

\section{Improved bounds when \texorpdfstring{$Q=I_d$}{Q equals Id}}\label{quant:identity-section}

When \(Q=I_d\), we improve the face estimates by retaining the
conditional probability that the coordinates outside each face
satisfy their intervals. Replacing this probability by one would
lose its decay as the number of effective constraints grows.
Relative density comparisons in at most nine coordinates give a
bound for this probability, which leads to
\eqref{main:eq:identity-bounded} and
\eqref{main:eq:identity-marginal}. This argument uses the independent
coordinates of the Gaussian target; the martingale coordinates
need not be independent.

\begin{proposition}\label{quant:prop:identity-bounded}
Let \(n\ge2\), and let \(\xi_1,\ldots,\xi_N\) be martingale
differences in \(\R^d\), where \(N\le3nd\). Suppose
\[
 \|\xi_i\|_\infty\le a,\qquad n^{-1/2}\le a\le1,\qquad
 \sum_{i=1}^N\E_{i-1}\xi_i\xi_i^\top=I_d
 \quad\text{almost surely}.
\]
Then
\begin{equation}\label{quant:eq:identity-bounded}
 \rho_{\cR}\!\left(\sum_{i=1}^N\xi_i,G_{I_d}\right)
 \le C\min\!\left\{1,\,
 a\sqrt{\log(2end)}
       \bigl[\log(2end)+\{\log(en)\}^2\bigr]\right\}.
\end{equation}
The constant is absolute and no restriction is imposed on the
individual conditional covariance matrices.
\end{proposition}

\begin{proof}
We first compare low-dimensional marginal densities at predictable
indices. These comparisons control the conditional probability of
the remaining rectangle constraints. We then sum the third-order
replacement terms using Schur-complement identities and remove the
Gaussian smoothing by scalar moderate-deviation bounds.

Put \(C_i=\E_{i-1}\xi_i\xi_i^\top\),
\(P_i=\sum_{j\le i}C_j\), and \(S_i=\sum_{j\le i}\xi_j\).
Choose a sufficiently large absolute \(D\) and set
\[
 T=D\sqrt{\log(2end)},\qquad
 \Delta=aT\bigl[T^2+\{\log(en)\}^2\bigr].
\]
We may assume that \(\Delta\) is smaller than a fixed constant:
otherwise the unit bound proves the assertion. In particular
\(aT^3\) is small, \(T\le Cn^{1/6}\), and
\(\log(e/\Delta)\le C\log(en)\).
Scalar martingale Bernstein bounds and a Gaussian union bound
allow clipping at \([-T,T]^d\), at cost \((2nd)^{-20}\).
If one interval misses this cube, the probability of the rectangle
is already bounded by this tail estimate. We first treat rectangles
with positive interval lengths; all endpoint conventions follow
at the end by monotone limits.

For a clipped rectangle \(A=\prod_j I_j\), write
\[
 \theta=\sum_j\Pp(G_j\notin I_j),\qquad
 W_A=\sum_{\text{finite faces }b}\phi(b),
 \qquad G_j\sim N(0,1).
\]
The Gaussian hazard bound gives \(W_A\le CT\theta\).
Initially suppose that
\(\theta\le K\log(e/\Delta)\), where \(K\ge3\) is fixed.
Choose
\begin{equation}\label{quant:eq:identity-width}
 \varepsilon=\frac1{100(1+\theta)^3},\qquad
 h=\frac{C_0a[T^2+\{\log(en)\}^2]}{\varepsilon}.
\end{equation}
First take \(C_0\) sufficiently large, and then make the allowed
upper bound on \(\Delta\) smaller. This ensures
\[
 hT\sqrt{\log(e+\theta)}\le c,\qquad h\le1,\qquad h\ge Ca .
\]
Throughout the proof we also permit endpoint movements of size
\(C[h\sqrt{\log(e+\theta)}+T^{-1}]\), with any fixed numerical
coefficient. For nonempty resulting intervals, the sum of the Gaussian
exit probabilities is comparable to \(\theta\), and the sum of the
endpoint densities is at most \(CW_A\).
These facts follow by integrating the one-dimensional Gaussian
hazard, whose size at \(|x|\le2T\) is at most \(C(1+T)\).

\paragraph{Relative densities at predictable indices.}
Let
\[
 A_i=I_d-P_i+h^2I_d,\qquad q_i=(A_i)_{JJ},
 \qquad 1\le |J|=k\le9 .
\]
For every predictable \(\sigma\in\{1,\ldots,N\}\) and
\(\|z\|_\infty\le3T+1\), we claim
\begin{equation}\label{quant:eq:identity-relative-density}
 \E\phi_{q_\sigma}(z-S_{\sigma-1,J})
 =\{1+O(\varepsilon+\Delta)\}
             \phi_{(1+h^2)I_k}(z).
\end{equation}
Stop the projected martingale at \(\sigma-1\), and complete its
total conditional covariance to \(I_k\) by independent signs along
columns of the remaining covariance square root. Subdividing each column into \(M\) signs
with coefficient \(M^{-1/2}\) makes their coordinate bound at most
\(a\) when \(M\ge a^{-2}\). The completed sum \(U_M\) converges
weakly to
\[
 U=S_{\sigma-1,J}
       +N\bigl(0,(I_d-P_{\sigma-1})_{JJ}\bigr),
\]
where the last Gaussian is conditional on the stopped history.
Every completed array has total conditional covariance \(I_k\).

Since \(C_{\sigma,JJ}\preceq ka^2I_k\), the kernel on the left
of \eqref{quant:eq:identity-relative-density} is the conditional
convolution defining \(U\), followed by covariance
\(h^2I_k-C_{\sigma,JJ}\). For positive matrices
\(L\preceq M\preceq U_0\),
\[
 \sqrt{\det L/\det U_0}\,\phi_L
 \le\phi_M\le
 \sqrt{\det U_0/\det L}\,\phi_{U_0}.
\]
Thus it suffices to compare
\(\E\phi_{t^2I_k}(z-U)\) at
\(t^2=h^2-ka^2\) and \(t^2=h^2\); the determinant factors are
\(1+O(a^2/h^2)\).

For finite \(M\), tilt the completed bounded array with the
deterministic parameter \(\lambda=z/(1+t^2)\).
Coordinate reflections allow either sign of its entries.
Lemma~\ref{thr:lem:sparse-tilt} gives the pathwise compensator
\[
 \sum_i\log\E_{i-1}e^{\lambda^\top x_i}
       =|\lambda|^2/2+O(aT^3).
\]
Under its tilted law, it couples \(U_M-\lambda\) to a standard
\(k\)-normal with displacement at most
\(Ca[T^2+\{\log(en)\}^2]\), except on an event of probability
\(n^{-200}\). The bound is uniform in the number of added signs.
By increasing \(C_0\), the Euclidean displacement is at most
\(\varepsilon t\).

Completing the square in the likelihood gives
\begin{align}
 \E\phi_{t^2I_k}(z-U_M)
 ={}&
 e^{-\lambda^\top z+(1+t^2)|\lambda|^2/2}
 \{1+O(\Delta)\}\notag\\
 &\quad{}\times
 \E_\lambda\phi_{t^2I_k}(z-U_M-t^2\lambda).
 \label{quant:eq:identity-density-tilt}
\end{align}
Here a pathwise multiplicative factor bounded between
\(e^{-CaT^3}\) and \(e^{CaT^3}\) is written as \(1+O(\Delta)\).
For \(|v|\le\varepsilon t\),
\[
 (1-C_k\varepsilon)\phi_{(1-\varepsilon)t^2I_k}(x)
 \le \phi_{t^2I_k}(x+v)
 \le (1+C_k\varepsilon)\phi_{(1+\varepsilon)t^2I_k}(x).
\]
The inequalities follow from completing squares and comparing
determinants. In \eqref{quant:eq:identity-density-tilt}, the
Gaussian comparison is evaluated at
\(z-(1+t^2)\lambda=0\), where its smoothed density is bounded
below by a positive constant. The exceptional coupling event
costs at most \(Ct^{-9}n^{-200}=O(\varepsilon)\), since
\(t\ge cn^{-1/2}\) and \(\theta\le C\log(en)\).
Changing \(t^2\) to \((1\pm\varepsilon)t^2\) at this center
has relative cost \(O(\varepsilon)\).
Pass to \(U\) by weak convergence, using the bounded continuous
Gaussian kernel. Finally replacing \(t^2\) by \(h^2\) costs
\(O(a^2(1+T^2))=O(\Delta)\), and
\(a^2/h^2=O(\varepsilon)\).
This proves \eqref{quant:eq:identity-relative-density}.

\paragraph{The probability of the remaining coordinate constraints.}
For ordered distinct coordinates \(J\), \(1\le |J|\le3\), fix
endpoints \(b_J\). Let \(L_i(J,b)\) be the conditional probability
that \(S_{i-1}+N(0,A_i)\) satisfies all remaining rectangle
constraints, given the past and its coordinates \(J\) being
equal to \(b\). The remaining intervals may have any of the
permitted movements. We show that every predictable \(\sigma\)
satisfies
\begin{equation}\label{quant:eq:identity-survival}
 \E\{\phi_{q_\sigma}(b-S_{\sigma-1,J})L_\sigma(J,b)\}
 \le \frac{C}{(1+\theta)^3}\prod_{j\in J}\phi(b_j).
\end{equation}

For each coordinate outside \(J\), select the more probable of its two
Gaussian exits, omitting a coordinate without a finite constraint,
and truncate that exit to \([-3T,3T]\).
Let \(p_j\) be the selected probability under independent
\(N(0,1+h^2)\) coordinates, and let \(\mu=\sum p_j\).
The hazard bound, the permitted endpoint movements, and deletion
of the at most three coordinates in \(J\) give
\[
 c\theta-C\le\mu\le C(1+\theta).
\]
Truncating the selected exits removes at most
\(2d\bar\Phi(3T/\sqrt{1+h^2})\).
If \(\theta\) is bounded, the claimed estimate follows directly
from \eqref{quant:eq:identity-relative-density}.

Otherwise let \(Z\) count the selected exits of the conditional
Gaussian vector given its values \(b\) on \(J\). Satisfying all
remaining interval constraints requires \(Z=0\), so the indicator
of this event is bounded by \((Z-\mu)^6/\mu^6\).
Every joint moment of at most six distinct exit indicators,
multiplied by the marginal density on \(J\) and then averaged at
\(\sigma\),
is an integral of a density in at most \(|J|+6\le9\) coordinates.
Gaussian disintegration and
\eqref{quant:eq:identity-relative-density} therefore identify
that moment with its independent Gaussian reference, with
relative error \(O(\varepsilon+\Delta)\).

Repeated indicators in the sixth-degree expansion reduce to a
single indicator. Under the independent reference,
\[
 \E(Z-\mu)^6\le C(\mu+\mu^2+\mu^3),\qquad
 \E(Z+\mu)^6\le C(1+\mu)^6.
\]
For the first bound, any nonzero term in the expansion in centered
independent indicators uses each index at least twice; there are
at most three distinct indices, and each centered moment of order
at least two is at most \(2p_j\) in absolute value.
The second bound follows by expanding in falling factorials.
It bounds the sum of the absolute reference terms, so the moment
comparison and the indicator bound give
\begin{align*}
 &\E\{\phi_{q_\sigma}(b-S_{\sigma-1,J})L_\sigma(J,b)\}\\
 &\quad\le C\phi_{(1+h^2)I_{|J|}}(b)
 \left\{\frac{\mu+\mu^2+\mu^3}{\mu^6}
 +(\varepsilon+\Delta)\frac{(1+\mu)^6}{\mu^6}\right\}.
\end{align*}
In the present large-\(\theta\) case, \(\mu\asymp\theta\).
Also
\(\varepsilon+\Delta\le C(1+\theta)^{-3}\), because
\(\theta\le K\log(e/\Delta)\). These bounds give the required
factor \((1+\theta)^{-3}\).
The reference marginal density on \(J\) is at most
\(C\prod_{j\in J}\phi(b_j)\), because \(hT\) is small.
This proves \eqref{quant:eq:identity-survival}.
By Lemma~\ref{core:lem:clock}, the same estimate can be multiplied
and summed against any nonnegative predictable weights of
deterministically bounded total mass.

\paragraph{The third-derivative coefficients.}
Write \(F_D(x)=\Pp(x+G_D\in A')\), where \(A'\) is any permitted
rectangle. For a positive definite covariance \(D\), a set \(I\) of already
fixed coordinates, and a new coordinate \(l\), put
\[
 v_{l\mid I}=v_l-D_{lI}D_{II}^{-1}v_I,\qquad
 d_{l\mid I}=D_{ll}-D_{lI}D_{II}^{-1}D_{Il}.
\]
For an ordered pair \(j,k\), write
\[
 q=D_{jj},\quad \gamma=D_{kj}/q,\quad
 d_2=D_{kk}-\gamma^2q,\quad
 z_j=\frac{b_j-x_j}{\sqrt q},\quad
 z_{k\mid j}=\frac{b_k-x_k-\gamma(b_j-x_j)}{\sqrt{d_2}}.
\]
The absolute third directional derivative is bounded by a sum
over ordered faces. Each term is the product of the marginal density on the face, the
conditional probability of the remaining interval constraints, and
the corresponding coefficient in the following table:
\begin{equation}\label{quant:eq:identity-flags}
\begin{array}{c|l}
 J & \text{coefficient}\\ \hline
 (j)& |v_j|^3|z_j^2-1|/q\\[2pt]
 (j,k)&
 2|v_j|^2|v_{k\mid j}||z_j|/\sqrt q
 +|v_j||v_{k\mid j}|^2|z_{k\mid j}|/\sqrt{d_2}\\[2pt]
 (j,k,l)&|v_jv_{k\mid j}v_{l\mid jk}|.
\end{array}
\end{equation}
To verify this representation, whiten the Gaussian and integrate
the third derivative onto successive faces, as in
\citet[arXiv v1, Lemmas~2.2--2.3]{FKLZ2026}.
The first integration gives \(v_j/\sqrt q\).
The other derivatives split into this normal and the orthogonal
tangent space. Two normal derivatives give \(z_j^2-1\); one
normal and one tangent give \(2z_j\); two tangent derivatives
produce a further conormal followed by either a normal derivative
or one more face integration. The successive conormal coefficients
are \(v_{k\mid j}/\sqrt{d_2}\) and
\(v_{l\mid jk}/\sqrt{d_{l\mid jk}}\).
The affine-span Jacobian is the product of the corresponding
standard deviations, namely \(\sqrt{\det D_{JJ}}\).
Multiplying these factors gives exactly
\eqref{quant:eq:identity-flags}, with endpoint signs removed only
after the signed integration identity. Unbounded faces follow by
Gaussian decay.

\paragraph{Removing the Taylor displacements.}
In a within-step remainder, \(D=A_i\),
\(x=S_{i-1}+tv\), \(0\le t\le1\), and \(v\) is either \(\xi_i\)
or a conditional \(N(0,C_i)\) vector.
Restrict first to
\(\|v\|_\infty\le D_1a\sqrt{\log(2end)}\).
After fixing the coordinates \(J\) at their endpoints, the change
in the conditional mean of a coordinate \(l\notin J\) is
\[
 t\{v_l-(A_i)_{lJ}q_i^{-1}v_J\}.
\]
Its absolute value is at most
\[
 Ca\sqrt{\log(2end)}/h\le C/T,
\]
because \((A_i)_{ll}\le2\) and
\(|(A_i)_{lJ}q_i^{-1/2}|\le\sqrt{(A_i)_{ll}}\).
Widening the intervals outside \(J\) by this deterministic amount
therefore removes their dependence on \(v\).
On
\[
 |q_i^{-1/2}(b-S_{i-1,J})|
       \le D_2\sqrt{\log(2end)},
\]
the ratio of shifted to unshifted density on \(J\) is at most
\(\exp\{Ca\log(2end)/h\}\le C\).
The standardized offsets in \eqref{quant:eq:identity-flags} have
uniformly bounded additive changes. Indeed,
\[
 |q_i^{-1/2}tv_J|
 \le \sqrt{|J|}\,\|v_J\|_\infty/h
 \le Ca\sqrt{\log(2end)}/h\le C.
\]
The successive standardized residuals are orthogonal projections
of this whitened vector, so each offset changes by some \(e\)
with \(|e|\le C\). The scalar inequalities
\[
 |z+e|\le C'(1+|z|),\qquad
 |(z+e)^2-1|\le C'(1+z^2),\qquad |e|\le C,
\]
therefore replace the shifted Hermite factors by the displayed
polynomial bounds at the unshifted point. The regression
coefficients and conditional variances do not depend on \(x\).

The discarded regions have negligible total contribution.
On the complement of the restriction on
\(|q_i^{-1/2}(b-S_{i-1,J})|\), a Gaussian density times a polynomial
of degree at most two is bounded by
\(Ch^{-|J|}(2end)^{-cD_2^2}\), including the permitted shifts.
All regression coefficients and inverse variances in a set of at
most three coordinates are bounded by fixed powers of \(h^{-1}\).
Since \(h^{-1}\le C\sqrt n\) and there are at most \(CNd^3\)
terms, a sufficiently large \(D_2\) makes their total at most
\((2nd)^{-20}\).
Outside the restriction on \(\|v\|_\infty\), the original increment
contributes nothing and the conditional Gaussian satisfies
\(\Pp_{i-1}(\|v\|_\infty>t)\le2d e^{-t^2/(2a^2)}\).
The unrestricted derivative bound
\[
 |D_v^3F_{A_i}|\le Ch^{-3}d^{3/2}\|v\|_\infty^3
\]
follows by differentiating the full Gaussian density in its
whitened direction. Tail integration, \(N\le3nd\), and a large
\(D_1\) again give \((2nd)^{-20}\).

\paragraph{Summing the conditional second moments.}
All expectations in this paragraph are conditional on
\(\mathcal F_{i-1}\). For a fixed coordinate set \(I\), \(|I|\le2\), and
coordinate \(l\notin I\), compute the regression from \(A_i\)
and put
\[
 u_{l\mid I,i}=\E_{i-1}v_{l\mid I}^2,\qquad
 d_{l\mid I,i}=(A_i)_{ll}
       -(A_i)_{lI}(A_i)_{II}^{-1}(A_i)_{Il}.
\]
The same quantities arise for the original and Gaussian increments, since
both conditional covariances are \(C_i\).
The decreasing covariance path \(A_{i-1}=A_i+C_i\) satisfies
\begin{equation}\label{quant:eq:identity-energies}
 \sum_i u_{l\mid I,i}\le C,\qquad
 \sum_i\frac{u_{l\mid I,i}}{\sqrt{d_{l\mid I,i}}}\le C,\qquad
 \sum_i\frac{(C_i)_{jj}}{(A_i)_{jj}}\le C\log(en).
\end{equation}
Here is the matrix calculation. With
\(e_i=\E_{i-1}(v_Iv_{l\mid I})\), minimization of the updated
regression loss gives
\[
 d_{l\mid I,i-1}-d_{l\mid I,i}
 =u_{l\mid I,i}
  -e_i^\top\{(A_i)_{II}+(C_i)_{II}\}^{-1}e_i.
\]
Conditional Cauchy--Schwarz gives
\(e_ie_i^\top\preceq u_{l\mid I,i}(C_i)_{II}\).
The subtracted term is therefore at most
\[
 u_{l\mid I,i}
 \tr\{(A_i)_{II}^{-1}(C_i)_{II}\}
 \le Ca^2h^{-2}u_{l\mid I,i}
 \le u_{l\mid I,i}/2.
\]
Also \(A_{i-1,JJ}\preceq2A_{i,JJ}\) on every relevant submatrix,
so consecutive Schur complements have ratio at most two.
The update formula therefore bounds \(u_{l\mid I,i}\) by twice
the decrease of \(d_{l\mid I,i}\). It also gives
\[
 \frac{u_{l\mid I,i}}{\sqrt{d_{l\mid I,i}}}
 \le 2(1+\sqrt2)
 \left(\sqrt{d_{l\mid I,i-1}}-\sqrt{d_{l\mid I,i}}\right).
\]
Summing these two inequalities proves the first two bounds of
\eqref{quant:eq:identity-energies}, since the initial Schur
complement is at most \(1+h^2\). Telescoping scalar logarithms
proves the third.
For the pair regression \(\gamma_i=(A_i)_{kj}/(A_i)_{jj}\),
\(\gamma_i v_j=v_k-v_{k\mid j}\) also gives
\begin{equation}\label{quant:eq:regression-energy}
 \sum_i\gamma_i^2(C_i)_{jj}\le C,
\end{equation}
using \(\sum_i(C_i)_{kk}=1\) and the first bound above.

\paragraph{Centering the Gaussian offsets.}
For one or two fixed coordinates, use the deterministic vector
\(\lambda=b\) and set
\[
 z_i^\circ=q_i^{-1/2}
       (b-S_{i-1,J}-q_i\lambda).
\]
On \(|z_i^\circ|\le D_3\sqrt{\log(en)}\), the singleton
coefficient satisfies
\[
 1+|z_j|^2\le
 C\{\log(en)+qT^2\}.
\]
For a pair,
\[
 1+|z_j|\le
 C\{\sqrt{\log(en)}+\sqrt q\,T(1+|\gamma_i|)\},\qquad
 1+|z_{k\mid j}|\le
 C\{\sqrt{\log(en)}+\sqrt{d_2}\,T\}.
\]
Indeed the mean projections onto the first two conormals are
\(\sqrt q(b_j+\gamma_i b_k)\) and \(\sqrt{d_2}b_k\);
the centered projections have norm at most \(|z_i^\circ|\).

After conditional expectation, the resulting coefficient envelopes
are predictable. Their total masses, for a fixed set of endpoints,
are respectively
\begin{equation}\label{quant:eq:identity-weight-budgets}
 Ca\bigl[\{\log(en)\}^2+T^2\bigr],\qquad
 Ca\{\log(en)+T\},\qquad Ca
\end{equation}
for faces with one, two and three fixed coordinates.
For the singleton use
\(\E_{i-1}|v_j|^3\le Ca(C_i)_{jj}\) and
\eqref{quant:eq:identity-energies}.
For a pair write \(c_i=(C_i)_{jj}\) and
\(u_i=u_{k\mid j,i}\). Both original and Gaussian increments satisfy
\[
 \E_{i-1}|v_j|^2|v_{k\mid j}|
       \le Ca\sqrt{c_i u_i},\qquad
 \E_{i-1}|v_j||v_{k\mid j}|^2\le Ca u_i .
\]
For the original increment these follow from \(|v_j|\le a\) and
Cauchy--Schwarz. For the Gaussian increment use its second and fourth
moments and \(\sqrt{c_i}\le a\).
The first pair coefficient is summed using
\[
 \sqrt{\log(en)}
 \sqrt{\sum_i c_i/(A_i)_{jj}}\sqrt{\sum_i u_i}
 +T\sqrt{\sum_i c_i(1+|\gamma_i|)^2}\sqrt{\sum_i u_i}
 \le C\{\log(en)+T\}.
\]
The other pair coefficient uses
\(\sum_i u_i/\sqrt{d_2}\le C\).
For three fixed coordinates,
\[
 \E_{i-1}|v_jv_{k\mid j}v_{l\mid jk}|
       \le Ca\sqrt{u_{k\mid j,i}u_{l\mid jk,i}},
\]
by boundedness and Cauchy--Schwarz, or by Gaussian moments.
Another application of Cauchy--Schwarz to the sum proves the third
bound in \eqref{quant:eq:identity-weight-budgets}.

The centered-offset complement can also be discarded without a
factor \(N\). Outside \(|z_i^\circ|\le D_3\sqrt{\log(en)}\),
the same predictable coefficient envelopes can be multiplied by
\(C(1+|z_i^\circ|^2)\).
Tilt the original array by the same deterministic \(\lambda=b\)
at every time. For a nonnegative predictable weight \(w_i\)
and nonnegative function \(q_0\), completing the square and
bounding the prefix compensator gives
\begin{align*}
 &\E\{w_iq_0(z_i^\circ)
              \phi_{q_i}(b-S_{i-1,J})\}\\
 &\quad\le
 e^{-\lambda^\top b+(1+h^2)|\lambda|^2/2+CaT^3}
 \E_\lambda\{w_iq_0(z_i^\circ)
              \phi_{q_i}(b-S_{i-1,J}-q_i\lambda)\}.
\end{align*}
The omitted covariance term is
\(-\lambda^\top C_{i,JJ}\lambda/2\le0\).
The prefactor is at most \(C\prod_{j\in J}\phi(b_j)\).
On the discarded centered region, the remaining density times
\(1+|z_i^\circ|^2\) is at most
\[
 Ch^{-|J|}\exp\{-cD_3^2\log(en)\}.
 \]
The pathwise bounds \eqref{quant:eq:identity-weight-budgets}
hold under this one tilted measure, so their sum brings no factor
\(N\). Summing endpoint products costs at most \(W_A^{|J|}\),
which is polynomial in \(n\), since
\(\theta\le C\log(en)\) and \(T\le Cn^{1/6}\).
Increasing \(D_3\) makes the total error \(O(n^{-20})\).

\paragraph{The smoothed comparison.}
The Gaussian heat identity for \(A_{i-1}=A_i+C_i\), followed by
matching the first two conditional moments, telescopes the
within-step remainders just estimated. Apply
\eqref{quant:eq:identity-survival} to their predictable weights
after the deterministic widening that removed the Taylor shifts.
Lemma~\ref{core:lem:clock} and
\eqref{quant:eq:identity-weight-budgets} give
\begin{align}
 &\left|\Pp(S_N+hG'\in A')
         -\Pp(\sqrt{1+h^2}G\in A')\right|\notag\\
 &\quad\le
 \frac{Ca}{(1+\theta)^3}
 \left(
 [\{\log(en)\}^2+T^2]W_A
 +[\log(en)+T]W_A^2+W_A^3
 \right)+Cn^{-20}\notag\\
 &\quad\le C\Delta.
 \label{quant:eq:identity-smoothed}
\end{align}
Here \(G,G'\) are independent standard Gaussian vectors, with
\(G'\) independent of the martingale. The final inequality uses
\(W_A\le CT\theta\) and
\(T^2\log(en)\le[T^3+T\{\log(en)\}^2]/2\).
The estimate holds for all the permitted nonempty moved
rectangles \(A'\).

\paragraph{Removing the smoothing.}
The required scalar input is a direct consequence of
\citet[Theorems 2.1--2.2]{fan2013cramer}. For every coordinate
and \(|t|\le2T\),
\begin{equation}\label{quant:eq:identity-scalar}
 |\Pp(S_{N,j}\le t)-\Phi(t)|
 \le Ca\{1+t^2+\log(en)\}\phi(t).
\end{equation}
To check the hypotheses, the scalar predictable quadratic variation
is exactly one and,
for every integer \(r\ge3\),
\[
 |\E_{i-1}\xi_{ij}^r|
 \le a^{r-2}\E_{i-1}\xi_{ij}^2
 \le (r!/2)a^{r-2}\E_{i-1}\xi_{ij}^2.
\]
Thus the variance-relative Bernstein parameter is \(a\), even
when an individual conditional variance is zero. The cited
theorems assume a trivial initial sigma-field. For arbitrary
\(\mathcal F_0\), prepend a zero increment and reveal
\(\mathcal F_0\) at that step. This preserves the sum, predictable quadratic variation
and Bernstein parameter, and the constants are independent of the
number of increments. The cited logarithmic tail-ratio error is at most
\(Ca[|t|^3+(1+|t|)\log(1/a)]\).
It is small on the present range, and
\(\log(1/a)\le\frac12\log n\).
Multiplication by the appropriate normal tail and Mills'
inequality prove \eqref{quant:eq:identity-scalar}.
Reflection treats the other sign. Strict endpoint conventions
follow by taking one-sided endpoint limits.

Choose
\[
 \gamma=\sqrt{C_4\log(e+\theta)},\qquad r=h\gamma .
\]
A crossing from \(A\) to outside its \(r\)-enlargement after
adding \(hG'\) requires, for some finite endpoint \(b\),
\[
 |G'_j|>\gamma,\qquad |S_{N,j}-b|\le h|G'_j|.
\]
The same conditions are necessary for a crossing from outside
\(A\) into its \(r\)-contraction. If the contraction is empty,
the latter event is empty. For
\(|G'_j|\le D_4\sqrt{\log(2end)}\), the strip endpoints have
absolute value at most \(2T\). Formula
\eqref{quant:eq:identity-scalar} and the Gaussian density bound
give its probability, conditionally on \(G'_j\), at most
\[
 C\phi(b)e^{Th|G'_j|}
 \{h|G'_j|+a[1+T^2+\log(en)]\}.
\]
The complement over all coordinates costs \((2nd)^{-20}\)
for a large \(D_4\). Gaussian tail integration over
\(|G'_j|>\gamma\), followed by summation over endpoints, bounds
the crossing probability by
\[
 C\{h+a[1+T^2+\log(en)]\}
 W_Ae^{-c\gamma^2}+(2nd)^{-20}\le C\Delta
\]
when \(C_4\) is large enough. The same calculation applies with
the original sum replaced by \(G\).

Apply \eqref{quant:eq:identity-smoothed} to the enlarged and
contracted rectangles, and use these crossing bounds on both
sides. It remains to bound the Gaussian effect of moving all
endpoints by at most \(2r\). Along this movement the sum of the Gaussian exit probabilities
is at least \(c\theta\). The derivative of a Gaussian product
probability with respect to one endpoint is at most
\(C\phi(b)e^{1-c\theta}\): deleting one coordinate removes at
most one unit of exit probability. Integration therefore gives
\[
 CrW_Ae^{-c\theta}\le C\Delta.
\]
If a contraction collapses an interval, integrate only up to
the first collapse; its product probability is then zero.
This proves the desired bound when
\(\theta\le K\log(e/\Delta)\).

If \(\theta>K\log(e/\Delta)\), remove whole coordinate constraints
to obtain a containing cylinder for which the sum of the Gaussian
exit probabilities lies in
\([K\log(e/\Delta)-1,K\log(e/\Delta)]\).
Each coordinate contributes at most one, so such a cylinder
exists. Its Gaussian probability is at most \(e\Delta^K\).
The result just proved bounds its martingale probability by
\(C\Delta\), and therefore also bounds both probabilities of
the original rectangle. This completes the proof.
\end{proof}

\begin{corollary}\label{quant:cor:identity-light-tails}
Suppose \(Q=I_d\) almost surely and
\(\log(ed)\ge\{\log(en)\}^2\).
In the bounded model,
\[
 \rho_{\cR}(S_n,G_{I_d})
 \le C_B\min\{1,\{\log(ed)\}^{3/2}/\sqrt n\}.
\]
Under marginal coordinate exponential moments of order
\(\alpha\in\{1,2\}\),
\[
 \rho_{\cR}(S_n,G_{I_d})
 \le C_K\min\{1,\{\log(ed)\}^{3/2+1/\alpha}/\sqrt n\}.
\]
\end{corollary}

\begin{proof}
The bounded assertion follows from
Proposition~\ref{quant:prop:identity-bounded}, since
\(\log(2end)\le C\log(ed)\) in the stated dimension range.
We give the exact covariance completion for the marginal case.
Use the globally truncated, conditionally centered increments from
\eqref{quant:eq:marginal-truncation}, with cutoff
\(v=K\{C\log(2end)\}^{1/\alpha}\), and write their normalized
covariances as \(\widehat C_i\). Their sum is at most \(I_d\),
because the covariance loss is positive semidefinite. Let
\[
 D=I_d-\sum_i\widehat C_i\succeq0.
\]
Its maximum diagonal has expectation at most
\(C_K(2nd)^{-80}\). Put \(a=2v/\sqrt n\). In the nontrivial
range of the asserted bound \(a<1\); also \(a\ge n^{-1/2}\).
After the original increments append \(M=\lceil a^{-2}\rceil\)
independent signs along each column \(D^{1/2}e_j/\sqrt M\).
Every added coordinate is bounded by \(a\), since \(D_{jj}\le1\).
The completed array has total conditional covariance \(I_d\) and at most
\(n+d(n+1)\le3nd\) steps.

Conditionally on the original terminal data, the added sum has
coordinate exponential moments bounded by
\(\exp(t^2D_{jj}/2)\). On
\(\max_jD_{jj}\le(2nd)^{-40}\), its maximum norm is at most
\((2nd)^{-10}\), except with probability at most
\((2nd)^{-30}\), by the conditional exponential bound and a
coordinate union bound. Markov's inequality and
\eqref{quant:eq:marginal-truncation} show that, on the same
probability space, the completed bounded sum differs from
\(S_n\) by at most \(C(2nd)^{-10}\) outside an event of
probability at most \(C_K(2nd)^{-30}\).
Gaussian rectangle anti-concentration
\citep[Theorem 1]{CCKNazarov2017}, followed by
Proposition~\ref{quant:prop:identity-bounded}, gives
\[
 \rho_{\cR}(S_n,G_{I_d})
 \le C_K\min\!\left\{1,\,
 \frac{\{\log(2end)\}^{1/\alpha}}{\sqrt n}
 \sqrt{\log(2end)}
 [\log(2end)+\{\log(en)\}^2]\right\}+C_Kn^{-2}.
\]
The dimension assumption reduces this to the stated bound.
\end{proof}

\section{The exact light-tail convergence range}\label{thr:section}

The factors involving \(\log n\) in the quantitative bounds leave some
sequences below the dimension threshold untreated. To obtain the full
convergence range, we compare intersections of a fixed number of block
exits, with errors measured relative to products of Gaussian tail
probabilities. A fixed collection of blocks is treated jointly; these are
the core blocks. The remaining blocks have small exit probabilities.
The completion and exponential-tilting arguments of
Section~\ref{core:section} give the required comparisons on fixed
collections of blocks. The correlation condition
\eqref{thr:eq:row-power} then controls the sum over all such collections.
The same argument will be used under finite-moment assumptions.

\subsection{Removing a small conditional completion}

The next lemma removes a small perturbation while preserving an error
bound relative to the reference event probability. We will apply it to
Gaussian and sign-vector completions, and later to logistic smoothing.

\begin{lemma}\label{thr:lem:completion-transfer}
Let \(U=Y+R\), where \(Y\) is measurable with respect to a
sigma-field \(\mathcal F\). Let \(D_w\) be a positive diagonal
matrix, and suppose that
\[
 \Pp(\|D_wR\|_\infty>z\mid\mathcal F)
       \le2k\exp\{-z^2/(4c)\},\qquad z>0.
\]
Let \(E^{+z},E^{-z}\) denote the weighted expansion and
contraction of a fixed intersection of coordinate constraints.
Suppose that \(c=c_n\to0\), that \(g=g_n>0\), and that
\[
 \Pp(U\in E^{+z})\le Cg e^{Cz}\quad(z\ge0).
\]
If, uniformly for \(z=o(1)\), the probabilities for \(U\) and a
reference vector \(G\) on \(E^{\pm z}\) differ by \(o(g)\), and
\[
 \Pp(G\in E^{+z}\setminus E^{-z})\le Cgz e^{Cz},
\]
then
\[
 |\Pp(Y\in E)-\Pp(G\in E)|=o(g).
\]
The same conclusion holds if the conditional noise bound is replaced
by
\[
 \Pp(\|D_wR\|_\infty>z\mid\mathcal F)
       \le C_0e^{-c_0z/\epsilon_n},\qquad \epsilon_n\to0,
\]
with fixed positive \(C_0,c_0\).
\end{lemma}

\begin{proof}
The case \(c=0\) is immediate. Put \(z_0=c^{1/4}\) and
\(\pi_j=2k\exp\{-j^2/(4\sqrt c)\}\).
Conditioning on \(\mathcal F\) gives
\[
 (1-\pi_1)\Pp(Y\in E)\le\Pp(U\in E^{+z_0}).
\]
If \(U\in E^{-z_0}\), \(Y\notin E\), and
\(jz_0<\|D_wR\|_\infty\le(j+1)z_0\), then
\(Y\in E^{+jz_0}\).
The conditional noise bound and the first inequality, now applied
to an expanded event, give
\begin{align*}
 &\Pp(Y\in E^{+jz_0},\,\|D_wR\|_\infty>jz_0)
       \le\pi_j\Pp(Y\in E^{+jz_0}),\\
 &(1-\pi_1)\Pp(Y\in E^{+jz_0})
       \le\Pp(U\in E^{+(j+1)z_0}).
\end{align*}
Summing the shells therefore gives
\[
 \Pp(Y\in E)\ge\Pp(U\in E^{-z_0})
 -\frac1{1-\pi_1}
      \sum_{j\ge1}\pi_j\Pp(U\in E^{+(j+1)z_0}).
\]
The sum is \(o(g)\), by the arbitrary-expansion bound and the
Gaussian decay of \(\pi_j\). The two small expansions of \(G\)
differ from its original probability by \(o(g)\).
These inequalities prove the result, also when a contraction is
empty.
For the exponential bound take \(z_0=\sqrt{\epsilon_n}\) and
\(\pi_j=C_0\exp\{-c_0j/\sqrt{\epsilon_n}\}\). The same two
conditional inequalities apply, and
\(\sum_{j\ge1}\pi_j e^{C(j+1)z_0}\to0\).
\end{proof}

\subsection{Relative comparison on a fixed collection of blocks}

The next comparison keeps the Gaussian tail probability in each error
bound. Let \(T=T_n\ge1\), \(T\le n\).
For an array with deterministic predictable quadratic variation in a
fixed dimension, suppose that every
tilt \(\|\lambda\|_\infty\le T\) has cumulant error at most
\(\epsilon_n\to0\), and admits a coupling
\begin{equation}\label{thr:eq:tilt-interface}
 \Pp_\lambda\{\|U-\Gamma\lambda-Z_\Gamma\|_\infty>\tau_n\}
       \le n^{-A},\qquad T\tau_n\to0,
\end{equation}
for every fixed \(A>0\). Constants may depend on the fixed
dimension and on \(A\). Lemmas~\ref{thr:lem:sparse-tilt}
and~\ref{thr:lem:conditional-tilt} establish the cumulant and coupling
bounds used here.

\begin{lemma}\label{thr:lem:exact-block-comparison}
Use the partition, weights, and event \(E\) of
Lemma~\ref{thr:lem:weighted-gaussian-shifts}. Call its first group of
coordinates the core, and suppose that all rare thresholds are at most
\(T\). Suppose that the array with sum \(U\) has
predictable quadratic variation \(\Gamma=V_0+cD_w^{-2}\),
\(0\le c\le1\), almost surely, and that \(U\) satisfies
\eqref{thr:eq:tilt-interface} and the preceding cumulant bound.
Then, uniformly for weighted endpoint movements \(z=o(1)\),
\begin{equation}\label{thr:eq:exact-relative}
 |\Pp(U\in E^{\pm z})-\Pp(G_\Gamma\in E^{\pm z})|=o(g).
\end{equation}
For every \(z\ge0\),
\begin{equation}\label{thr:eq:exact-expansion}
 \Pp(U\in E^{+z})\le C_k g e^{C_kz}.
\end{equation}
No tilt is applied to the core coordinates. Their rectangle may be
arbitrary, and the covariance within each block may be singular.
\end{lemma}

\begin{proof}
In each rare block choose the coordinate with largest selected threshold
\(x_B\). Its variance is
\(\sigma_B^2=1+c/(1+b_B)^2\). Tilt only these coordinates, with
\(\lambda_B=x_B/\sigma_B^2\), and put all other tilt coordinates
equal to zero. The block diagonal target gives
\[
 a=\exp\{-\tfrac12\lambda^\top\Gamma\lambda\}
   =\prod_B e^{-x_B^2/(2\sigma_B^2)},\qquad
 \frac a g\le C_k(1+T)^k.
\]
Here and below the number of rare blocks is at most the fixed total
dimension \(k\). If \(f(u)\) is the inverse likelihood weight
\[
 f(u)=\exp\left\{-\sum_B\lambda_B(u_B-x_B)\right\},
\]
then
\[
 e^{-\epsilon_n}a\,\E_\lambda[f(U)\ind_{\{U\in F\}}]
 \le\Pp(U\in F)
 \le e^{\epsilon_n}a\,\E_\lambda[f(U)\ind_{\{U\in F\}}],
\]
for every event \(F\) imposed on the terminal sum.
For the tilted Gaussian the same identity holds with equality and
without \(e^{\pm\epsilon_n}\).

On the successful coupling event the weighted displacement is at most
\(\Delta_n=(1+T)\tau_n=o(1)\), and the likelihood weights differ by
a factor between \(e^{-kT\tau_n}\) and \(e^{kT\tau_n}\).
Thus the Gaussian terms can be bounded on the corresponding
expansion and contraction by \(\Delta_n\).
Lemma~\ref{thr:lem:weighted-gaussian-shifts} bounds the resulting
boundary error by \(C_kg\Delta_ne^{C_k(z+\Delta_n)}\).
On an expanded event the inverse weight is at most \(e^{2kz}\).
The failed-coupling contribution, after the prefactor is restored,
is therefore at most
\[
 C_kg\,n^{-A}(1+T)^ke^{C_k(z+\Delta_n)}.
\]
Choose \(A>k+2\). These estimates prove
\eqref{thr:eq:exact-relative} when \(z=o(1)\).
They also prove \eqref{thr:eq:exact-expansion} for every \(z\ge0\),
using the Gaussian expansion bound and keeping the original tilts.
No lower bound on the probability of the entire event is used:
the normalizing factor \(g\) is the product of the selected
one-coordinate Gaussian tail probabilities.
\end{proof}

We next explain the completion needed when the target is only
approximately block diagonal. Consider a fixed signed coordinate
projection with sum \(Y\), predictable quadratic variation \(P\), and
a deterministic comparison correlation matrix \(V\). We do not
require \(V=\E P\). Let \(V_0\) retain the entire core block and
each rare block of \(V\), and set the covariances between these
groups to zero. Suppose that
\[
 \|P-V\|_{\max}\le\eta,\qquad
 \|D_w(V-V_0)D_w\|_{\max}\le\zeta .
\]
The diagonal entries of \(D_w\) are at most \(1+T\).
For a sufficiently large fixed-dimensional constant choose
\begin{equation}\label{thr:eq:weighted-completion}
 c=C_k\{\eta(1+T)^2+\zeta\},\qquad
 \Gamma=V_0+cD_w^{-2}.
\end{equation}
Then \(\Gamma-P\succeq0\), \(\Gamma-V\succeq0\), and their
diagonal entries are at most \(2c/w_j^2\).
These assertions follow by bounding the operator norm of a
\(k\)-dimensional matrix by \(k\) times its maximum entry, after
conjugation by \(D_w\).

Given the original terminal history, append \(n\) fresh sign-vector
increments
\begin{equation}\label{thr:eq:local-sign-completion}
 n^{-1/2}(\Gamma-P)^{1/2}\varepsilon_\ell,
 \qquad 1\le\ell\le n,
\end{equation}
where the coordinates of each \(\varepsilon_\ell\) are independent
signs. The completed array has predictable quadratic variation \(\Gamma\).
Each new increment has coordinate bound \(C_k/\sqrt n\) when \(c\le1\),
and
\begin{equation}\label{thr:eq:weighted-completion-tail}
 \Pp(\|D_w(U-Y)\|_\infty>z\mid\mathcal F_N)
      \le2k e^{-z^2/(4c)}.
\end{equation}
Indeed, the conditional moment generating function of any coordinate
is bounded by that of a Gaussian with the same variance, by
\(\cosh t\le e^{t^2/2}\).
A Gaussian completion from \(G_V\) to \(G_\Gamma\) has the same
tail bound.

\begin{corollary}\label{thr:cor:weak-comparison}
In this setting, assume \(c\to0\), and assume that the completed
array satisfies the hypotheses of
Lemma~\ref{thr:lem:exact-block-comparison}. Then
\begin{equation}\label{thr:eq:weak-relative}
 |\Pp(Y\in E)-\Pp(G_{V_0}\in E)|
 +|\Pp(G_V\in E)-\Pp(G_{V_0}\in E)|=o(g).
\end{equation}
The conclusion is uniform over the indicated fixed-dimensional
projections and their core rectangles.
\end{corollary}

\begin{proof}
Use Lemma~\ref{thr:lem:completion-transfer} with reference
\(G_\Gamma\), invoking
\eqref{thr:eq:exact-relative}--\eqref{thr:eq:exact-expansion}
and the Gaussian weighted boundary bound.
Then use \eqref{thr:eq:weighted-ridge} to replace \(G_\Gamma\)
by \(G_{V_0}\). The Gaussian part follows from the same
completion argument without a distributional comparison error.
\end{proof}

\begin{lemma}\label{thr:lem:strong-comparison}
Fix \(k,\kappa>0\). Let \(Y\) be a \(k\)-coordinate martingale
projection with predictable quadratic variation \(P\), and let \(V\) be a
deterministic correlation matrix satisfying
\(|V_{ij}|\le1-\kappa\) for \(i\ne j\) and
\(\|P-V\|_{\max}\le\eta\). We do not require \(V=\E P\).
Let \(1\le x_i\le T\). Choose a graph such that
\[
 |V_{ij}|x_ix_j\le\zeta
 \quad\text{between distinct components},\qquad
 \zeta\to0,\quad \eta T^2\to0.
\]
Complete the projection to
\(\Gamma=V+k\eta I_k\) by \eqref{thr:eq:local-sign-completion}.
If this completed array satisfies \eqref{thr:eq:tilt-interface}
and has cumulant error \(o(1)\), then, with
\(D_h=1+(h-1)(1-\kappa)\),
\begin{equation}\label{thr:eq:strong-relative-upper}
 \Pp(Y_i>x_i,\ i\le k)
 \le C_{k,\kappa}
 \prod_{i\ {\rm isolated}}\bar\Phi(x_i)
 \prod_{\substack{C\ {\rm component}\\|C|\ge2}}
 \exp\left\{-\frac{\sum_{i\in C}x_i^2}{2D_{|C|}}\right\}.
\end{equation}
\end{lemma}

\begin{proof}
Since \(\Gamma-P\succeq0\), and its diagonal entries are at most
\((k+1)\eta\), the conditional probability that the completion
moves every coordinate by at most \(t=C_k\sqrt\eta\) is at least
\(1/2\). Hence
\[
 \Pp(Y>x)\le2\Pp(U>x-t\mathbf1).
\]
For \(i\) in a component of size \(h\), take
\(\lambda_i=x_i/D_h\). The same eigenvalue bound as in
Lemma~\ref{thr:lem:strong-gaussian} gives
\[
 -\lambda^\top x+\tfrac12\lambda^\top\Gamma\lambda
 \le-\sum_C\frac{\sum_{i\in C}x_i^2}{2D_{|C|}}+C_k.
\]
The cross-component terms are bounded by \(C_k\zeta\), and the
added diagonal covariance contributes at most \(C_k\eta T^2\).

Under the predictable tilt, couple \(U\) to
\(\Gamma\lambda+Z_\Gamma\).
On the successful event, every nonisolated inverse-likelihood
factor is at most \(e^{\lambda_i(t+\tau_n)}\) on the enlarged
coordinate event. Discard these factors with this bound.
The marginal covariance of the isolated Gaussian coordinates has
smallest eigenvalue at least \(1/2\) eventually, since their
off-diagonal covariances are at most \(\zeta\).
Its density is bounded by \(C_k\), independently of its mean.
Integration of the remaining inverse-likelihood factors is therefore
at most
\[
 C_k e^{(t+\tau_n)\sum_i\lambda_i}
            \prod_{i\ {\rm isolated}}x_i^{-1}.
\]
The exponential is bounded, because
\((\sqrt\eta+\tau_n)T\to0\).
The failed event contributes at most
\(n^{-A}e^{t\sum_i\lambda_i}\) before restoring the common
likelihood prefactor. Taking \(A>k+2\) absorbs it into the product
of inverse thresholds, as \(T\le n\).
The cumulant error is \(o(1)\).
Restoring the prefactor and applying Mills' inequalities to the
isolated coordinates proves \eqref{thr:eq:strong-relative-upper}.
\end{proof}

\subsection{From rare intersections to rectangles}

The following proposition reduces rectangle convergence to comparisons
on fixed collections of blocks. Its assumptions can be verified by
exponential tilting under light tails and by positive smoothing under
finite moments.

\begin{proposition}\label{thr:assembly}
Let \(Y_n\in\R^{d_n}\), and let \(\Sigma_n\) satisfy
\eqref{thr:eq:row-power}, with fixed \(s,M\).
Use the block partition of Lemma~\ref{thr:block-partition},
with block size at most \(b_0\), separation \(\kappa>0\), and
cross-block correlations \(r_{BC}\).
Put \(T_n=4\sqrt{\log(ed_n)}\).
Suppose the following three properties hold uniformly.

\begin{enumerate}
\item For every fixed number of selected coordinates, their
rectangle distance from the corresponding Gaussian projection
tends to zero. In addition, for either sign,
\begin{equation}\label{thr:eq:scalar-relative-input}
 \max_{j\le d_n}\ \sup_{0\le x\le T_n}
 \left|\frac{\Pp(\pm Y_{n,j}>x)}{\bar\Phi(x)}-1\right|
 \longrightarrow0 .
\end{equation}
The same assertion holds with a weak inequality at the endpoint.

\item Fix a number of core blocks and a number of distinct rare
blocks. Impose an arbitrary rectangle on the core, and a nonempty
compatible set of signed upper exits on each selected rare block.
Let the rare thresholds in block \(B\) belong to
\([b_B,2b_B]\cap[1,T_n]\), where
\(b_B\ge\beta_n\) for some sequence \(\beta_n\to\infty\).
Set \(w_B=1+b_B\), \(\zeta_n=(1+\beta_n)^{-1}\), and suppose
\[
 r_{BC}w_B\le\zeta_n
       \quad(B\text{ rare},\,C\text{ core}),\qquad
 r_{BC}w_Bw_C\le\zeta_n
       \quad(B,C\text{ distinct rare blocks}).
\]
Write \(x_B\) for the largest selected threshold in \(B\).
If \(G_0\) has covariance \(\Sigma_n\) restricted to the entire core
and to each rare block separately, with these groups independent, then
\begin{equation}\label{thr:eq:weak-input}
 |\Pp(Y_n\in E)-\Pp(G_0\in E)|
       =o(1)\prod_B\bar\Phi(x_B).
\end{equation}

\item In the preceding rare-block setting, without the core,
join \(B,C\) when \(r_{BC}w_Bw_C>\zeta_n\).
Select one signed coordinate in each block, with threshold
\(x_B\in[b_B,2b_B]\), and
write \(D_h=1+(h-1)(1-\kappa)\). Then
\begin{equation}\label{thr:eq:strong-input}
 \begin{aligned}
 \Pp(\text{all selected coordinate exits})
 &\le C_k
 \prod_{B\ {\rm isolated}}\bar\Phi(x_B)
 \\
 &\quad{}\times
 \prod_{\substack{C\ {\rm component}\\|C|\ge2}}
 \exp\left\{-\frac{\sum_{B\in C}x_B^2}{2D_{|C|}}\right\}.
 \end{aligned}
\end{equation}
Here \(k\) is the fixed number of selected rare blocks.
\end{enumerate}
Then
\[
 \sup_{A\in\mathcal R_{d_n}}
   |\Pp(Y_n\in A)-\Pp(G_{\Sigma_n}\in A)|\longrightarrow0 .
\]
The three properties are sequential: they are required for every
sequence of the indicated fixed-dimensional configurations.
No uniformity in the number of selected blocks is required.
\end{proposition}

\begin{proof}
Bounded-dimensional subsequences are covered by the first property,
so assume \(d_n\to\infty\). A signed face is an upper exit of
either a coordinate or its negative.
Remove every face whose positive threshold exceeds \(T_n\).
By \eqref{thr:eq:scalar-relative-input}, the probability of any
removed exit, for either law, is at most
\[
 2d_n\{1+o(1)\}\bar\Phi(T_n)=o(1).
\]
For a remaining block let \(q_B\) be its Gaussian probability
of at least one exit and let \(b_B\) be its smallest signed
threshold. By \eqref{thr:eq:block-tail},
\[
 \bar\Phi(b_B)\le q_B\le2b_0\bar\Phi(b_B).
\]
Write \(\Lambda=\sum_Bq_B\), omitting unconstrained blocks.

It suffices first to prove comparison when \(\Lambda\le A+1\), for any
fixed \(A>0\). Indeed, if \(\Lambda>A+1\), retain successive whole
blocks until the sum of their Gaussian exit probabilities first exceeds
\(A\). They define a containing cylinder for which this sum lies in
\([A,A+1]\).
Lemma~\ref{thr:gaussian-cylinder} bounds its Gaussian probability
by \(C(1+A)^{-1/(s+1)}\). Once comparison is established for
rectangles with a bounded sum of Gaussian exit probabilities, the same
bound, up to \(o(1)\), holds for \(Y_n\).
The original rectangle is contained in this cylinder, so letting
\(A\to\infty\) completes this reduction.

Fix \(A\), and take an arbitrary sequence of rectangles after the removal
of thresholds above \(T_n\),
with \(\Lambda\le A+1\).
Order their block probabilities \(q_{j,n}\) decreasingly, pad with
zeros, and extract a subsequence on which
\(q_{j,n}\to a_j\) for every fixed \(j\).
Then \(\sum_ja_j\le A+1\).
By diagonalization there are \(K_n\to\infty\) for which
\[
 \sum_{j\le K_n}|q_{j,n}-a_j|\to0,\qquad
 \max_{j>K_n}q_{j,n}\to0.
\]
The second assertion also follows from
\(q_{j,n}\le(A+1)/j\).
For fixed \(J\), retain the first \(J\) blocks as a core and remove
blocks \(J<j\le K_n\). Their total exit probability under the
Gaussian law is at most \(\sum_{j>J}a_j+o(1)\).
Under the law of \(Y_n\) it is at most a constant times this
quantity, plus \(o(1)\).
To see this, use the bounded number of faces in each block and
\eqref{thr:eq:scalar-relative-input} for positive thresholds;
a nonpositive threshold has Gaussian tail at least \(1/2\)
and satisfies the same constant upper bound trivially.

All remaining blocks are rare: their minimum thresholds are at
least some \(\beta_n\to\infty\).
Remove their secondary faces whose thresholds exceed \(2b_B\).
This has total probability \(o(1)\) under both laws, since
\[
 \sup_{b\ge\beta_n}
       \frac{\bar\Phi(2b)}{\bar\Phi(b)}\longrightarrow0,
 \qquad \sum_Bq_B\le A+1 .
\]
Now every retained rare threshold belongs to \([b_B,2b_B]\),
and remains at most \(T_n\).
Set \(w_B=1+b_B\) and \(\zeta_n=(1+\beta_n)^{-1}\).
Remove a rare block if \(r_{BC}w_B>\zeta_n\) for some core block
\(C\). Its total \(q_B\)-mass is at most
\begin{equation}\label{thr:eq:core-deletion}
 C Jb_0M\,\zeta_n^{-s}
       \sup_{b\ge\beta_n}(1+b)^s\bar\Phi(b)=o(1).
\end{equation}
Assign each removed block to one core block for which the displayed
correlation threshold is exceeded, and apply \eqref{thr:eq:block-correlation}.
The marginal relative upper bound gives the same negligible
deletion probability for \(Y_n\).

Keep \(A,J\), and an intersection order \(k\), fixed.
For each \(k\)-tuple of retained rare blocks, use the graph in
the third property. A tuple is weak if it has no edge, and strong otherwise.
Within each selected block, inclusion--exclusion expresses the
block-exit indicator as a signed sum of intersections of a
nonempty set of its faces. Incompatible opposite exits on the
same coordinate have empty intersection and may be omitted.
There are at most \(C_{k,b_0}\) terms for a fixed tuple.
For each compatible term choose the coordinate with largest threshold
in each block and put \(g=\prod_B\bar\Phi(x_B)\).
Then \(g\le\prod_Bq_B\).
The second property therefore bounds the total weak-tuple
discrepancy from the independent-block Gaussian reference by
\[
 o(1)\sum_{\substack{\text{distinct }B_1,\ldots,B_k\\
                         \text{unordered}}}
                   \prod_{\ell=1}^kq_{B_\ell}
 \le o(1)\frac{(A+1)^k}{k!}=o(1).
\]

The same weak comparison holds for the actual Gaussian vector.
Indeed, in a selected projection its weighted covariance differs
from the independent-block covariance by at most \(C_{J,k}\zeta_n\)
in operator norm. Complete both to
\(V_0+C_{J,k}\zeta_nD_w^{-2}\).
The Gaussian part of Corollary~\ref{thr:cor:weak-comparison},
or directly Lemmas~\ref{thr:lem:weighted-gaussian-shifts}
and~\ref{thr:lem:completion-transfer}, gives an \(o(g)\) error.

For strong tuples we bound the sum of the probabilities themselves.
Consider a nontrivial graph component of size \(h\).
Since \(D_h<h\), there is a constant \(c_{k,\kappa}>0\) such that
\begin{equation}\label{thr:eq:component-gain}
 \frac{\sum_{i\in C}x_i^2}{2D_h}
 \ge\frac12\min_{i\in C}x_i^2
                  +c_{k,\kappa}\max_{i\in C}x_i^2 .
\end{equation}
For completeness, write \(u=\min x_i\), \(v=\max x_i\).
The left side is at least \(\{(h-1)u^2+v^2\}/(2D_h)\);
after subtracting \(u^2/2\), its minimum relative to \(v^2\)
is bounded below by
\(\min\{1,h-D_h\}/(2D_h)>0\).

Choose a spanning tree in each nontrivial component, rooted at a
minimum-threshold vertex. Every required tree edge satisfies
\[
 1\le\zeta_n^{-s}r_{B_iB_j}^{s}w_{B_i}^{s}w_{B_j}^{s}.
\]
Mills' inequality bounds \(e^{-u^2/2}\) by
\(C(1+u)q_{B_{\rm root}}\).
All polynomial endpoint factors are absorbed by half the
exponential gain in \eqref{thr:eq:component-gain}.
Consequently the bound in \eqref{thr:eq:strong-input}, multiplied
by the graph indicator, is dominated by a finite sum of rooted
forest terms of the form
\begin{equation}\label{thr:eq:strong-forest}
 C_k\zeta_n^{-s(k-1)}e^{-c_{k,\kappa}\beta_n^2/2}
 \left(\prod_{ij\ {\rm tree\ edge}}r_{B_iB_j}^{s}\right)
 \prod_{i\ {\rm root\ or\ isolated}}q_{B_i}.
\end{equation}
There is at least one nontrivial component.
Sum first over leaves and then toward the roots.
Summation over each leaf gives a factor at most \(b_0M\), by
\eqref{thr:eq:block-correlation}; summation over each root gives a factor
at most \(A+1\). There are finitely many rooted forests at fixed order.
Since \(\zeta_n^{-1}=1+\beta_n\), the sum of
\eqref{thr:eq:strong-forest} tends to zero.
This proves that all strong tuples contribute \(o(1)\).

For the actual Gaussian vector, apply
Lemma~\ref{thr:lem:strong-gaussian}: covariances between selected
coordinates in distinct components obey
\(|V_{ij}|x_ix_j\le4\zeta_n\).
The independent-block Gaussian reference obeys the same upper bound,
because the selected coordinates belong to independent blocks.
Thus the strong-tuple sum is \(o(1)\) for all three laws.
The core event can be discarded for these upper bounds.

At order zero, the core rectangle probabilities converge by the
fixed-dimensional part of the first property.
At each fixed positive order, the preceding weak and strong
estimates show that the sums of joint rare-block exits, multiplied
by the indicator of the core rectangle, converge to those for the
reference vector in which the core and all rare blocks are independent.
Multiply the two adjacent Bonferroni inequalities for no rare
exits by the same indicator. The independent reference's
remainder at order \(k\) is at most
\((A+1)^{k+1}/(k+1)!\).
Taking \(n\to\infty\) first and then \(k\to\infty\) proves the
comparison after the deletions.
Finally let \(J\to\infty\), removing the error
\(C\sum_{j>J}a_j\).

Every original subsequence admits the further subsequence just
analyzed. Thus the comparison is uniform when the sum of Gaussian exit
probabilities is bounded:
a sequence violating uniformity would retain a fixed positive
error on every further subsequence. The containing-cylinder
reduction then proves the proposition.
\end{proof}

\subsection{Sufficiency under light-tail assumptions}

\begin{proposition}\label{thr:prop:light-validity}
Let \(S_n,Q,\Sigma\) be as in
\eqref{main:eq:normalization}. Assume \(\Sigma_{jj}=1\) and
\(\max_j\sum_k|\Sigma_{jk}|^s\le M\), with fixed \(s,M\). Suppose that
\begin{equation}\label{thr:eq:bracket-in-probability}
 \log(ed)\|Q-\Sigma\|_{\max}\longrightarrow0
                    \quad\text{in probability}.
\end{equation}
Then the rectangle distance from \(S_n\) to \(G_\Sigma\) tends
to zero under either of the following conditions:
\begin{enumerate}
\item \(\|X_i\|_\infty\le B\) almost surely, or
\(\E_{i-1}\exp\{(|X_{ij}|/K)^\alpha\}\le2\) almost surely
for every \(i,j\), with fixed \(\alpha\in\{1,2\}\), provided
\[
 \{\log(ed)\}^{3}/n\longrightarrow0;
\]
\item \(\E\exp\{(|X_{ij}|/K)^\alpha\}\le2\) for every \(i,j\),
with fixed \(\alpha\in\{1,2\}\), provided
\[
 \{\log(ed)\}^{3+2/\alpha}/n\longrightarrow0.
\]
\end{enumerate}
The conclusions are uniform if the convergence in probability in
\eqref{thr:eq:bracket-in-probability} is uniform.
In particular, they hold uniformly under the bound on
\(\|\|Q-\Sigma\|_{\max}\|_p\) in
\eqref{main:eq:setting} whenever \(\delta_n\log(ed_n)\to0\).
\end{proposition}

\begin{proof}
Choose \(\eta_n\) satisfying
\eqref{thr:eq:localization-choice}. We first consider bounded
increments. Lemma~\ref{thr:lem:feasible-completion}, with its sign
completion, reduces the problem to a sum \(Y_n\) whose predictable
quadratic variation \(P\) satisfies
\[
 \|P-\Sigma\|_{\max}\le\eta_n,\qquad
 \eta_n\log(ed)\to0.
\]
The completed array has finitely many increments, each bounded by
\(C_B/\sqrt n\). It agrees with \(S_n\) outside an event of
probability \(o(1)\). Its covariance need not equal
\(\Sigma\); the argument below uses the displayed almost-sure bound
on \(P\) and the original deterministic comparison matrix
\(\Sigma\).

Put \(T=4\sqrt{\log(ed)}\).
For every fixed-dimensional projection, the local completions
\eqref{thr:eq:local-sign-completion} have the same increment
bound, after a fixed change of constant.
The deterministic predictable quadratic variation of each completed
projection has diagonal entries at most two for all sufficiently large
\(n\).
Lemma~\ref{thr:lem:sparse-tilt} supplies
\eqref{thr:eq:tilt-interface}, since
\begin{equation}\label{thr:eq:bounded-validity-errors}
 \frac{T^3+T\{\log(en)\}^2}{\sqrt n}\longrightarrow0.
\end{equation}
Indeed the hypothesis gives \(T=o(n^{1/6})\).
It controls \(T^3/\sqrt n\), and
\(T\{\log(en)\}^2/\sqrt n\to0\) follows directly.
The remaining term in the coupling displacement is bounded using
\[
 2T^2\log(en)\le T^3+T\{\log(en)\}^2.
\]
The cumulant error is \(O(T^3/\sqrt n)=o(1)\).

We verify all three properties of Proposition~\ref{thr:assembly}.
For any fixed coordinate projection, complete to the corresponding
principal submatrix of the original target \(\Sigma\), plus a
vanishing scalar ridge.
Apply the zero tilt in Lemma~\ref{thr:lem:sparse-tilt}, and remove
the completion with Lemma~\ref{thr:lem:completion-transfer},
using the Gaussian boundary bound with no rare blocks.
This gives fixed-dimensional rectangle convergence, including
singular projections.
For a single coordinate and \(1\le x\le T\), use one rare block,
weight \(1+x\), and scale \(g=\bar\Phi(x)\).
Corollary~\ref{thr:cor:weak-comparison} gives a relative error
\(o(\bar\Phi(x))\), uniformly, because
\(\eta_n(1+T)^2\to0\).
For \(0\le x<1\), fixed-dimensional convergence suffices since
\(\bar\Phi(x)\) is bounded below. The same proof treats either
sign and either endpoint convention.
This proves \eqref{thr:eq:scalar-relative-input}.

For each weak configuration in the proposition, its weighted
cross-group covariance is at most \(\zeta_n\), and
\eqref{thr:eq:weighted-completion} has \(c\to0\).
Corollary~\ref{thr:cor:weak-comparison} gives
\eqref{thr:eq:weak-input}.
For a strong configuration retain the coordinate with largest threshold
in each block. Across distinct blocks the absolute correlations
are at most \(1-\kappa\); between graph components,
\[
 |V_{ij}|x_ix_j\le4\zeta_n.
\]
Lemma~\ref{thr:lem:strong-comparison} gives
\eqref{thr:eq:strong-input}.
Thus Proposition~\ref{thr:assembly} proves convergence for \(Y_n\),
and \(Y_n=S_n\) with probability tending to one, which gives the
same conclusion for \(S_n\).
The argument for bounded dimensions is already included in the
fixed-dimensional step.

Now suppose that the coordinate exponential moments hold
conditionally. Use the Gaussian version of the global feasible
completion. The resulting \(2n\) increments have conditional
coordinate tail scale \(C_K/\sqrt n\) in normalized units.
Each fixed projection is then completed by precisely \(n\)
sign-vector increments as in \eqref{thr:eq:local-sign-completion}.
These increments are bounded by \(C_k/\sqrt n\).
The completed projection consequently has at most \(3n\) increments and
satisfies Lemma~\ref{thr:lem:conditional-tilt}.
The cumulant error is \(O(n^{-1/2}T^3)\), and the coupling
displacement times \(T\) is bounded by
\[
 Cn^{-1/2}
 \left\{T\{\log(en)\}^{2+1/\alpha}
              +T^2\log(en)+T^3\right\}=o(1).
\]
Again \(T=o(n^{1/6})\), and the middle term is bounded by the
preceding elementary inequality.
The smallness condition
\(T\{\log(en)\}^{1/\alpha}/\sqrt n\to0\) also holds.
All three properties of Proposition~\ref{thr:assembly} follow
from the same relative block arguments.
The larger truncation threshold used for the coupling under the tilt
does not enlarge the bound \(O(n^{-1/2}T^3)\) on the cumulant error.

It remains to treat marginal exponential moments.
Set
\[
 u=K\{128\log(2end)\}^{1/\alpha},\qquad
 M_i=\|X_i\|_\infty,\qquad
 Y_i=X_i\ind_{\{M_i\le u\}}
             -\E_{i-1}[X_i\ind_{\{M_i\le u\}}],
\]
and let \(S'\) and \(Q'\) be the normalized sum and predictable
quadratic variation of
the \(Y_i\)'s. The coordinate union bound gives
\[
 \Pp(M_i>t)\le2d e^{-(t/K)^\alpha}.
\]
For \(r=1,2\), integration of this bound above \(u\) shows
\[
 \E[M_i^r\ind_{\{M_i>u\}}]\le C_K(2nd)^{-30}.
\]
For example, integrate \(rt^{r-1}2d e^{-(t/K)^\alpha}\)
and include the endpoint term
\(u^r\Pp(M_i>u)\); the fixed polynomial factor is absorbed by
half of the exponent at \(u\).
Consequently
\[
 \E\|S_n-S'\|_\infty
 \le2\sqrt n\max_i\E[M_i\ind_{\{M_i>u\}}].
\]
If \(m_i=\E_{i-1}[X_i\ind_{\{M_i\le u\}}]\), then
\begin{equation}\label{thr:eq:marginal-covariance-loss}
 Q-Q'=\frac1n\sum_i
 \left\{\E_{i-1}[X_iX_i^\top\ind_{\{M_i>u\}}]
                                    +m_im_i^\top\right\}\succeq0.
\end{equation}
Since \(m_i=-\E_{i-1}[X_i\ind_{\{M_i>u\}}]\),
conditional Cauchy--Schwarz gives
\[
 \E\|Q-Q'\|_{\max}
 \le2\max_i\E[M_i^2\ind_{\{M_i>u\}}].
\]
In particular,
\begin{equation}\label{thr:eq:marginal-truncation}
 \E\{\|S_n-S'\|_\infty+\|Q-Q'\|_{\max}\}
       \le C_K(2nd)^{-20}.
\end{equation}
These bounds require no conditional tail assumption.

Let \(a_n=(2nd)^{-2}\). Markov's inequality in
\eqref{thr:eq:marginal-truncation} implies that the two differences
are at most \(a_n\), outside an event of probability \(o(1)\).
Apply the bounded feasible completion to the \(Y_i\)'s at tolerance
\(\eta_n+a_n\) about the original deterministic \(\Sigma\).
Its sum agrees with \(S'\) except on an event of probability
\(o(1)\). All its increments, and all fixed-dimensional local
completion increments, have normalized coordinate bound
\[
 b=C_K\{\log(2end)\}^{1/\alpha}/\sqrt n.
\]
The number of sign-vector increments used for the full-dimensional
completion may depend on \(d\);
Lemma~\ref{thr:lem:sparse-tilt} permits this.
The relative comparison needs exactly
\begin{equation}\label{thr:eq:marginal-validity-errors}
 b\{T^3+T\{\log(en)\}^2\}\longrightarrow0.
\end{equation}
If \(\log(ed)\ge\{\log(en)\}^2\), its left side is bounded by
\[
 C_Kn^{-1/2}\{\log(ed)\}^{3/2+1/\alpha}\longrightarrow0.
\]
If \(\log(ed)<\{\log(en)\}^2\), it is bounded by
\(C_Kn^{-1/2}\{\log(en)\}^{3+2/\alpha}\), which also tends to
zero. Furthermore
\((\eta_n+a_n)(1+T)^2\to0\).
The bounded argument therefore proves rectangle convergence
for \(S'\).

Finally expand and contract each rectangle by \(a_n\) in every
coordinate. The Gaussian boundary bound for unit marginal
variances is
\[
 \sup_{A\in\mathcal R_d}
 \Pp(G_\Sigma\in A^{+a_n}\setminus A^{-a_n})
       \le C a_n\sqrt{\log(ed)}=o(1).
\]
This follows from Nazarov's inequality and does not require a
nonsingular \(\Sigma\)
\citep{Nazarov2003,CCKNazarov2017}.
Together with \eqref{thr:eq:marginal-truncation}, it returns
the conclusion from \(S'\) to \(S_n\).
All constants used at any fixed projection order depend only on
that order and the fixed model parameters.
The subsequence argument in Proposition~\ref{thr:assembly}
therefore proves the asserted uniform versions.
\end{proof}

Propositions~\ref{low:cubic} and~\ref{low:bracket} provide the
matching necessity statements for these light-tail classes.
Their covariance matrices are identities, so they satisfy
\eqref{thr:eq:row-power} for every \(s>0\) and \(M\ge1\).
Together with the upper comparison, these constructions give
the light-tail parts of Theorem~\ref{main:validity}.
Section~\ref{main:proofs} passes from the sequential comparisons
to the class suprema and combines the two necessary conditions.

\section{Proofs under finite moments}\label{poly:section}

This section proves the finite-moment parts of
Theorems~\ref{main:joint}--\ref{main:identity}. We retain
\eqref{main:eq:moment-envelope} and use the normalized moment sum
\[
 v_q=\sum_{i=1}^n\E\|\xi_i\|_\infty^q
 \le B^qn^{1-q/2},\qquad \xi_i=X_i/\sqrt n.
\]
The moment assumption is marginal, so the predictable moment sum
need not be bounded on every path. We first compare positive smoothed
test functions along the sequence of conditional covariance matrices.
We then stop the array when the conditional moment sum is too large
and truncate the retained increments. This gives almost-sure bounds
with deterministic constants. A single exceptional event controls
truncation and stopping for the whole vector; the displacement due to
centering is bounded separately. The convergence argument applies
for every \(q>2\);
the quantitative upper bounds and the matching adaptive lower
construction use \(2<q\le3\).

\subsection{A positive comparison in a fixed projection}

We use the Gaussian operators \(P_A\) and logistic densities \(p_h\)
defined in Section~\ref{quant:section}. The following elementary
comparison permits the covariance loss caused by truncation.

\begin{lemma}\label{poly:lem:positive-comparison}
Let \(y_1,\ldots,y_N\) be \(k\)-dimensional martingale differences,
where \(k\) is fixed. Let \(C_i\) be predictable positive semidefinite
matrices such that
\[
 \widehat C_i=\E_{i-1}y_i y_i^\top\preceq C_i,\qquad
 P_i=\sum_{\ell\le i}C_\ell\preceq\Gamma
 \quad\text{almost surely},
\]
where \(\Gamma\) is deterministic. Suppose that
\(\|y_i\|_\infty\le c_1h\) and
\(\|C_i\|_{\max}\le c_1h^2\). Put
\[
 d_i=\|C_i-\widehat C_i\|_{\max},\qquad
 b_i=\E_{i-1}\|y_i\|_\infty^3,
\]
and assume the pathwise bound
\begin{equation}\label{poly:eq:relative-budget}
 \sum_{i=1}^N
 \left\{\frac{d_i}{h^2}+\frac{b_i}{h^3}
                   +\frac{\|C_i\|_{\max}^2}{h^4}\right\}\le A.
\end{equation}
Let
\[
 W=\sum_{i=1}^N y_i+G_{\Gamma-P_N},
\]
where the last Gaussian is generated conditionally on the terminal
history, using independent randomness. If a positive \(C^4\) function
\(f\) satisfies
\begin{align}
 |\partial_{j_1\cdots j_\ell}f(x)|
   &\le c_2h^{-\ell}f(x),\qquad 1\le\ell\le4,\notag\\
 e^{-c_2\|z\|_\infty/h}f(x)
   &\le f(x+z)\le e^{c_2\|z\|_\infty/h}f(x),
 \label{poly:eq:positive-test}
\end{align}
then
\begin{equation}\label{poly:eq:relative-upper}
 \E f(W)\le e^{C A}\E f(G_\Gamma).
\end{equation}
Here \(C\) depends only on \(k,c_1,c_2\). If \(A\) is sufficiently
small, then also
\begin{equation}\label{poly:eq:relative-two-sided}
 \E f(W)=\{1+O(A)\}\E f(G_\Gamma).
\end{equation}
The assertions are valid for singular \(\Gamma\), arbitrary initial
sigma-fields, and every predictable prefix of the given finite array.
\end{lemma}

\begin{proof}
First take deterministic \(N\), and set
\[
 F_i(x)=P_{\Gamma-P_i}f(x),\qquad Y_i=\sum_{\ell\le i}y_\ell.
\]
The remaining covariance is positive semidefinite at every prefix.
The functions \(F_i\) inherit \eqref{poly:eq:positive-test} and
satisfy \(F_{i-1}=P_{C_i}F_i\) pathwise. Both \(F_i\) and \(C_i\)
are \(\mathcal F_{i-1}\)-measurable. Thus the comparison can be made
conditionally before increment \(i\), even though the laws of later
increments depend on the observed past. Since \(k\) is fixed and
\(\|C_i\|_{\max}\le c_1h^2\), the logarithmic Lipschitz bound gives
\begin{equation}\label{poly:eq:relative-base}
 cF_i(x)\le P_{tC_i}F_i(x)\le CF_i(x),
 \qquad 0\le t\le1.
\end{equation}
Indeed the upper bound follows from
\(\E\exp(c_2\|G_{tC_i}\|_\infty/h)\le C\); the lower bound follows
from Jensen's inequality applied to the reverse shift bound.

Taylor expansion to order two gives
\[
 \E_{i-1}F_i(x+y_i)
 =F_i(x)+\frac12\widehat C_i:\nabla^2F_i(x)+R_i(x),
 \qquad |R_i(x)|\le C h^{-3}b_iF_i(x).
\]
The integral remainder is bounded using
\(\|y_i\|_\infty/h\le c_1\). On the other hand, the heat equation and
\eqref{poly:eq:relative-base} give
\[
 P_{C_i}F_i(x)
 =F_i(x)+\frac12C_i:\nabla^2F_i(x)+R_i'(x),
 \qquad
 |R_i'(x)|\le Ch^{-4}\|C_i\|_{\max}^2F_i(x).
\]
Consequently, at \(x=Y_{i-1}\),
\[
 \left|\E_{i-1}F_i(Y_i)-F_{i-1}(Y_{i-1})\right|
 \le a_iF_{i-1}(Y_{i-1}),\qquad
 a_i=C\left\{\frac{d_i}{h^2}+\frac{b_i}{h^3}
                      +\frac{\|C_i\|_{\max}^2}{h^4}\right\}.
\]
Dividing \(F_i(Y_i)\) by \(\prod_{\ell\le i}(1+a_\ell)\)
produces a nonnegative supermartingale. The product is at most
\(e^{CA}\), and \(F_0(0)=\E f(G_\Gamma)\), proving
\eqref{poly:eq:relative-upper}. If \(CA<1/2\), division by
\(\prod_{\ell\le i}(1-a_\ell)\) instead gives a submartingale.
Its product lies between \(e^{-2CA}\) and one, which proves
\eqref{poly:eq:relative-two-sided}.

The argument only convolves by positive semidefinite covariances.
No density of \(G_\Gamma\) is required. Conditioning on the initial
sigma-field does not change \(F_0(0)\). A predictable terminal index
is handled by suppressing the increment at that index and all later
increments, and then adding the displayed Gaussian completion. This shows that
the constants do not depend on the number of increments.
\end{proof}

Here and below a predictable index \(\tau\) means that
\(\{\tau=i\}\in\mathcal F_{i-1}\); its retained prefix ends at
\(\tau-1\). The lemma applies, in particular, to positive products of
logistic densities and to logistic convolutions of nonnegative
bounded Borel functions.

\begin{corollary}\label{poly:cor:relative-smoothing}
Let \(y_i\) be \(k\)-dimensional martingale differences with
\(\|y_i\|_\infty\le a\), and suppose
\[
 P_N=\sum_i\E_{i-1}y_i y_i^\top\preceq\Gamma,\qquad
 \sum_i\E_{i-1}\|y_i\|_\infty^3\le v
 \quad\text{almost surely}.
\]
For independent standard logistic coordinates \(H\), put
\(W=\sum_i y_i+G_{\Gamma-P_N}\). If \(a\le\sigma\) and
\(v/\sigma^3\) is sufficiently small, then, uniformly over Borel
\(0\le g\le1\),
\begin{equation}\label{poly:eq:relative-smoothing}
 \E g(W+\sigma H)
 =\{1+O_k(v/\sigma^3)\}\E g(G_\Gamma+\sigma H).
\end{equation}
The same assertion holds at any predictable prefix.
\end{corollary}

\begin{proof}
For \(f(x)=\E g(x+\sigma H)\), logarithmic differentiation of
\(p\) gives \eqref{poly:eq:positive-test}, with \(h=\sigma\) and
constants depending only on \(k\). If \(g\) is zero almost
everywhere both sides vanish; otherwise \(f\) is positive.
Take \(C_i=\widehat C_i\). Conditional Jensen and boundedness give
\[
 \|C_i\|_{\max}^2
 \le\E_{i-1}\|y_i\|_\infty^4
 \le a\,\E_{i-1}\|y_i\|_\infty^3.
\]
Thus the sum on the left side of \eqref{poly:eq:relative-budget} is at most
\(v/\sigma^3+av/\sigma^4\le2v/\sigma^3\).
\end{proof}

\subsection[Sufficiency for every q greater than two]
{Sufficiency for every \(q>2\)}

\begin{proposition}\label{poly:prop:validity}
Fix \(q>2\), \(s>0\), \(M\ge1\), and \(B\ge2\). Suppose
\[
 \max_i\E\|X_i\|_\infty^q\le B^q,\qquad
 \max_j\sum_k|\Sigma_{jk}|^s\le M,\qquad \Sigma_{jj}=1.
\]
For an array sequence, assume
\begin{equation}\label{poly:eq:validity-assumptions}
 \log(ed)=o\!\left(n^{\min\{1/3,\,1-2/q\}}\right),
 \qquad
 \log(ed)\|Q-\Sigma\|_{\max}\longrightarrow0
 \quad\text{in probability}.
\end{equation}
Then \(\rho_{\cR}(S_n,G_\Sigma)\to0\).
The convergence is uniform whenever the convergence in probability
in \eqref{poly:eq:validity-assumptions} is uniform.
\end{proposition}

\begin{proof}
The proof requires relative probabilities for fixed collections of
rare exits. An absolute approximation error, even if it tends to
zero, would not suffice when these probabilities are summed over
many coordinates.

\paragraph{Localization under the marginal moment assumption.}
Set \(T_n=4\sqrt{\log(ed)}\). The dimension condition implies
\begin{equation}\label{poly:eq:window-moments}
 n^{-1/2}T_n^3\longrightarrow0,\qquad
 n^{1-q/2}T_n^q\longrightarrow0.
\end{equation}
Choose \(u_n\downarrow0\) slowly enough that
\(B^qn^{1-q/2}T_n^q/u_n^q\to0\), and put
\(a_n=u_n/T_n\).
Truncate the normalized increments at \(a_n\) and center them
conditionally:
\[
 y_i=\xi_i\ind_{\{\|\xi_i\|_\infty\le a_n\}}
       -\E_{i-1}[\xi_i\ind_{\{\|\xi_i\|_\infty\le a_n\}}].
\]
Writing \(Q_y=\sum_i\E_{i-1}y_i y_i^\top\), tail integration gives
\begin{align}
 \Pp(\exists i:\|\xi_i\|_\infty>a_n)
   &\le v_qa_n^{-q}=o(1),\notag\\
 T_n\E\left\|S_n-\sum_i y_i\right\|_\infty
   &\le2T_nv_qa_n^{1-q}=o(1),\notag\\
 T_n^2\E\|Q-Q_y\|_{\max}
   &\le2T_n^2v_qa_n^{2-q}=o(1).
 \label{poly:eq:window-truncation}
\end{align}
The covariance loss is positive semidefinite, by the same
calculation used in \eqref{poly:eq:identity-moment-budgets}.
Let \(V_y=\sum_i\E_{i-1}\|y_i\|_\infty^3\).
For \(2<q<3\),
\[
 \E V_y\le C_qv_qa_n^{3-q};
\]
for \(q\ge3\), conditional centering and marginal Jensen give
\[
 \E V_y\le C_q B^3 n^{-1/2}.
\]
In either case, \eqref{poly:eq:window-moments} yields
\(T_n^3\E V_y\to0\).

Choose deterministic \(v_n,\eta_n>0\) such that
\[
 T_n^3v_n+T_n^2\eta_n\longrightarrow0,\qquad
 \Pp(V_y>v_n\ \text{or}\ \|Q_y-\Sigma\|_{\max}>\eta_n)
 \longrightarrow0.
\]
Retain the \(y_i\) until the next increment would make the sum of
conditional third moments exceed \(v_n\), or would violate the
feasibility condition for a covariance matrix within entrywise distance
\(\eta_n\) of \(\Sigma\).
Lemma~\ref{thr:lem:feasible-completion} supplies a feasible terminal
matrix \(A\) and a conditional Gaussian completion.
Denote the resulting sum by \(\widetilde S\), and the sum of the
conditional covariance matrices of its retained increments by \(P\). Thus
\begin{equation}\label{poly:eq:localized-window}
 \widetilde S=\sum_i I_i y_i+G_{A-P},\qquad
 \|A-\Sigma\|_{\max}\le\eta_n,\qquad
 \sum_i I_i\E_{i-1}\|y_i\|_\infty^3\le v_n .
\end{equation}
The original truncated sum and \(\widetilde S\) agree off an event of
probability \(o(1)\). Each retained increment has maximum norm at
most \(2a_n\). The Gaussian completion is treated exactly by the heat
operator, so it need not be included in the third-moment sum.

\paragraph{Relative comparison after smoothing.}
Choose \(\sigma_n>0\) so that
\begin{equation}\label{poly:eq:window-smoothing}
 a_n/\sigma_n\longrightarrow0,\qquad
 v_n/\sigma_n^3\longrightarrow0,\qquad
 T_n\sigma_n\longrightarrow0 .
\end{equation}
Such a choice is possible because \(T_na_n\to0\) and
\(T_n^3v_n\to0\).
Consider any fixed-dimensional projection of
\eqref{poly:eq:localized-window}, with signs changed when needed.
If a deterministic covariance \(\Gamma\) dominates the projected
\(A\), add a conditional normal of covariance \(\Gamma-A\).
Conditionally on the history of the retained increments, the two added
Gaussian vectors combine into \(G_{\Gamma-P}\).
Corollary~\ref{poly:cor:relative-smoothing} therefore gives
\begin{equation}\label{poly:eq:window-relative}
 \E g(U+\sigma_nH)
 =\{1+o(1)\}\E g(G_\Gamma+\sigma_nH),
 \qquad 0\le g\le1,
\end{equation}
uniformly over Borel tests and over all such projections of a fixed
dimension. Neither \(\Gamma\) nor its within-block restrictions
need be invertible.

We spell out how this smoothed comparison retains the probability
scale needed when summing over collections of blocks.
Use the weights and event \(E\) of
Lemma~\ref{thr:lem:weighted-gaussian-shifts}; their largest weight is
at most \(CT_n\). Suppose the corresponding Gaussian probability
scale is \(g_*>0\). A logistic displacement changes a weighted
endpoint by at most
\[
 T_n\sigma_n C\|H\|_\infty.
\]
The Gaussian arbitrary-expansion and strip bounds in that lemma
imply
\[
 \Pp(G_\Gamma+\sigma_nH\in E^{+z})
 \le Cg_*e^{Cz},\qquad z\ge0,
\]
and, uniformly for \(z=o(1)\),
\[
 \left|\Pp(G_\Gamma+\sigma_nH\in E^{\pm z})
             -\Pp(G_\Gamma\in E^{\pm z})\right|=o(g_*).
\]
Indeed the required expectation is bounded by a constant times
\[
 g_*\,\E\!\left[
 T_n\sigma_n\|H\|_\infty
       e^{CT_n\sigma_n\|H\|_\infty}\right]=o(g_*).
\]
The relative error in \eqref{poly:eq:window-relative}, multiplied by
these Gaussian probabilities, is still \(o(g_*)\).
The exponential-tail version of
Lemma~\ref{thr:lem:completion-transfer} now removes the logistic
noise, on this same relative scale. The upper-bound part of the
argument only needs the Gaussian arbitrary-expansion estimate.

\paragraph{Comparisons on fixed collections of blocks.}
We verify the inputs of Proposition~\ref{thr:assembly}.
For a weak selection of a fixed number of rare blocks and a fixed
core, let \(\Gamma_0\) be its block-diagonal Gaussian covariance,
and let \(D_w\) be its diagonal weight matrix.
The definition of a weak selection gives a deterministic
\(\epsilon_n\to0\) bounding all weighted cross-block correlations.
Choose
\[
 c_n=C_k(\epsilon_n+\eta_nT_n^2),\qquad
 \Gamma=\Gamma_0+c_nD_w^{-2},
\]
with \(C_k\) large enough. On the selected \(k\) coordinates,
\[
 \|D_w(\Sigma-\Gamma_0)D_w\|_{\mathrm{op}}
       \le C_k\epsilon_n,\qquad
 \|D_w(A-\Sigma)D_w\|_{\mathrm{op}}
       \le C_k\eta_nT_n^2.
\]
Hence \(\Gamma\succeq A\), and
the added normal has weighted marginal variances at most \(C_kc_n\).
Let
\[
 g_*=\prod_B\bar\Phi(x_B),
\]
where \(x_B\) is the largest selected threshold in rare block \(B\).
The preceding relative smoothing argument, followed by
Lemma~\ref{thr:lem:completion-transfer}, removes both the logistic
noise and the conditional Gaussian completion. Finally
\eqref{thr:eq:weighted-ridge} removes the deterministic ridge.
It follows that
\[
 \Pp(\widetilde S\in E)
   =\Pp(G_{\Gamma_0}\in E)+o(g_*).
\]
The estimate is uniform over every weak selection of a fixed order,
including its core rectangle and the weighted expansions used in
Proposition~\ref{thr:assembly}.

For a strong tuple, let \(V\) be the signed principal submatrix
of the original comparison matrix \(\Sigma\) on its selected
coordinates. This matrix need not equal the covariance of the
truncated and stopped projection.
Complete instead to \(\Gamma=V+k\eta_nI\).
The marginal variances of the extra conditional normal are
at most \(2k\eta_n\). Thus its displacement on the scale \(T_n\)
has size \(o(1)\).
Lemma~\ref{thr:lem:strong-gaussian} bounds every expanded Gaussian
event by \(Ce^{Cz}\) times
\[
 \prod_{i\ \mathrm{isolated}}\bar\Phi(x_i)
 \prod_{\substack{C\ \mathrm{component}\\|C|\ge2}}
 \exp\left\{-\frac{\sum_{i\in C}x_i^2}
                 {2\{1+(|C|-1)(1-\kappa)\}}\right\}.
\]
The tilt in that lemma is fixed at the unshifted thresholds, so its
arbitrary-expansion bound integrates against the logistic noise.
Apply \eqref{poly:eq:window-relative} and the upper half of the
conditional-completion argument. This gives exactly the strong-tuple
bound required by Proposition~\ref{thr:assembly}.

The one-coordinate instance gives relative normal tails uniformly
through \(T_n\), for either sign and for both strict and weak
inequalities. For thresholds in \([0,1]\), the fixed-dimensional
absolute comparison suffices since the normal tail is bounded
below. For thresholds in \([1,T_n]\), use the single rare block
with weight \(1+x\); the variance inflation is negligible because
\(\eta_nT_n^2\to0\).
Projections of fixed dimension give rectangle
convergence: complete them to
\(\Sigma_{JJ}+|J|\eta_nI\), use
\eqref{poly:eq:window-relative}, and remove the vanishing smoothing
and completion. Their Gaussian boundary probabilities tend to zero
uniformly because every marginal variance is one.
This also proves the result directly on subsequences of bounded
dimension. On subsequences where \(d\to\infty\), the one-coordinate
tail estimate allows thresholds above \(T_n\) to be removed with an
error tending to zero, since
\(d\bar\Phi(T_n)\to0\).
All three assumptions of Proposition~\ref{thr:assembly} are now
verified, so
\[
 \rho_{\cR}(\widetilde S,G_\Sigma)\longrightarrow0.
\]

\paragraph{Return to the original array.}
Equation \eqref{poly:eq:window-truncation} gives an
\(o_{\Pp}(T_n^{-1})\) displacement. The stopped and completed sum in
\eqref{poly:eq:localized-window} differs from the truncated sum only on
an event of probability \(o(1)\).
The Gaussian rectangle boundary bound
\[
 \sup_{A\in\cR}
 \Pp(G_\Sigma\in A^{+t}\setminus A^{-t})
 \le Ct\sqrt{\log(ed)}
\]
follows by applying Nazarov's inequality to the \(2d\) signed
coordinates; see \citet[Theorem~1]{CCKNazarov2017}.
It transfers convergence back to \(S_n\).
All choices used deterministic upper bounds and uniform convergence
in probability, which proves the uniform assertion.
\end{proof}

\subsection[An upper bound of optimal order when Q equals I-d]
{An upper bound of optimal order when \(Q=I_d\)}

\begin{proposition}\label{poly:prop:identity-upper}
Suppose \(Q=I_d\) almost surely and \(2<q\le3\). Then
\begin{equation}\label{poly:eq:identity-upper}
 \rho_{\cR}(S_n,G_{I_d})
 \le C_q\min\left\{1,
          \bigl(v_q\{\log(ed)\}^{q/2}\bigr)^{1/(q+1)}\right\}.
\end{equation}
The bound allows arbitrary predictable covariance paths.
\end{proposition}

\begin{proof}
Set
\[
 T=4\sqrt{\log(ed)},\qquad
 \varepsilon=(v_qT^q)^{1/(q+1)},\qquad h_0=\varepsilon/T.
\]
It is enough to consider \(\varepsilon\) smaller than a fixed
constant. Notice the identity \(v_q=\varepsilon h_0^q\).
Stop before the first increment that would make the sum of conditional
\(q\)th moments exceed \(h_0^q\), and suppress that increment and all
later ones. Write \(I_i\) for the indicator that increment \(i\) is
retained, and put
\[
 z_i=I_i\xi_i,\qquad
 y_i=z_i\ind_{\{\|z_i\|_\infty\le h_0\}}
             -\E_{i-1}[z_i\ind_{\{\|z_i\|_\infty\le h_0\}}],
 \qquad C_i=\E_{i-1}z_i z_i^\top .
\]
The probability of stopping is at most \(\varepsilon\).
The probability that any original increment is truncated is also at most
\(v_qh_0^{-q}=\varepsilon\). On every retained path,
\begin{equation}\label{poly:eq:identity-moment-budgets}
 \begin{split}
 \|y_i\|_\infty&\le2h_0,\qquad
 \sum_i\E_{i-1}\|y_i\|_\infty^3\le C_qh_0^3,\\
 d_i:=\|C_i-\E_{i-1}y_i y_i^\top\|_{\max},\qquad
 \sum_i d_i&\le2h_0^2,\qquad
 \sum_i\|C_i\|_{\max}^2\le C_qh_0^4 .
 \end{split}
\end{equation}
For the covariance claim, the loss equals
\[
 \E_{i-1}\!\left[z_i z_i^\top
                 \ind_{\{\|z_i\|_\infty>h_0\}}\right]
 +\E_{i-1}\!\left[z_i\ind_{\{\|z_i\|_\infty>h_0\}}\right]
  \E_{i-1}\!\left[z_i\ind_{\{\|z_i\|_\infty>h_0\}}\right]^\top .
\]
It is positive semidefinite and has maximum diagonal at most
\(2h_0^{2-q}\E_{i-1}\|z_i\|_\infty^q\). The third-moment bound
follows by truncation and conditional Jensen. Finally,
\[
 \|C_i\|_{\max}^2
 \le2\E_{i-1}\|y_i\|_\infty^4+2d_i^2
 \le4h_0\E_{i-1}\|y_i\|_\infty^3+2d_i^2.
\]
This proves the last inequality in
\eqref{poly:eq:identity-moment-budgets}. On the event of no stopping
and no truncation, the accumulated centering displacement is at most
\[
 h_0^{1-q}\sum_i\E_{i-1}\|z_i\|_\infty^q\le h_0.
\]

\paragraph{Probabilities near prescribed endpoints.}
Write \(P_i=\sum_{\ell\le i}C_\ell\) and
\(Y_i=\sum_{\ell\le i}y_\ell\). Then \(P_i\preceq I_d\).
For any predictable index \(\tau\), complete \(Y_{\tau-1}\)
conditionally by \(N(0,I_d-P_{\tau-1})\). Let \(W_J\) be the
projection onto \(k\le4\) distinct coordinates. For
\[
 u(b)=\operatorname{sgn}(b)\min(|b|,T),\qquad
 \omega(b)=\exp\{-u(b)b+u(b)^2/2\},
\]
we claim, whenever \(h\ge h_0\) and \(hT\le c\), that
\begin{equation}\label{poly:eq:pin-bound}
 \Pp\{\|W_J-b\|_\infty\le c_2h\}
 \le C_{q,c_2}h^k\prod_{j=1}^k\omega(b_j).
\end{equation}
Apply Lemma~\ref{poly:lem:positive-comparison} to
\(f(x)=\prod_{j=1}^k p_h(x_j-b_j)\).
The sum in \eqref{poly:eq:relative-budget} is bounded by a constant
independent of \(n,d\), by
\eqref{poly:eq:identity-moment-budgets}; also
\(\|C_i\|_{\max}\le Ch_0^2\le Ch^2\).
On the indicated box, \(f\ge c h^{-k}\). For a one-dimensional
standard normal \(G\),
\[
 \E p_h(G-b)
 =\int\phi(b+hz)p(z)\,dz
 \le C\omega(b)\int e^{|u(b)|h|z|}p(z)\,dz
 \le C\omega(b).
\]
The first inequality uses
\(\phi(x)e^{u x}\le(2\pi)^{-1/2}e^{u^2/2}\).
The target coordinates are independent, so the product of these
one-dimensional bounds proves \eqref{poly:eq:pin-bound}.

\paragraph{Summing local errors.}
Append \(G_{I_d-P_n}\) to \(Y_n\), and denote the resulting vector
by \(W_0\). For a rectangle \(A=\prod_j I_j\), let
\[
 w=\sum_j\Pp(G_j\notin I_j),\qquad
 W_A=\sum_b\omega(b),
\]
where the second sum runs over its finite signed endpoints.
Empty rectangles and the whole space require no argument.
First suppose \(w\le4\log(e/\varepsilon)\), and use
\[
 h=h_0(1+w).
\]
Then \(hT=\varepsilon(1+w)\) is small. The normal hazard bound
on \(|b|\le T\), and
\(\sum_{|b|>T}\omega(b)\le2de^{-T^2/2}\le1\), give
\begin{equation}\label{poly:eq:identity-weight}
 W_A\le CT(1+w).
\end{equation}

Smooth each interval by a fixed smooth probability density supported
on \([-h,h]\), and take the product of the smoothed intervals.
The absolute value of an ordered derivative of order \(\ell\le4\)
involving \(k\) distinct coordinates is at most \(C_\ell h^{-\ell}\)
times the sum of the
indicators that these coordinates lie within \(h\) of one of their
endpoints. The undifferentiated factors are at most one.
Telescope with the Gaussian covariances \(I_d-P_i\).
The mismatch between \(C_i\) and the covariance of \(y_i\) produces
a second derivative with coefficient at most \(d_i\).
The Taylor remainder for the martingale increment has coefficient
\(\E_{i-1}\|y_i\|_\infty^3\). Expansion of the Gaussian heat
operator through first order produces a fourth derivative with
coefficient at most \(C\|C_i\|_{\max}^2\).
For this Taylor remainder, enlarge the interval around each selected
endpoint by \(2h_0\). Since \(\|y_i\|_\infty\le2h_0\), the indicator of a shifted
boundary interval is bounded by that of the corresponding enlarged
interval at the unshifted point. The conditional moment of \(y_i\)
can then be taken separately.

The coefficients may select prefixes close to the boundary, so
the bound on probabilities near endpoints must remain valid after that
selection. For each nonnegative predictable coefficient \(a_i\), write its sum as
an integral over levels of its cumulative sum. The first index crossing
a fixed level is predictable, and \eqref{poly:eq:pin-bound}
applies at this selected prefix.
At that index the Gaussian vector in the derivative has covariance
\(I_d-P_i+tC_i\), \(0\le t\le1\). Adding an independent conditional
normal of covariance \((1-t)C_i\) restores \(I_d-P_{i-1}\).
On at most four coordinates this added normal is within \(Ch\)
with conditional probability at least \(1/2\), because
\(\|C_i\|_{\max}\le Ch_0^2\). Thus
\eqref{poly:eq:pin-bound} bounds the probability of each small box
around the selected endpoints, uniformly at the crossing index.
A level never reached contributes zero; assigning it
the dummy index \(n+1\) gives the same upper bound.
Integrating the levels and using
\eqref{poly:eq:identity-moment-budgets}, the smooth discrepancy is
at most
\begin{equation}\label{poly:eq:identity-smooth-error}
 C_q\sum_{\ell=2}^4 h_0^\ell
                  \sum_{k=1}^{\ell}h^{k-\ell}W_A^k .
\end{equation}
By \eqref{poly:eq:identity-weight}, a summand is bounded by
\[
 C_q\varepsilon^k(1+w)^{2k-\ell}.
\]
For \(k=1\) its last exponent is nonpositive. For \(k\ge2\),
\(w\le4\log(e/\varepsilon)\) makes it \(O_q(\varepsilon)\).

\paragraph{Removing the compact smoothing.}
The compact support sandwiches \(A\) between its smooth expansions
and contractions by \(h\), using the same bandwidth. Since
\(\log\omega\) is \(T\)-Lipschitz, endpoint shifts of size \(Ch\)
preserve \eqref{poly:eq:identity-weight}. For independent Gaussian
coordinates, the derivative of a rectangle probability under a
simultaneous endpoint displacement is at most
\[
 C W_{A'}\exp(1-w')
\]
at an intermediate rectangle \(A'\), where \(w'\) is the sum of
its one-coordinate Gaussian exit probabilities. This follows by
differentiating the product of its one-dimensional probabilities.
The factors other than the
differentiated coordinate have product at most \(e^{1-w'}\).
The hazard bound on \(|b|\le T\) gives \(w'\ge cw-C\).
The outward tails beyond \(T\) have total at most
\(2d\bar\Phi(T)\le1\); the opposite tails remain at least \(1/2\).
Thus the Gaussian probability of the boundary region is at most
\[
 ChT(1+w)e^{-cw}\le C\varepsilon.
\]
If a contraction becomes empty, stop the integration at its first
collapse. Together with \eqref{poly:eq:identity-smooth-error},
this proves the asserted comparison for \(W_0\) on these rectangles.

Except on an event of probability \(2\varepsilon\), all original
increments are retained and untruncated. Then \(P_n=I_d\), the appended
normal is zero, and \(\|S_n-W_0\|_\infty\le h_0\).
One more expansion and contraction gives the same bound for \(S_n\).
Finally, if \(w>4\log(e/\varepsilon)\), retain whole coordinate
intervals until the sum of their Gaussian exit probabilities first
exceeds \(2\log(e/\varepsilon)\). For the resulting containing
cylinder, this sum is at most \(2\log(e/\varepsilon)+1\), and its Gaussian
probability is at most \(\varepsilon^2\).
The bound already proved for this cylinder bounds both probabilities
of the original rectangle by \(C_q\varepsilon\).
\end{proof}

\subsection{A joint bound for random predictable quadratic variation}

\begin{proposition}\label{poly:prop:joint}
Suppose \(2<q\le3\),
\(\max_j\sum_k|\Sigma_{jk}|^s\le M\), and
\(\|\|Q-\Sigma\|_{\max}\|_{L^p}\le\delta\), where
\(1\le p\le\infty\). Then
\begin{align}
 \rho_{\cR}(S_n,G_\Sigma)
 \le \min\bigg\{1,\,
 C_{q,s,M}\{\log(en)\}^{2+3\max\{2,s\}}
 \bigg[
 &(v_q\{\log(ed)\}^{q/2})^{1/(q+1)}
 \notag\\[-2pt]
 &+(\delta\log(ed))^{1/(2+1/p)}
 \bigg]\bigg\}.
 \label{poly:eq:joint-upper}
\end{align}
\end{proposition}

\begin{proof}
Let
\[
 T=D\sqrt{\log(2nd)},\qquad r=\lceil100\log(en)\rceil,\qquad
 \eta=c_0h^2/r,
\]
where \(h>0\) will be selected below. The constants \(D\) and \(c_0\)
are sufficiently large and sufficiently small, respectively.
We first work with \(hT\) smaller than a fixed constant.

\paragraph{One localization for the whole vector.}
Let
\(\mathcal K_\eta=\{A\succeq0:\|A-\Sigma\|_{\max}\le\eta\}\).
Retain increments until the next increment would violate either
\[
 \sum_i I_i\E_{i-1}\|\xi_i\|_\infty^q\le c_0h^q/r
\]
or domination of the retained original prefix covariance by a
member of \(\mathcal K_\eta\). Here the test is made before retaining
the increment, so \(I_i\) is predictable. For the final retained covariance
\[
 P=\sum_i I_i\E_{i-1}\xi_i\xi_i^\top,
\]
choose \(A(P)\in\mathcal K_\eta\) with \(P\preceq A(P)\), taking
\(A(P)=P\) when \(P\in\mathcal K_\eta\).
Lemma~\ref{thr:lem:feasible-completion} provides the measurable
selection.

Truncate the retained increments at \(h\) and center them conditionally:
\[
 y_i=I_i\xi_i\ind_{\{\|\xi_i\|_\infty\le h\}}
       -\E_{i-1}[I_i\xi_i\ind_{\{\|\xi_i\|_\infty\le h\}}].
\]
Append a conditional normal with covariance \(A(P)-P\), and call the
resulting sum \(W\), with predictable quadratic variation \(Q'\).
The same positive semidefinite covariance loss calculation as above
gives
\begin{equation}\label{poly:eq:joint-pathwise}
 \|Q'-\Sigma\|_{\max}\le\eta':=3c_0h^2/r
 \quad\text{almost surely}.
\end{equation}
If the original \(Q\) belongs to \(\mathcal K_\eta\), the total
sum of conditional \(q\)th moments does not exceed \(c_0h^q/r\),
and no increment exceeds \(h\) in maximum norm, every increment is
retained and the appended normal is zero.
On this event
\(\|S_n-W\|_\infty\le c_0h/r\).
The probability of its complement is at most
\begin{equation}\label{poly:eq:joint-exception}
 Crv_qh^{-q}+\Pp(\|Q-\Sigma\|_{\max}>c_0h^2/r).
\end{equation}
This is a single exceptional event for the vector.

\paragraph{The signed covariance correction.}
For a product-logistic smoothing \(f\) of a rectangle, use the finite
correction \(\mathcal U_r(U,E;x)\) of
Lemma~\ref{quant:lem:finite-correction}.
Since \(r\eta'/h^2=3c_0\), its absolute derivative sums through
order four are bounded by \(C_jh^{-j}\).
If \(P_i'\) is the sum of conditional covariance matrices through
step \(i\) of the modified array,
use the post-step pair
\[
 U_i=\E_{i-1}Q'-P_i'\succeq0,\qquad
 E_i=\Sigma-\E_{i-1}Q'.
\]
The pre-step pair is \((U_i+C_i',E_i)\), where \(C_i'\) is
the conditional covariance of that increment. The sum of this pair is
\(\Sigma-P_{i-1}'\), equal to the sum of the preceding post-step
pair. Thus the positive Gaussian bases may change while the
total covariance remains fixed. Conditional Jensen and
\eqref{poly:eq:joint-pathwise} give
\(\|E_i\|_{\max}\le\eta'\).
The switching bound \eqref{quant:eq:heat-switch}, including the
initial and terminal switches, therefore contributes \(O(n^{-10})\).

For each truncated martingale increment, Taylor expansion to order two has
remainder bounded both by
\(Ch^{-2}\|y_i\|_\infty^2\) and by
\(Ch^{-3}\|y_i\|_\infty^3\). Their minimum is at most
\(C_qh^{-q}\|y_i\|_\infty^q\).
For the corresponding Gaussian heat step, the remainder after the
linear covariance term is bounded by
\[
 C\min\{h^{-2}\|C_i'\|_{\max},
                 h^{-4}\|C_i'\|_{\max}^2\}
 \le C_qh^{-q}\|C_i'\|_{\max}^{q/2}
 \le C_qh^{-q}\E_{i-1}\|y_i\|_\infty^q.
\]
The first two moments agree. Conditional centering costs at most a
factor \(2^q\) in the moment sum. The appended Gaussian increment is
treated exactly by the heat operator. Thus
\begin{equation}\label{poly:eq:joint-corrected}
 \left|\E\mathcal U_r(0,\Sigma-Q';W)-\E f(G_\Sigma)\right|
 \le C_q(v_qh^{-q}+n^{-10}).
\end{equation}
This comparison uses absolute derivatives of a signed function; no
positivity is asserted for \(\mathcal U_r\).

\paragraph{Scalar strips for the terminal correction.}
For a coordinate \(j\), append a further conditional normal of
variance \(1+\eta'-Q'_{jj}\in[0,2\eta']\).
Its sum with \(W_j\), denoted by \(W_j^*\), has deterministic
predictable quadratic variation \(1+\eta'\). The truncated increments
satisfy
\[
 \|y_i\|_\infty\le2h,\qquad
 \sum_i\E_{i-1}\|y_i\|_\infty^3\le C_qc_0h^3/r.
\]
Both terminal Gaussians combine, conditionally on the original
terminal history, into the exact remaining variance.
Apply Lemma~\ref{poly:lem:positive-comparison} in dimension one at
bandwidth \(C_1h\), for a sufficiently large fixed \(C_1\), using
the actual conditional covariances of the truncated increments.
It follows that the
density of \(W_j^*+C_1hH\) is bounded by
\[
 C(\phi_{1+\eta'}*p_{C_1h})(x),
\]
where \(H\) is standard logistic and \(\phi_v\) is the
\(N(0,v)\) density. One can obtain this directly by taking
\(f(x)=p_{C_1h}(x-b)\) in the lemma.

Conditionally on the terminal history and \(W_j\), the additional
normal and logistic displacement is at most \(C_2h\) with probability
at least \(1/2\). Since \(\eta'T^2\le C(hT)^2/r\), integration of the
last density bound gives, for \(|b|\le2T\) and \(t\ge0\),
\begin{align}
 \Pp(|W_j-b|\le t)&\le C(t+h)\phi(b)e^{CTt},\notag\\
 \Pp(\pm W_j>T-h)&\le Ce^{-T^2/2}.
 \label{poly:eq:joint-strip}
\end{align}
For the first inequality use
\(\phi_{1+\eta'}(b+z)
 \le C\phi(b)e^{CT|z|}\) and the logistic exponential moment.
For the second, apply exponential Markov at parameter \(T\) to
\(N(0,1+\eta')+C_1hH\), at threshold \(T-C_3h\).

By \eqref{quant:eq:terminal-correction}, the difference between the
first expectation in \eqref{poly:eq:joint-corrected} and \(\E f(W)\)
is at most
\[
 C\sum_j\E\min\{f_j(W_j),1-f_j(W_j)\}.
\]
For each interval,
\(\min(f_j(x),1-f_j(x))\le
C\sum_{b\in\partial I_j}e^{-|x-b|/h}\).
Integrating \eqref{poly:eq:joint-strip} on exponential shells
bounds this cost by \(Ch\sum_b\phi(b)\).
The same sum controls unsmoothing. Indeed a crossing of face \(b\)
after adding \(hH_j\) requires
\(|W_j-b|\le h|H_j|\), and hence its probability is at most
\[
 Ch\phi(b)\E[(1+|H_j|)e^{ChT|H_j|}]
 \le Ch\phi(b).
\]
The Gaussian-side crossing satisfies the same estimate. On the
common good event in \eqref{poly:eq:joint-exception}, the displacement
from \(W\) to \(S_n\) is \(O(h/r)\), so another use of the strip
bound suffices.

\paragraph{Gaussian geometry and optimization.}
Intersect the rectangle with \([-T,T]^d\). By the second bound in
\eqref{poly:eq:joint-strip} and \eqref{poly:eq:joint-exception}, this
changes its probability by at most the common exceptional probability
plus \(Cn^{-8}\). For the Gaussian vector the change is at most
\(Cn^{-8}\).
Proposition~\ref{thr:prop:quantitative-cylinder} either retains the
restricted rectangle or replaces it by a containing cylinder of
Gaussian probability \(O(n^{-8})\). In both cases
\[
 \sum_b\phi(b)
 \le C\{\log(en)\}^u\sqrt{\log(ed)}
 \le C\{\log(en)\}^uT,
\]
where \(u=\max\{2,s\}\). If a cylinder is selected, the bound
on its discrepancy controls both probabilities of the original
rectangle. We have proved
\begin{align}
 \rho_{\cR}(S_n,G_\Sigma)
 \le C\{\log(en)\}^u hT
       +Crv_qh^{-q}
       +\Pp(\|Q-\Sigma\|_{\max}>c_0h^2/r)+Cn^{-8}.
 \label{poly:eq:joint-before-optimization}
\end{align}

Put \(\gamma=(2+1/p)^{-1}\), with \(1/\infty=0\), and choose
\begin{equation}\label{poly:eq:joint-bandwidth}
 h=\max\left\{(rv_q/T)^{1/(q+1)},
                   (r\delta/c_0)^\gamma T^{2\gamma-1}\right\}.
\end{equation}
The first term makes \(rv_qh^{-q}\le hT\). For finite \(p\),
Markov's inequality and
\[
 \big\{(r\delta/c_0)^\gamma T^{2\gamma-1}\big\}^{2p+1}
 =(r\delta/c_0)^p/T
\]
bound the probability involving \(Q-\Sigma\) by \(hT\).
For \(p=\infty\), the second candidate is
\(\sqrt{r\delta/c_0}\), so the strict exceedance probability is zero.
This also covers \(\delta=0\).
Unit marginal variances imply
\(v_q\ge n^{1-q/2}>0\), so the bandwidth is positive.
If \(hT\) is not small, the asserted bound follows by its unit cap.
In the remaining case these same inequalities imply
\(\log(2nd)\le Cn\), as required by the cylinder estimate.

Finally, \(r\le C\log(en)\) and
\(T^2\le C\log(en)\log(ed)\). Substitution of
\eqref{poly:eq:joint-bandwidth} in
\eqref{poly:eq:joint-before-optimization} gives
\eqref{poly:eq:joint-upper}; the displayed power of \(\log(en)\)
is at least as large as either power obtained by substitution. The term
\(n^{-8}\) is absorbed by the term involving \(v_q\), since \(2<q\le3\).
The separate powers are recorded in Appendix~\ref{app:logs}.
\end{proof}

\subsection{An adaptive lower bound with deterministic conditional variances}
\label{poly:subsec:lower}

The scalar construction that changes the increment law near a boundary
appears in
\citet[Example~2, pp.~148--149]{Joos1993}. A large symmetric step
is selected near a boundary, and the remaining variance is completed
on the complementary event. The construction below selects one of many
coordinates and distributes the exceptional jump over the remaining
steps. This preserves deterministic one-step conditional covariances
while bounding \(\E\|X_i\|_\infty^q\) uniformly in \(i,n,d\).

\begin{proposition}\label{poly:prop:identity-lower}
Fix \(2<q\le3\). For every sufficiently large \(n\) and every
\(d\ge1\), there is a martingale difference array such that
\[
 Q=I_d,\qquad
 \max_i\E\|X_i\|_\infty^q\le2^q,
\]
every one-step conditional covariance is deterministic and diagonal,
and every conditional increment law is centrally symmetric, while
\begin{equation}\label{poly:eq:identity-lower}
 \rho_{\cR}(S_n,G_{I_d})
 \ge c_q\min\left\{1,
   \bigl(n^{1-q/2}\{\log(ed)\}^{q/2}\bigr)^{1/(q+1)}\right\}.
\end{equation}
\end{proposition}

\begin{proof}
We construct two arrays with the same Gaussian target and compare
their probabilities of one rectangle. A difference between these
probabilities forces at least one Gaussian approximation error
to be half as large. The late increments will have the same
deterministic conditional variances in both arrays; on a selected
prefix event, only one coordinate's conditional law changes. Let
\(n_0=\lfloor n/2\rfloor\), \(n_1=n-n_0\).
First suppose \(d\) exceeds a sufficiently large fixed constant and
\(\log(ed)\le n^{1-2/q}\). Put
\[
 x=\Phi^{-1}(1-1/d),\qquad
 h=c_q n^{-(q-2)/(2(q+1))}
                    \{\log(ed)\}^{-1/(2(q+1))},\qquad
 \lambda=1/16 .
\]
Choose \(c_q>0\) small. Then
\begin{equation}\label{poly:eq:lower-scales}
 h\sqrt{\log(ed)}\le c_q,\qquad
 h\sqrt n\ge c_q n^{1/q}\longrightarrow\infty,\qquad
 \frac{\{\log(ed)\}^2}{n}\longrightarrow0
\end{equation}
uniformly in the indicated range.

In normalized units, use independent signs of magnitude
\[
 a_0=\sqrt{(1-h^2)/n_0}
\]
in the first \(n_0\) increments, independently across coordinates.
Write \(U_j\) for the prefix in coordinate \(j\).
In the reference array complete every coordinate with \(n_1\)
independent signs of magnitude \(h/\sqrt{n_1}\); let its sum be
\(S^0\). In the second array, define the disjoint prefix events
\[
 E_j=\{x+h/16\le U_j\le x+h/8,\quad U_k\le x\ (k\ne j)\}.
\]
On \(E_j\), replace coordinate \(j\) in every remaining increment by
\[
 \begin{cases}
  h/\sqrt\lambda,&\text{with probability }\lambda/(2n_1),\\
  -h/\sqrt\lambda,&\text{with probability }\lambda/(2n_1),\\
  0,&\text{otherwise},
 \end{cases}
\]
using fresh independent draws. All other coordinates retain their
reference completion. Outside the union of the \(E_j\), the arrays
agree.

Both completions have conditional variance \(h^2/n_1\) in every
coordinate at every remaining step. Their cross-covariances are zero.
Thus each array has deterministic diagonal one-step conditional
covariances and \(Q=I_d\). Conditional on each history, its coordinate
laws are independent and symmetric, so the whole conditional law
is centrally symmetric. The choice of which coordinate law to use
is nevertheless adapted to the prefix.

\paragraph{The boundary band.}
Stirling's formula for the symmetric binomial distribution gives,
at its lattice points \(z\) with \(|z|\le C\sqrt{\log(ed)}\),
\[
 \Pp(U_j=z)
 =\frac{2a_0}{\sqrt{1-h^2}}\,
   \phi\!\left(\frac z{\sqrt{1-h^2}}\right)
   \left\{1+O\left(\frac{1+(z/\sqrt{1-h^2})^4}{n_0}\right)\right\}.
\]
This follows by expanding the binomial entropy through its quadratic
term; the next term is \(O(z^4/n_0)\), and the Stirling prefactor
has relative error \(O((1+z^2)/n_0)\).
By \eqref{poly:eq:lower-scales} the error tends to zero uniformly.
The band contains order \(h\sqrt n\) lattice points. Its normal
density is within fixed factors of \(\phi(x)\), because
\(hx+h^2x^2\le C_q\). Mills' inequalities therefore give
\begin{equation}\label{poly:eq:lower-band}
 \Pp(x+h/16\le U_j\le x+h/8)
 \asymp_q h\phi(x)
 \asymp_q h\sqrt{\log(ed)}/d.
\end{equation}

We also need the uniform tail bounds
\begin{equation}\label{poly:eq:lower-tails}
 \Pp(U_j>x)+\Pp(S_j^0>x)\le C/d .
\end{equation}
For either of these sign sums \(V\), tilt its independent signs by
\(e^{xV}/\E e^{xV}\). Its first \(n_0\) tilted signs have a common
success probability in \([1/4,3/4]\), because \(xa_0=o(1)\).
Stirling's formula bounds their largest tilted point mass by
\(C/\sqrt n\). Conditional on the remaining signs, their probability
of an interval of length \(1/x\) is at most
\[
 Cn^{-1/2}\{1+(xa_0)^{-1}\}\le C/x.
\]
Since \(\E e^{xV}\le e^{x^2/2}\), summing the intervals
\([x+\ell/x,x+(\ell+1)/x)\), with their inverse-tilt weights
at most \(e^{-\ell}\), yields
\[
 \Pp(V>x)\le C e^{-x^2/2}/x\le C'/d.
\]
This proves \eqref{poly:eq:lower-tails}.

\paragraph{The difference between the two rectangle probabilities.}
On the selected band, the reference suffix has probability at least
\(0.4\) of finishing coordinate \(j\) at or below \(x\), for all
large \(n\). Indeed its standardized value is a normalized sum of
\(n_1\) signs and the required threshold lies between \(-1/8\) and
\(-1/16\); summing the same central binomial estimate gives this
bound. In the alternative array, finishing below \(x\) requires at
least one nonzero late jump. Its probability is at most \(\lambda\).
After integrating the unaffected coordinates, the difference of
the two probabilities of \((-\infty,x]^d\) is at least
\[
 (0.4-\lambda)d\,
 \Pp(x+h/16\le U_1\le x+h/8)\,
 \{\Pp(U_1\le x,S_1^0\le x)\}^{d-1}.
\]
For all sufficiently large fixed \(d\) the last factor is bounded
below uniformly by \eqref{poly:eq:lower-tails}. Thus the difference
is at least \(c_qh\sqrt{\log(ed)}\).
The two arrays have the same Gaussian target, so one of their
Gaussian approximation errors is at least half this quantity.

\paragraph{The moment bound and the remaining dimensions.}
Return to the unnormalized increments \(X_i=\sqrt n\,\xi_i\).
An initial increment has maximum norm \(\sqrt n\,a_0<2\), for large
\(n\). In each remaining increment, the ordinary signs have magnitude
at most \(2h\), while
there is at most one active exceptional coordinate. Hence
\[
 \E\|X_i\|_\infty^q
 \le (2h)^q+
 C_q n^{q/2-1}h^q\,\Pp\Big(\bigcup_j E_j\Big)
 \le (2h)^q+C_q n^{q/2-1}h^{q+1}\sqrt{\log(ed)}
 \le2^q
\]
after decreasing \(c_q\). The last expression is bounded because
\[
 n^{q/2-1}h^{q+1}\sqrt{\log(ed)}=c_q^{q+1}.
\]
The reference array satisfies the same bound.

For \(d\) in a fixed finite range, use one coordinate, threshold zero,
and the prefix band \([h/10,h/5]\), with
\(h=c_q n^{-(q-2)/(2(q+1))}\). Its probability is comparable to
\(h\), and the same completion comparison gives a difference at
least \(c_qh\). Add independent Rademacher coordinates in the
unnormalized array and test a cylinder. This preserves the moment
bound, since the estimate for each remaining increment becomes \(1+C_qc_q^{q+1}\le2^q\).

If \(\log(ed)>n^{1-2/q}\), perform the large-dimension construction
on
\(\lfloor\exp\{n^{1-2/q}/2\}\rfloor\) coordinates and pad the
remaining coordinates in the same manner. Its discrepancy is
bounded below by a positive constant. Projection preserves the
rectangle discrepancy and all the stated covariance properties.
This proves \eqref{poly:eq:identity-lower} for every \(d\ge1\).
\end{proof}

\begin{remark}\label{poly:rem:scalar-scale}
At \(d=1\), the rate in
Propositions~\ref{poly:prop:identity-upper}--\ref{poly:prop:identity-lower}
is \(n^{-(q-2)/(2(q+1))}\). This exponent agrees with the scalar
moment term in Haeusler's inequality \citep{haeusler1988rate},
restated in \citet[Theorem~1.1]{mourrat2013rate}: take its moment
parameter equal to \(q/2\). When the predictable quadratic variation
is one, that term is \(v_q^{1/(q+1)}\), and the uniform moment bound gives
\(v_q\le B^qn^{1-q/2}\). Scalar sharpness with predictable quadratic
variation equal to one already appears in \citet[Example~2]{Joos1993}.
For \(q=3\) the
rate is \(n^{-1/8}\).
There is no conflict with results under a variance-relative
conditional moment bound. In the active suffix of the scalar
construction,
\[
 \E_{i-1}X_i^2=\frac{nh^2}{n_1},\qquad
 \E_{i-1}|X_i|^3=\frac{n^{3/2}h^3}{\sqrt\lambda\,n_1},
 \qquad
 \frac{\E_{i-1}|X_i|^3}{\E_{i-1}X_i^2}
       =\frac{\sqrt n\,h}{\sqrt\lambda}.
\]
At \(q=3\), \(h\asymp n^{-1/8}\); the conditional variance of each
remaining increment tends to zero, and the last ratio diverges.
The active event has probability
of order \(h\), which keeps the marginal third moment bounded.
Deterministic conditional variances give neither a positive lower bound
on the individual variances nor independence of the increment laws
from the past.
\end{remark}

\section{Matching lower bounds}\label{low:section}

Proposition~\ref{poly:prop:identity-lower} gives the finite-moment lower
bound when \(Q=I_d\). We now prove the corresponding lower bounds for
bounded increments and marginal light tails, followed by the lower
bound involving \(Q-\Sigma\) common to all classes. Every construction
has covariance \(I_d\). The examples that change the conditional law
near a boundary develop the scalar arguments of
\citet{bolthausen1982exact} and \citet{mourrat2013rate}.
We verify the increment assumptions and the identity \(\E Q=I_d\)
for each example. Section~\ref{main:proofs} then combines the lower bounds.

\subsection{Lower bounds for bounded increments and marginal light tails}

Write \(\bar\Phi=1-\Phi\), and let \(\phi\) be the standard normal
density. We first record a scalar estimate with an explicit support
parameter.

\begin{lemma}\label{low:relative}
Let \(Y_1,\ldots,Y_N\) be independent copies of a centered
unit-variance random variable \(Y\) with \(|Y|\le b\), \(b\ge1\).
There are absolute constants \(c,C>0\) such that, when
\(x\ge0\) and \(b(1+x)/\sqrt N\le c\),
\begin{equation}\label{low:eq:relative}
 \log\frac{\Pp(N^{-1/2}\sum_{i=1}^NY_i>x)}{\bar\Phi(x)}
 =
 \frac{\E Y^3\,x^3}{6\sqrt N}
 +O\left(\frac{b^2x^4}{N}+\frac{b(1+x)}{\sqrt N}\right).
\end{equation}
The remainder has absolute value at most \(C\) times the displayed
quantity. The same assertion holds for the non-strict tail.
\end{lemma}

\begin{proof}
Apply \citet[Lemma~F.1]{CCKK2026} to the independent variables
\(\eta_i=Y_i/\sqrt N\). Its variance sum is one, its support
parameter is \(b/\sqrt N\), and its third-cumulant sum is
\(\E Y^3/\sqrt N\). The condition in that lemma is precisely
\(b(1+x)/\sqrt N\le c\), and its strict-tail conclusion is
\eqref{low:eq:relative}, with absolute constants.
For the non-strict tail at \(x>0\), let strict thresholds increase
to \(x\) and use continuity of the normal tail and the displayed
remainder bound. At \(x=0\), apply the strict-tail result to
\(-Y\) and take complements; both probabilities are
\(1/2+O(b/\sqrt N)\), which gives the same logarithmic bound.
\end{proof}

\begin{proposition}\label{low:cubic}
There is an absolute \(c>0\) such that, for all \(n\ge2,d\ge1\),
the class with \(\|X_i\|_\infty\le2\) and \(Q=I_d\) satisfies
\begin{equation}\label{low:eq:cubic}
 D^{I,\mathrm{bd}}_{n,d}
 \ge c\left[1\wedge n^{-1/2}\{\log(ed)\}^{3/2}\right].
\end{equation}
For \(\alpha\in\{1,2\}\), the marginal class with scale \(K=4\)
and \(Q=I_d\) satisfies
\begin{equation}\label{low:eq:marginal}
 D^{I,\mathrm{marg},\alpha}_{n,d}
 \ge c\left[1\wedge
 n^{-1/2}\{\log(ed)\}^{3/2+1/\alpha}\right].
\end{equation}
The examples can have independent coordinate processes, trivial
initial information and identity conditional covariance at every step.
\end{proposition}

\begin{proof}
For the bounded example take independent entries with
\[
 \Pp(X=-2)=\frac15,\qquad \Pp(X=1/2)=\frac45.
\]
Then \(\E X=0\), \(\E X^2=1\) and \(\E X^3=-3/2\).
Let \(d\to\infty\) and \(u=\Phi^{-1}(1-1/d)\). If
\(u^3/\sqrt n\to0\), Lemma~\ref{low:relative} gives
\[
 \log\frac{\Pp(S_{n,1}>u)}{\bar\Phi(u)}
 =
 -\frac{u^3}{4\sqrt n}
 +O\left(\frac{u^4}{n}+\frac{1+u}{\sqrt n}\right)
 =
 -\frac{u^3}{4\sqrt n}+o(u^3/\sqrt n).
\]
The derivative of \(z\mapsto\{1-(1+z)/d\}^d\) tends to
\(-e^{-1}\), uniformly for \(z=o(1)\). The rectangle
\((-\infty,u]^d\) therefore has error at least
\((e^{-1}/4+o(1))u^3/\sqrt n\).

For the marginal example let \(m=\lfloor n/2\rfloor\), \(k=n-m\),
\(a=\sqrt{m/n}\), \(c_1=\sqrt{k/n}\), and
\[
 b=(u/8)^{2/\alpha}.
\]
It suffices initially to take \(u\ge8\). The prefix consists of
\(m\) independent standard normals, with standardized sum \(Z\).
On the prefix event
\[
 A=\{Z>u/(2\sqrt2)\},
\]
use \(k\) fresh independent suffix variables with law
\[
 \Pp(Y=b)=\frac1{1+b^2},\qquad
 \Pp(Y=-b^{-1})=\frac{b^2}{1+b^2}.
\]
Off \(A\), use fresh standard normals. Both suffix laws are centered
with variance one, and \(\E Y^3=b-b^{-1}\).
Independent copies of the entire scalar construction define the
coordinate processes.
Since
\[
 \Pp(A)\le e^{-u^2/16}=e^{-4b^\alpha},
\]
the marginal exponential moment at scale four is bounded by
\[
 \E e^{(|N(0,1)|/4)^\alpha}
 +e^{-(4-4^{-\alpha})b^\alpha}<2
\]
for all sufficiently large \(u\). Thus the stated marginal condition
holds at every step and coordinate.

Take a fresh standard normal \(H\), put \(T=aZ+c_1H\), and write
\[
 x=\frac{u-aZ}{c_1},\qquad
 R(x)=\frac{\Pp(k^{-1/2}\sum_{i=1}^kY_i>x)}{\bar\Phi(x)}.
\]
Conditioning on the prefix gives the exact identity
\begin{equation}\label{low:eq:inactive-cancellation}
 \frac{\Pp(S_{n,1}>u)}{\bar\Phi(u)}-1
 =
 \E\bigl[\ind_A\{R(x)-1\}\mid T>u\bigr].
\end{equation}
On the complement of \(A\), the suffix is Gaussian in both arrays,
so those histories cancel in this identity.

Suppose \(bu^3/\sqrt n\to0\). On \(A\cap\{x\ge0\}\),
one has \(x\le Cu\), so Lemma~\ref{low:relative} gives uniformly
\[
 R(x)-1
 =
 \frac{(b-b^{-1})x^3}{6\sqrt k}
 +O\left(\frac{b^2u^4}{n}
       +\frac{b(1+u)}{\sqrt n}
       +\frac{b^2u^6}{n}\right).
\]
The remainder is \(o(bu^3/\sqrt n)\).
For \(x<0\), the scalar Berry--Esseen inequality and
\(\bar\Phi(x)\ge1/2\) give
\(|R(x)-1|\le Cb/\sqrt k=o(bu^3/\sqrt n)\).
The Gaussian representation
\[
 Z=aT+c_1V,
\]
where \(V\sim N(0,1)\) is independent of \(T\), shows that,
conditionally on \(T>u\), \(T/u\to1\), \(\Pp(A\mid T>u)\to1\),
and \(x/u\to c_1\) in probability. On \(A\), \(x_+/u\) is
bounded. Hence
\[
 \E[\ind_Ax_+^3\mid T>u]\sim c_1^3u^3.
\]
Substitution in \eqref{low:eq:inactive-cancellation} yields
\[
 \frac{\Pp(S_{n,1}>u)}{\bar\Phi(u)}-1
 \sim\frac{bu^3}{12\sqrt n}.
\]
Taking products over coordinates gives the lower order in
\eqref{low:eq:marginal} in the subcritical range.

We next obtain the lower bounds at a fixed critical scale.
For the bounded construction, suppose that
\(u^3/\sqrt n\to\lambda\in(0,\infty)\).
The remainder in Lemma~\ref{low:relative} tends to zero, since
\(u^4/n=O(u^{-2})\) and \((1+u)/\sqrt n=O(u^{-2})\).
Consequently
\[
 \frac{\Pp(S_{n,1}>u)}{\bar\Phi(u)}
       \longrightarrow e^{-\lambda/4}.
\]

For the marginal construction, suppose instead that
\(bu^3/\sqrt n\to\lambda\in(0,\infty)\).
On \(A\cap\{x\ge0\}\), we have \(x\le Cu\), and
Lemma~\ref{low:relative} gives
\[
 \sup_{0\le x\le Cu}
 \left|\log R(x)-\frac{(b-b^{-1})x^3}{6\sqrt k}\right|
 \le C\left\{\frac{b^2u^4}{n}
              +\frac{b(1+u)}{\sqrt n}\right\}=o(1).
\]
Conditionally on \(T>u\), the representation
\(Z=aT+c_1V\) gives \(\Pp(A\mid T>u)\to1\) and
\(x/u\to1/\sqrt2\) in probability. Since \(k/n\to1/2\)
and \(b\to\infty\), the cubic term in the last display tends
to \(\lambda/12\) in this conditional probability.
Moreover, \(R(x)\) is uniformly bounded on
\(A\cap\{x\ge0\}\), because its logarithm is at most
\(Cbu^3/\sqrt n+o(1)\).
On \(\{x<0\}\), Berry--Esseen gives
\[
 |R(x)-1|\le Cb/\sqrt k=o(1).
\]
These bounds justify passage to the conditional expectation in
\eqref{low:eq:inactive-cancellation}. Hence
\[
 \frac{\Pp(S_{n,1}>u)}{\bar\Phi(u)}
       \longrightarrow e^{\lambda/12}.
\]
The cancellation on the complement of \(A\) therefore also applies
at the critical scale.

Since \(d\bar\Phi(u)=1\), the Gaussian rectangle probability
tends to \(e^{-1}\). The bounded and marginal rectangle
probabilities tend, respectively, to
\[
 \exp\{-e^{-\lambda/4}\}
 \quad\text{and}\quad
 \exp\{-e^{\lambda/12}\}.
\]
Each limit differs from \(e^{-1}\).
If the claimed lower bound is bounded away from zero, choose a smaller
number of coordinates so that the corresponding critical
parameter tends to a sufficiently small fixed \(\lambda>0\).
Add independent Rademacher coordinates and test the cylinder on the
selected coordinates. This preserves the increment assumptions
and the identity conditional covariance at every step, and yields
a rectangle error bounded away from zero.

Finally, these constructions give constants uniform in \(n,d\).
Otherwise there would be a sequence along which the ratio of the
maximal rectangle error to the claimed lower-bound expression tends
to zero. Bounded \(n\) is impossible: a scalar Rademacher sum has
an atom of mass at least \(c/\sqrt n\), and its distance from a
continuous Gaussian distribution is at least half that mass.
After passing to a subsequence we may therefore assume
\(n\to\infty\). Bounded dimensions are ruled out by the same
atom bound, since the dimension factors are then bounded.
For \(d\to\infty\), if the lower-bound expression tends to zero,
the relative expansions below the critical scale apply. If it is
bounded away from zero, use the critical construction and add independent
Rademacher coordinates. Passing to one of these two subsequences
exhausts the possibilities and proves both uniform lower bounds.
\end{proof}

\subsection{A lower bound for the quadratic-variation error}

\begin{proposition}\label{low:bracket}
For \(n\ge2\), \(d\ge1\), \(\delta\ge0\) and \(1\le p\le\infty\),
the class of arrays with
\[
 \|X_i\|_\infty\le2,\qquad
 \E Q=I_d,\qquad
 \big\|\|Q-I_d\|_{\max}\big\|_p\le\delta
\]
has maximal rectangle error at least
\begin{equation}\label{low:eq:bracket}
 c\left[1\wedge\left\{
 n^{-1/2}\{\log(ed)\}^{3/2}
 +\{\delta\log(ed)\}^{1/(2+1/p)}
 \right\}\right].
\end{equation}
The conditional covariance matrices can be diagonal, the coordinate
processes independent, and the initial sigma-field trivial.
\end{proposition}

\begin{proof}
The first term is Proposition~\ref{low:cubic}. For the second
we use a boundary-triggered variance change, as in
\citet[Example~3, pp.~149--150]{Joos1993}; see also
\citet{mourrat2013rate} for sharpness of the scalar covariance
exponent. We give a bounded-increment construction with the required
dimension factor and normalization \(\E Q=I_d\), at scales not already
covered by the first term.
Let \(m=\lfloor n/2\rfloor\), \(k=n-m\), and, for each coordinate, take
independent Rademachers with standardized sums
\[
 P=m^{-1/2}\sum_{i=1}^m\varepsilon_i,\qquad
 Z=k^{-1/2}\sum_{i=1}^k\varepsilon'_i.
\]
The baseline array has terminal sum \(P\). For \(0<h\le h_*\),
where \(h_*\) is a sufficiently small absolute constant, put
\[
 u=\Phi^{-1}(1-1/(2d)),\quad
 v=\frac{h^2}{\log(ed)},\quad
 A=\{u-\sqrt v/2\le P\le u\},\quad
 \pi=\Pp(A).
\]
The adaptive terminal sum is
\[
 \frac{P+\ind_A\sqrt v\,Z}{\sqrt{1+\pi v}}.
\]
The first \(m\) increments give the prefix; the remaining increments
give the suffix when \(A\) occurs. This is a martingale array, since
\(A\) is determined before the suffix and all fresh signs are independent.
The unnormalized increments are bounded by two, because \(n/m\le3\),
\(n/k\le2\) and \(v\le1\).
Its predictable quadratic variation is diagonal, with entries
\[
 Q_{jj}=\frac{1+v\ind_{A_j}}{1+\pi v},
\]
and therefore \(\E Q=I_d\).

Suppose first that
\(\log(ed)\le n^{1/3}\) and
\(h\ge C_*\sqrt{\log(ed)/n}\).
Uniformly in this range,
\begin{equation}\label{low:eq:binomial-band}
 \pi\asymp\sqrt v\,\phi(u)\asymp h/d,\qquad
 \Pp(P>u)\le C/d.
\end{equation}
Here \(1+u\asymp\sqrt{\log(ed)}\) and
\(\phi(u)\asymp\sqrt{\log(ed)}/d\), including \(d=1\).
The lattice spacing of \(P\) is \(2/\sqrt m\); the specified lower
bound on \(h\) makes the strip contain order \(\sqrt{mv}\)
lattice points. Stirling's inequalities give, at each such point \(x\),
\[
 \Pp(P=x)\asymp m^{-1/2}e^{-x^2/2},
\]
since \(x^4/m=O(1)\) and \(|x^2-u^2|=O(1)\).
For \(u\ge1\), successive binomial masses above the first point
above \(u\) have ratio at most \(e^{-2u/\sqrt m}\).
Their geometric sum is at most \(C\phi(u)/u\).
When \(u<1\), the dimension is bounded and the last inequality in
\eqref{low:eq:binomial-band} follows after changing \(C\).
This proves \eqref{low:eq:binomial-band}.

The suffix cannot bring back a prefix already above \(u\), because
it is then inactive. On \(A\), it creates an exit with probability
at least \(\Pp(Z>1/2)\ge c\). The last inequality follows from
\(\E Z^2=1\), \(\E Z^4\le3\), Paley--Zygmund applied to \(Z^2\),
and symmetry. Before normalization, the decrease in the
one-coordinate distribution function is therefore between
\(c\pi\) and \(\pi\).
Taking independent copies across coordinates makes the difference
of the two rectangle probabilities at \(u\), before normalization, at least
\(cd\pi\ge c'h\).
For large \(d\) this follows from the bound \(1-C/d\) for the
success probabilities. For bounded \(d\), symmetry gives a
positive baseline success probability, and \(h_*\) can be reduced
to preserve a fixed positive lower bound.

The corresponding normalized threshold for the adaptive array is
\(u/\sqrt{1+\pi v}\). The difference of the Gaussian product
probabilities at the two thresholds is at most
\[
 C\pi v\log(ed)=O(h^3/d).
\]
Indeed, the endpoint shift is at most \(Cu\pi v\), and the
density at the shifted endpoint is comparable to \(\phi(u)\)
because \(\pi vu^2=O(h^2)\).
For small \(h_*\), the triangle inequality puts an error of at
least \(ch\) on one of the two admissible arrays.

For finite \(p\), the quadratic-variation error satisfies
\[
 \big\|\|Q-I_d\|_{\max}\big\|_p
 \le v\{(d\pi)^{1/p}+\pi\}
 \le Cv h^{1/p},
\]
with the same assertion using \(h^{1/p}=1\) at \(p=\infty\).
Choose
\[
 h=c_0\left[1\wedge
 \{\delta\log(ed)\}^{1/(2+1/p)}\right]
\]
with a sufficiently small absolute \(c_0\). This makes
\(\|\|Q-I_d\|_{\max}\|_p\le\delta\).
If \(h<C_*\sqrt{\log(ed)/n}\), the first term in
\eqref{low:eq:bracket} absorbs it. If \(\log(ed)>n^{1/3}\),
the capped first term is already one. Finally, a scalar
Rademacher sum supplies a positive lower bound for each of the
finitely many remaining \(n\), uniformly in \(d\), by cylinder
padding. Combining the separately attained lower bounds proves
the result.
\end{proof}

\section{Proofs of the main theorems}\label{main:proofs}

We now combine the upper bounds and lower constructions in the order
of the main statements. The convergence theorem uses the separate
relative comparisons, so its proof does not inherit the
additional powers of \(\log n\) in the quantitative bounds.

\begin{proof}[Proof of Theorem~\ref{main:joint}]
For bounded and conditional exponential increments,
Corollary~\ref{quant:cor:joint-rates}, with
\(u=\max\{2,s\}\), gives the upper bounds in
\eqref{main:eq:rates}. In the bounded case the term due to the increments is
\[
 CBn^{-1/2}\{\log(en)\}^{1+3u}\{\log(ed)\}^{3/2}.
\]
Under conditional exponential moments, truncation introduces at most one further power of
\(\log(en)\). For marginal exponential moments, the same corollary
uses the truncation envelope
\(C_K\{\log(2end)\}^{1/\alpha}/\sqrt n\).
Since
\[
 \{\log(2end)\}^{1/\alpha}
 \le C\{\log(en)\}^{1/\alpha}\{\log(ed)\}^{1/\alpha},
\]
the common power \(2+3u\) of \(\log(en)\) gives
\eqref{main:eq:marginal-rates}. The term involving \(\delta\) in
that corollary has exponent
\((1+2u)/(2+1/p)\le2+3u\).
For \(2<q\le3\), substitute
\(v_q\le B^qn^{1-q/2}\) in Proposition~\ref{poly:prop:joint}
to obtain the upper bound in \eqref{main:eq:moment-rates}.

Proposition~\ref{low:cubic} gives the lower bounds for the
bounded and marginal-tail increment terms, with \(Q=I_d\).
The bounded example belongs to both conditional exponential
classes: its unnormalized increments are bounded by two and \(K\ge4\).
Proposition~\ref{poly:prop:identity-lower} gives the finite-moment
lower bound for all sufficiently large \(n\), uniformly in \(d\),
with \(q\)th moment of the largest absolute coordinate at most \(2^q\).

Proposition~\ref{low:bracket} supplies the covariance contribution.
Its unnormalized increments are bounded by two and the normalized sum
has covariance \(I_d\), so the array belongs to every class in the theorem.
The two terms in the lower bound may be attained by different arrays. For
\(x,y\ge0\),
\[
 \max\{1\wedge x,1\wedge y\}
 \ge\tfrac12\{1\wedge(x+y)\}.
\]
Taking the class supremum over these separate examples proves the
joint lower bounds. The light-tail constructions apply for every
\(n\ge2\); the finite-moment construction requires the fixed
threshold stated in part~\textup{(iii)}.
\end{proof}

\begin{proof}[Proof of Theorem~\ref{main:validity}]
Assume the dimension condition for the class under consideration and
\(\delta_n\log(ed_n)\to0\). For finite \(p\) and every
\(\varepsilon>0\), Markov's inequality gives
\[
 \sup_{\text{admissible arrays}}
 \Pp\{\log(ed_n)\|Q-\Sigma\|_{\max}>\varepsilon\}
 \le
 \left\{\frac{\delta_n\log(ed_n)}{\varepsilon}\right\}^{p}
 \longrightarrow0.
\]
For \(p=\infty\), the probability is eventually zero.
Proposition~\ref{thr:prop:light-validity} therefore applies to
every sequence of admissible bounded, conditional exponential,
or marginal exponential arrays. For every fixed \(q>2\),
Proposition~\ref{poly:prop:validity} gives the same conclusion
under the corresponding finite-moment dimension condition.

These sequential conclusions imply convergence of the class
suprema. Otherwise a subsequence would admit arrays and rectangles
whose errors exceed a fixed positive constant. Choosing one such
array at each index would contradict the appropriate convergence
proposition. This proves sufficiency, including arbitrarily slowly
subcritical dimension sequences.

For necessity, Proposition~\ref{low:bracket} requires
\(\delta_n\log(ed_n)\to0\) in every class. Its bounded example
satisfies all the increment assumptions under consideration.
Proposition~\ref{low:cubic} supplies the cubic dimension
obstruction for bounded and conditional exponential increments,
and the exponent \(1/(3+2/\alpha)\) for marginal exponential
increments.

For \(2<q\le3\), Proposition~\ref{poly:prop:identity-lower}
requires
\[
 n^{1-q/2}\{\log(ed_n)\}^{q/2}\to0,
 \qquad\text{equivalently}\qquad
 \log(ed_n)=o(n^{1-2/q}).
\]
In this range \(1-2/q\le1/3\).
For \(q>3\), the bounded example of Proposition~\ref{low:cubic}
belongs to the moment class and requires
\(\log(ed_n)=o(n^{1/3})\); now \(1-2/q>1/3\).
Every example used for the dimension condition has \(Q=I_d\), and hence is admissible
for every \(\delta\ge0\). These arguments give
all three equivalences.
\end{proof}

\begin{proof}[Proof of Theorem~\ref{main:identity}]
Corollary~\ref{quant:cor:identity-light-tails} gives the bounded
and marginal upper bounds when
\(\log(ed)\ge\{\log(en)\}^2\).
The corresponding constructions in Proposition~\ref{low:cubic}
satisfy \(Q=I_d\) almost surely and give the lower bounds. Their
envelopes fit every fixed \(B\ge2\) or \(K\ge4\), respectively.

For \(2<q\le3\), combine
Propositions~\ref{poly:prop:identity-upper}
and~\ref{poly:prop:identity-lower}, inserting
\(v_q\le B^qn^{1-q/2}\) in the upper bound.
The latter construction has \(q\)th moment of the largest absolute coordinate at most
\(2^q\) and works in every dimension for all sufficiently large
\(n\). The required lower bound on \(n\) and the comparison constants depend only on the
fixed moment parameters. This proves
\eqref{main:eq:identity-polynomial} throughout its stated range.
\end{proof}

\appendix
\section{Smaller powers of \texorpdfstring{\(\log n\)}{log n}}\label{app:logs}

The powers of \(\log(en)\) in the main results can be reduced
in several regimes. These refinements follow from the Gaussian
cylinder estimate and the separate optimization of the finite-moment
bound. They leave the powers of \(n\), \(\log d\), and \(\delta\) unchanged.

\begin{proposition}\label{main:rem:logs}
The power \(2+3\max\{2,s\}\) in the bounded and
exponential-moment parts of Theorem~\ref{main:joint}
can be replaced by \(5\) if
\(\|\Sigma\|_{\mathrm{op}}\le M_0\) for a fixed \(M_0\),
or if
\[
 \log(ed)\ge\{\log(en)\}^{4s+1}.
\]
For \(2<q\le3\), the power in \eqref{main:eq:moment-rates}
can always be replaced by \(1+\max\{2,s\}\), and by \(2\)
in either of these two regimes. Constants under the operator-norm
assumption may also depend on \(M_0\).
The correlation condition in \eqref{main:eq:setting} implies the operator-norm assumption when
\(0<s\le1\).
\end{proposition}

\begin{lemma}\label{app:lem:cylinder}
Under the hypotheses of Proposition~\ref{thr:prop:quantitative-cylinder},
its conclusions hold with \(u=1\) if
\(\|\Sigma\|_{\mathrm{op}}\le M_0\).
For every \(s>0\), the same reduction also holds with
\begin{equation}\label{thr:eq:high-dimensional-weight}
 W\le C\log(en)
 \left\{\sqrt{\log(2end)}+\{\log(en)\}^{2s+1/2}\right\},
\end{equation}
for the retained rectangle. When a proper cylinder is needed,
its Gaussian probability is at most \(Cn^{-8}\). Consequently
\eqref{thr:eq:quantitative-weight} holds with \(u=1\) when
\(\log(ed)\ge\{\log(en)\}^{4s+1}\).
\end{lemma}

\begin{proof}
We prove the two cylinder alternatives separately.

\paragraph{Bounded operator norm.}
The Gaussian H\"older inequality
\citep[Theorem~1]{CDP2015} gives
\[
 \Pp(G\in A)
 \le\prod_j(1-q_j)^{1/M_0}
 \le\exp\{-V(A)/M_0\},
 \qquad q_j=\bar\Phi(b_j).
\]
To check its matrix condition, write \(G_j=\langle v_j,Z\rangle\)
on the Gaussian support. The nonzero eigenvalues of
\(\sum_j v_jv_j^\top\) and \(\Sigma\) agree, so
\(M_0^{-1}\sum_jv_jv_j^\top\preceq I\).
If the sum \(V(A)\) exceeds \(10M_0\log(en)\), retain the selected coordinate constraints
until their tail-probability sum reaches this value, with overshoot at most one.
If \(s\le1\), then
\(\sum_j|\Sigma_{ij}|\le\sum_j|\Sigma_{ij}|^s\le M\);
the operator norm is at most \(M\).

\paragraph{The correlation summability assumption.}
For coordinates with \(b_j>A_0\sqrt{\log(en)}\), the proof of
Proposition~\ref{thr:prop:quantitative-cylinder} either bounds the sum
of their tail probabilities by \(C\log(en)\), or selects a subset of
constraints with this bound whose cylinder has Gaussian probability
at most \(Cn^{-8}\). For the remaining thresholds, we use a coloring
argument.

For \(b_j\le A_0\sqrt{\log(en)}\), join coordinates whose absolute
correlation exceeds \(\theta_0=\{8\log(en)^2\}^{-1}\).
The correlation summability bound limits the degree of this graph,
so it can be colored with at most \(C_{s,M}\{\log(en)\}^{2s}\) colors.
On each color, comparison with equicorrelation \(\theta_0\)
and the common-factor argument in that proposition give
\[
 \Pp(\text{no exit in that color})
 \le Cn^{-10}+\exp\{-c_{A_0}\textstyle\sum q_j\}.
\]
Here the common factor changes a threshold by at most
\(C_{A_0}/\sqrt{\log(en)}\); negative thresholds only improve
the scaling term. If the sum of exit probabilities for the low thresholds exceeds
\(C\{\log(en)\}^{2s+1}\), the rectangle is contained in a cylinder
of probability at most \(Cn^{-8}\). Retain all constraints in a group if its tail-probability sum is below
the stated bound, and otherwise use the selected cylinder. Combining
the bounds for the high and low thresholds gives
\eqref{thr:eq:high-dimensional-weight}: the smaller signed threshold
\(b_j\) controls both endpoint weights of its coordinate by
\(C\{1+(b_j)_+\}\bar\Phi(b_j)\).

Finally,
\(\log(ed)\ge\{\log(en)\}^{4s+1}\) implies
\(\log(2end)\le C\log(ed)\); substitution into
\eqref{thr:eq:high-dimensional-weight} gives the stated \(u=1\)
bound.
\end{proof}

\begin{proof}[Proof of Proposition~\ref{main:rem:logs}]
In the proofs of Proposition~\ref{quant:prop:pathwise} and
Corollary~\ref{quant:cor:joint-rates}, the dependence of the
powers of \(\log(en)\) on \(u\) enters through the cylinder weight
\(W\le C\{\log(en)\}^{u}\sqrt{\log(ed)}\).
Lemma~\ref{app:lem:cylinder} therefore permits \(u=1\)
in \eqref{quant:eq:pathwise} and \eqref{quant:eq:joint-rates}
in either stated regime. The additional factor for the light-tail increment term is then
at most \(\{\log(en)\}^{2+3u}=\{\log(en)\}^5\),
and the term involving \(\delta\) has a smaller power.

For finite moments, the same cylinder substitution applies in
\eqref{poly:eq:joint-before-optimization}. Use its bandwidth
\eqref{poly:eq:joint-bandwidth}, with
\(\gamma=(2+1/p)^{-1}\), \(r=\lceil100\log(en)\rceil\)
and \(T=D\sqrt{\log(2nd)}\).
\(T^2\le C\log(en)\log(ed)\), and
\[
 \begin{split}
 \{\log(en)\}^u(rv_qT^q)^{1/(q+1)}
 &\le C\{\log(en)\}^{u+(q+2)/(2(q+1))}
          (v_q\{\log(ed)\}^{q/2})^{1/(q+1)},\\
 \{\log(en)\}^u(r\delta T^2)^\gamma
 &\le C\{\log(en)\}^{u+2\gamma}
                         (\delta\log(ed))^\gamma .
 \end{split}
 \]
Both exponents of \(\log(en)\) are at most \(u+1\).
The term \(n^{-8}\) is absorbed by the increment-dependent term, since
\(2<q\le3\). This gives the power \(1+\max\{2,s\}\) in general and
the power \(2\) when \(u=1\).
\end{proof}

These powers of \(\log n\) give valid upper bounds; we do not assert
that they are optimal. Proposition~\ref{low:scalar}
shows that the scalar conditional sub-exponential class can have error
at least \(c\{\log(en)\}^2/\sqrt n\), so a uniform multiplicative bound
relative to the term \(n^{-1/2}\) over that full scalar class requires at least
two powers of the sample logarithm.

\section{A scalar logarithmic obstruction}\label{low:scalar-section}

\begin{proposition}\label{low:scalar}
For every sufficiently large even \(n\), there are scalar martingale
differences with
\[
 \E_{i-1}e^{|X_i|/4}\le2,\qquad
 \frac1n\sum_i\E_{i-1}X_i^2=1
\]
such that
\[
 \Pp\left(n^{-1/2}\sum_iX_i\le0\right)-\frac12
 \ge c\frac{(\log n)^2}{\sqrt n}.
\]
Their initial sigma-field is trivial and their conditional variances
are deterministic.
\end{proposition}

\begin{proof}
Let
\[
 m=n/2,\quad v=n^{-1/2},\quad V=mv,\quad
 w=\frac{n-V}{n-m},\quad A=\frac12\log n.
\]
Use \(n-m\) independent \(N(0,w)\) prefix increments, with sum \(P\).
The suffix has \(m\) increments, each of conditional variance \(v\).
Its default law is \(N(0,v)\). At suffix step \(j\), let \(H_{j-1}\)
be the preceding suffix sum and \(r=(m-j)v\).
Replace the default law by
\[
 \Pp(Y=A)=\frac{v}{A^2+v},\qquad
 \Pp(Y=-v/A)=\frac{A^2}{A^2+v}
\]
exactly on the predictable event
\[
 |P|\le\sqrt V,\qquad
 V^{3/4}\le r\le V^{4/5},\qquad
 |P+H_{j-1}|\le\sqrt r/2.
\]
Fresh independent random variables realize each conditional law.
Direct calculation gives
\[
 \E Y=0,\quad \E Y^2=v,\quad
 \E Y^3=v(A-v/A),\quad
 \E Y^4=v(A^2-v+v^2/A^2).
\]
The total conditional variance is exactly \(n\).
For the active law,
\[
 \E e^{|Y|/4}
 \le e^{v/(4A)}+\frac{v}{A^2}e^{A/4}=1+o(1).
\]
For a centered Gaussian of variance \(w'\le2\), the inequality
\(|x|/4\le1/8+x^2/8\) gives
\[
 \E e^{|N(0,w')|/4}
 \le e^{1/8}(1-w'/4)^{-1/2}<2.
\]
This verifies the conditional exponential assumption.

Condition on a prefix sum \(P\) with \(|P|\le\sqrt V\), and write
\(t=-P\). For Gaussian completion after the current
suffix step, use
\[
 f_r(x)=\Phi((t-x)/\sqrt r).
\]
At zero remaining variance set \(f_0(x)=\ind_{\{x\le t\}}\).
The last suffix step always has the default Gaussian law, so it
contributes no replacement error and the terminal sum has no atoms.
On the switching strip,
\(f_r'''(H_{j-1})\ge cr^{-3/2}\), while
\(\|f_r^{(4)}\|_\infty\le Cr^{-2}\).
Taylor expansion against the default \(N(0,v)\) step therefore
gives a conditional CDF gain at least
\[
 cvA r^{-3/2}-CvA^2r^{-2}\ge c'vA r^{-3/2}.
\]
The last inequality uses \(r\ge V^{3/4}\gg A^2\).
At default steps the comparison is exact.

We need a lower bound on the switching-strip probability.
Add an independent \(N(0,\sqrt V)\) variable to a partial suffix
sum \(H_\ell\). Its Gaussian replacement telescope, expanded to
second order, has error at most
\[
 CA\sum_{a=0}^{\ell-1}
 \frac{v}{(\sqrt V+av)^{3/2}}
 \le CA V^{-1/4}.
\]
Each possible suffix law has third absolute moment at most \(CAv\),
and the conditional variances are exactly \(v\), which justifies
this deterministic denominator sum.
At a switching time, the comparison Gaussian has variance
\(\ell v+\sqrt V\asymp V\). Its mass on
\([t-\sqrt r/4,t+\sqrt r/4]\) is at least \(c\sqrt{r/V}\),
uniformly in \(|t|\le\sqrt V\).
The error \(CA V^{-1/4}\) is negligible relative to this mass,
which is at least \(cV^{-1/8}\).
Removing the added Gaussian and widening to radius \(\sqrt r/2\)
costs at most \(2e^{-cV^{1/4}}\). Thus
\[
 \Pp(|H_{j-1}-t|\le\sqrt r/2\mid P)\ge c\sqrt{r/V}.
\]
Telescoping the one-step Gaussian replacements gives a conditional
CDF gain at least
\[
 \frac{cA}{\sqrt V}
 \sum_{V^{3/4}\le(m-j)v\le V^{4/5}}\frac{v}{(m-j)v}
 \ge \frac{cA\log V}{\sqrt V}.
\]
Off the activation event the law is entirely Gaussian.
The prefix \(P\sim N(0,n-V)\) has density at least \(c/\sqrt n\)
on \([-\sqrt V,\sqrt V]\).
Integrating the conditional gain over this interval proves
\[
 \Pp\left(\sum_iX_i\le0\right)-\frac12
 \ge c\,\frac{A\log V}{\sqrt n}
 \ge c'\frac{(\log n)^2}{\sqrt n}.
\]
\end{proof}

The example distinguishes an absolute conditional exponential
moment from a Bernstein bound relative to conditional variance.
For \(\xi_i=X_i/\sqrt n\), it satisfies the printed assumptions
\textup{(A1)}--\textup{(A2)} of \citet{fan2014generalization}
with \(c_0=1/4\), \(c_1=2\), and zero error in the predictable quadratic variation.
At zero tilt, their Lemma~3.4 asserts an error of order
\(\log n/\sqrt n\). Proposition~\ref{low:scalar} contradicts
that quantitative conclusion under those hypotheses. The
conclusions of their Theorem~2.1 at threshold zero give the same
incompatible rate, since the constructed terminal law has no atoms.
Its conditional variances are deterministic and bounded above by
\(2/n\) in normalized units, so such an upper bound does not
remove the conflict.

The variance-relative Bernstein hypothesis used by
\citet{fan2013cramer} is different. In the active suffix here,
\(\E_{i-1}X_i^3/\E_{i-1}X_i^2=A-v/A\to\infty\).
Thus the absolute exponential bound does not imply that hypothesis
at scale \(n^{-1/2}\). The contradiction concerns this quantitative rate. It does not rule out
qualitative moderate-deviation statements under an absolute conditional
exponential-moment bound.

\section*{AI Statement}
GPT-6 Pro was used to assist with the development of mathematical arguments and the preparation of this manuscript. The authors take full responsibility for the content of the paper.

\bibliographystyle{apalike}
\bibliography{Optimal_High_Dimensional_Martingale_CLTs}
\end{document}